\documentclass[leqno,12pt]{article}

\usepackage [dvips]{graphics,psfrag}
\usepackage{amssymb,amsmath}

\newtheorem{theorem}{Theorem}
\newtheorem{proposition}[theorem]{Proposition}
\newtheorem{lemma}[theorem]{Lemma}
\newtheorem{definition}[theorem]{Definition}
\newtheorem{corollary}[theorem]{Corollary}
\newenvironment{proof}{\noindent{\bf Proof. }}{\par}
\newcommand{\C}{{\mathbb C}}
\newcommand{\R}{{\mathbb R}}
\newcommand{\Z}{{\mathbb Z}}

\newcommand{\CB}{{\cal B}}

\newcommand{\CH}{{\cal H}}

\newcommand{\Ker}{\operatorname{Ker}}
\newcommand{\vol}{\operatorname{vol}}
\newcommand{\n}{\mathfrak{n}}
\renewcommand{\c}{\mathfrak{c}}
\renewcommand{\a}{\mathfrak{a}}
\title{Residues formulae for volumes and Ehrhart polynomials  of
convex polytopes.} 
\author{Welleda Baldoni-Silva  and Mich{\`e}le Vergne}
\date{January 2001}
\begin{document}
\maketitle

\section*{Introduction}

These notes on volumes and Ehrhart polynomials of convex polytopes 
are afterwards thoughts  on lectures already delivered in Roma, by 
the second author  in December 1999. The subject of these lectures 
was to explain  Jeffrey-Kirwan residue formula for volumes of 
convex polytopes [J-K]  and  the residue formula for Ehrhart 
polynomials  of rational polytopes [B-V 1, 2]. The  main concept 
used in these formulae  is the study of rational functions with 
poles on 
 an arrangement of hyperplanes and  the total residue 
 of such a function, which can be computed recursively.
  Now, what  about 
``concrete" polytopes ? The residue formula  leads in principle to 
an answer, but to be carried out in practice, without further 
thinking, requires too many steps.  Thus we  have been fascinated 
by  calculations of Chan-Robbins in [CR], and  the conjectures of 
Chan-Robbins-Yuen [C-R-Y] on the volume of the a certain polytope 
(proved by Zeilberger [Z]) and of some of its faces, as well as  
by other examples  given  by Pitman-Stanley [Pi-S], Stanley [S].  
All these examples are related with the root system $A_n$.  
 Therefore we will  study here in more details what can be 
said on residues formulae for this root system and subsets of this
system.

Let us sketch briefly the content of these notes. Short 
bibliographical comments are also given at the end of each 
section. We first recall the definition of total residue in 
Section \ref{total} and  Jeffrey-Kirwan residue formula for 
volumes  of general  convex polytopes (and give a proof) in 
Section \ref{jeffrey} and   Appendix (Section \ref{appendix1}). A  
formula  for change of variables in total residues is obtained in 
Section \ref{changeof}. 
  We go on considering flow polytopes (see Examples below, and Definition in 
Section \ref{flow}), and defining a special class of flow 
polytopes: the cascade polytopes . In Section \ref{chambers},  we 
show that  the Jeffrey-Kirwan residue formula for the volume of   
cascade polytopes is  an iterated constant term formula. We then 
give  a  proof of a divisibility property  between the volume of a 
cascade polytope and the volume of some particular face (this 
implies Conjecture 4 of [C-R-Y]) and a similar divisibility 
property for Ehrhart polynomial.  
 We   carry out the computation of the residue for the complete
 flow polytope, following the indications of Zeilberger  [Z] .
  It thus leads  to the known formula [Z] for the 
volume of the complete flow polytope. Our method  is algebraic 
though combinatorists would probably prefer another type of proof. 
For example, Chan-Robbins-Yuen  proved that their  polytope $P$ 
could be explicitly decomposed in a union of $K$ elementary 
simplices (each of volume equal $\frac{1}{(\dim P)!}$). Later $K$ 
was computed by Zeilberger. Here we apply directly the  residue  
formula for the volume,  and  we do not use any simplicial 
decomposition.

In Section \ref{eh}, we state  a residue formula   for Ehrhart 
polynomials 
 of flow polytopes (the proof is not given here, but is similar to the one given for volumes;
 it is a reformulation  of Khovanskii-Pukhlikov theorem [KP]). 

In Section \ref{nice}, as suggested by R. Stanley, we show   that 
formulae for  volumes  of  a family of  flow polytopes   can be  
transmuted  very explicitly  to formulae  for Ehrhart polynomials, 
 due to the formula of change of variables in residues.  
Sections \ref{chambers}, \ref{volume}, \ref{calculation} and 
\ref{CR}    can be skipped if the reader is interested  only in 
the nice formula of Section \ref{nice}. Finally, we state a 
symmetry property that must satisfy the Kostant partition 
function.

Let us now state here more precisely the  setting and the contents 
of these notes.

 Let $V$ be a $r$-dimensional real vector space and
$V^*$ its dual vector space. Let $\Phi=\{\alpha^1, 
\alpha^2,...,\alpha^N\}$ be a  sequence  of non-zero linear forms 
on $V^*$ all lying in an open half space (we do not assume the 
$\alpha^k$ to be distinct) and spanning $V^*$. We denote by 
$\Delta^+$ the set $\{\Phi\}$ (we mean $\Delta^+ $ and $\Phi$ are 
the same sets, but $\Phi$ may have multiplicities), and by 
$\Delta$ the set $\Delta^+\cup -\Delta^+$ . The closed convex cone 
$C(\Delta^+)$  generated by $\Delta^+$ is  decomposed  as a  union 
of  closure  of big chambers $\c$ as in [B-V 1].  
 For $a\in V^*$, we denote by $P_\Phi(a)\subset \R_+^N$ 
 the convex polytope consisting  of all solutions $(u_1,u_2, ..., 
u_N)$, in non negative real numbers $u_k$, of the equation 
$$\sum_{k=1}^N u_k \alpha^k=a$$ and by $\vol P_\Phi(a)$ its 
volume.  Any convex polytope can be realized in that way.

Consider the space $R_{\Delta}$ of rational functions of $x\in 
V_\C$ with poles on the union $\CH_{\C}$ of the hyperplanes 
$\alpha^k(x)=0$. A subset $\sigma$ of $\Delta$ is called a {\bf 
basis of} $\Delta$, if the elements $\alpha\in \sigma$ form a 
basis of $V^*$ . For such $\sigma$, set $$ 
f_{\sigma}(x):=\frac{1}{\prod_{\alpha\in \sigma}\alpha(x)}. $$ In 
appropriate coordinates $x_1,...,x_r$, the function $f_\sigma$ is 
simply $\frac{1}{x_1 x_2... x_r}$ and we denote  by $S_{\Delta}$ 
the subspace  of $R_{\Delta}$ spanned by  such  ``simple" elements 
$f_{\sigma}~$. The vector space $S_{\Delta}$ is contained in the  
homogeneous component of degree $-r$  of $R_{\Delta}$ and we have 
the direct sum decomposition $$R_\Delta=S_\Delta\oplus 
(\sum_{k=1}^r \partial_k R_\Delta).$$  

We call  the projection map $$ Tres_{\Delta}: R_{\Delta}\to 
S_{\Delta} $$  according to this decomposition the total residue 
map.  This projection vanishes outside the homogeneous component 
of degree $-r$ of $R_\Delta$. The total residue of a function is 
again a function.  It consists of computing  a rational function 
$f(x)$  up to derivatives, in other words,  $f(x) dx $ and 
$(Tres_\Delta (f))(x) dx$ are top dimensional  holomorphic forms 
representing  the same  cohomology class  on $V_\C-\CH_{\C}$. 

 We consider $a\in V^*$. Then
 $$J_\Phi(a)(x)=Tres_\Delta(\frac{e^{\langle a,x\rangle}}{\prod_{k=1}^N\alpha^k(x)})=
 \frac{1}{(N-r)!} Tres_\Delta (
\frac{\langle  a,x\rangle^{N-r}}{\prod_{k=1}^N\alpha^k(x)})$$ is 
an element of the vector space $S_\Delta$.

We denote by $$K_\Phi(a)= Tres_\Delta( \frac{e^{\langle  
a,x\rangle}}{\prod_{k=1}^N(1-e^{\langle  \alpha^k,x\rangle})})$$
the  ``periodic" version  of $J_{\Phi}$. This is also an element  
of the vector space $S_\Delta$.

Now, to  each big chamber $\c$ of the subdivision of $C(\Delta^+)$ 
is associated a linear form $f\to \langle  \langle  
\c,f\rangle\rangle$ on $S_\Delta$. It takes value $1$ or $0$ on a 
normalized multiple of $f_\sigma$ whether or not $\c$ is contained 
in $C(\sigma)$. This is explained in Section \ref{jeffrey} and in 
the Appendix (Section \ref{appendix1}).  Here are the two 
fundamental formulae we will use in these notes.  

{\bf Formula 1:} for $ a\in \overline{\c}$,  we have

$$\vol P_\Phi(a)=\langle  \langle  \c, J_\Phi(a)\rangle\rangle.$$

  If $\Phi$ spans a lattice in 
  $V^*$ and if $a$ belongs to this lattice, we denote by $k_\Phi(a)$ 
 the number of all solutions $(u_1,u_2, ..., 
u_N)$, in non negative {\bf integral}  numbers $u_k$ of the 
equation $$\sum_{k=1}^N u_k \alpha^k=a.$$ Similarly (under  
simplifying assumptions that we do not state in the introduction), 
the following formula holds.

{\bf Formula 2:} for $ a\in \overline{\c}$, we have

  $$k_\Phi(a)= \langle\langle \c, K_\Phi(a)\rangle\rangle.$$

  In particular, as well known, the functions $a\mapsto \vol P_\Phi(a)$
  and the Ehrhart function   $ a\mapsto k_\Phi(a)$
  are polynomials  on a big chamber $\c$. 
(The notion of Minkowski sums and of big chambers are intimately 
related:   if ${\bf a_1},..., {\bf a_m}$ are vectors in the 
closure of a big chamber $\c$, then the polytope $P_\Phi({\bf 
a_1}+\cdots+{\bf a_m})$ is isomorphic to the  Minkowski sum of the 
polytopes  $P_\Phi({\bf a_k})$.)

This shows that calculation of volumes  (or of Ehrhart 
polynomials) is done by an algebraic version of ``integration": it 
consists of computing  particular rational function $f(x)$ on 
$V_\C$  modulo derivatives. Furthermore this formula shows clearly 
that the  $S_\Delta$- valued  {\bf polynomial function} $a\to 
J_\Phi(a)$  determines entirely the {\bf locally polynomial 
function} $a\mapsto   \vol P_\Phi(a).$

 Thus the calculation of  the volume of a flow polytope  $P_\Phi(a), a \in \overline{\c}$, 
  is  divided in TWO  problems. 
 
 A) Compute the linear form $f\to \langle\langle \c,f\rangle\rangle$.
 
 B) Compute  the function $$\frac{e^{\langle a,x\rangle}}{\prod_{k=1}^N\alpha^k(x)}$$ up to derivatives.

  \bigskip

  We will now study some particular cases, where A) or B) can be solved.

The Chan-Robbins-Yuen  polytope $CRY_r$ consists of solutions  
$(x_{ij})\geq 0$  where $1\leq i< j\leq (r+1)$ of the linear 
equations: 

\begin{itemize}
\item

 $$\sum_{i=2}^{r+1}x_{1,i}=1, $$

\item

 $$\sum_{i=1}^{r}x_{i, r+1}=1,$$
 
 \item
 For $2\leq j\leq r$,
 $$\sum_{i=1}^{j-1}x_{ij}=\sum_{k=j+1}^{r+1}x_{jk}.$$
 
\end{itemize}

\newpage
The CRY polytope is  also called the complete flow polytope. 
Indeed, consider a graph $G$ with $r+1$ vertices  $1,2,3,..., 
r,r+1$ 
 and  edges $i\mapsto j$ ($1\leq i< j\leq r+1)$, for all $i< j$.

\begin{figure}
        \centering
        \psfrag{b}{$\text{flow in}$}
        \psfrag{a}{$\text {flow out}$}
        \psfrag{1}{$1$}
        \psfrag{2}{$2$}
        \psfrag{3}{$3$}
        \psfrag{4}{$4$}
\includegraphics{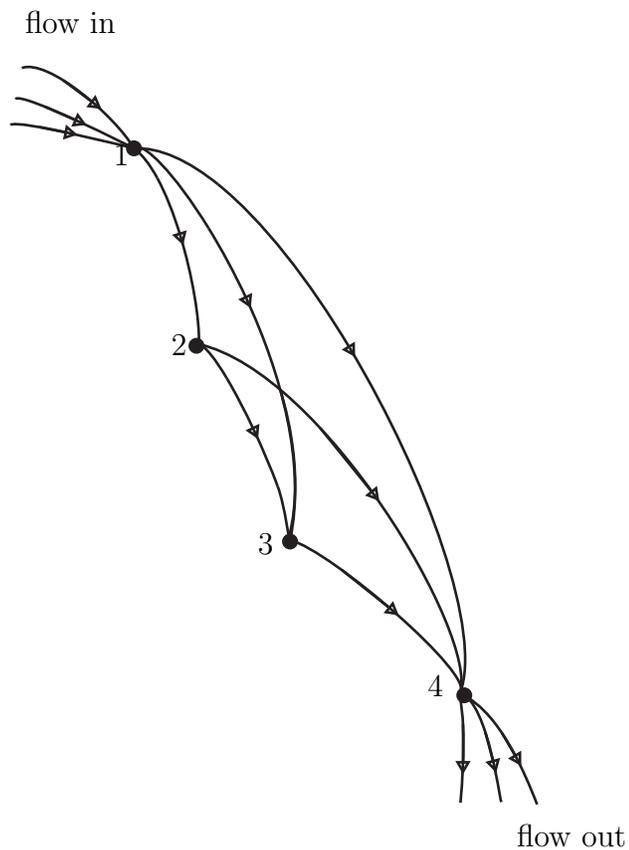}
        \caption{Graph for $CRY_3.$}
        \label{fig:1}
\end{figure}

 We can imagine that the  positive quantity $x_{ij}$ is the quantity of liquid at time $t$
 in the branch $i\mapsto j$ of this cascade. Thus the linear equations above
 reflect the constant flow
 of the  cascade  ( See Figure 1).

 We  can also consider other graphs with possibly multiple edges.
  The associated polytope  will be called a flow polytope. See the precise definition in Section
  \ref{flow}.
 For example, Pitman-Stanley  polytope is the polytope 
 consisting of solutions in non negative numbers
 of  the inequations
  $$y_{1,r+1}\leq a_1,$$ 
$$y_{1,r+1}+y_{2,r+1}\leq a_1+a_2,$$ 
$$y_{1,r+1}+y_{2,r+1}+y_{3,r+1}\leq a_1+a_2+a_3,$$ $$...$$ where 
$a_i$ are non negative real numbers.

It is   associated to a flow graph, with $r+1$ vertices, and  
edges from $i$ to $r+1$ and from $i$ to $i+1$ and last edge 
$r\mapsto r+1$ 
 of multiplicity $2$.

\begin{figure}
        \centering
        \psfrag{a1}{$a_1$}
        \psfrag{a2}{$a_2$}
        \psfrag{a3}{$a_3$}
        \psfrag{l}{$z_{12}$}
        \psfrag{m}{$z_{23}$}
        \psfrag{n}{$z_{34}$}
        \psfrag{y14}{$y_{14}$}
        \psfrag{y24}{$y_{24}$}
        \psfrag{y34}{$y_{34}$}
        \psfrag{a}{$y_{14}\leq 0$}
        \psfrag{b}{$y_{14}+y_{24}\leq a_1+a_2$}
        \psfrag{c}{$y_{14}+y_{24}+y_{34}\leq a_1+a_2+a_3$}
        
        \includegraphics{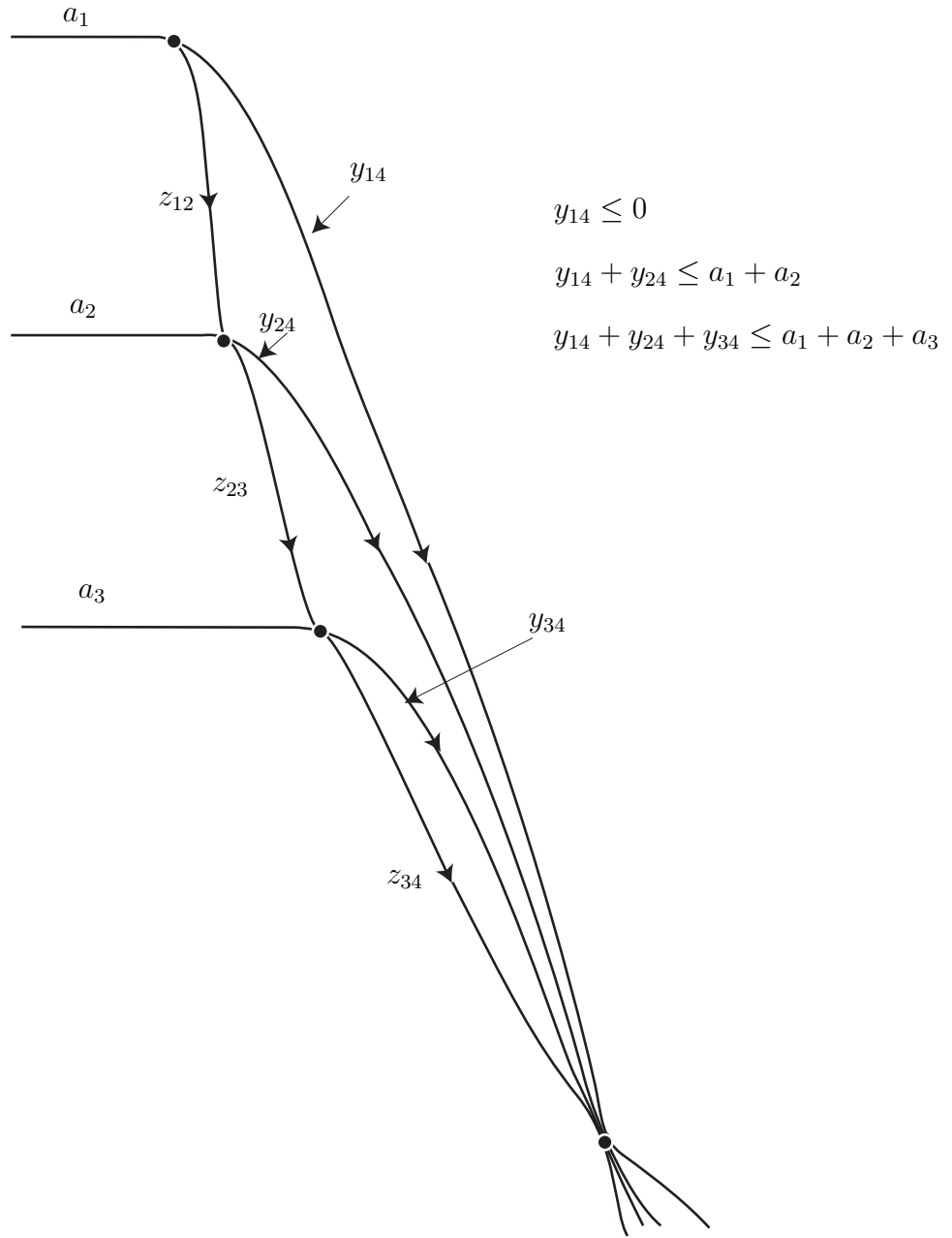}
        \caption{Graph for $\pi_3(a).$ }
        \label{fig:2}

\end{figure}

We can imagine that $a_1$, $a_2$,...$a_r$ are additional sources 
of water, so that the quantities of water $y_{i,r+1}$ in the 
streams $i\to r+1$ satisfies the inequalities stated above, with a 
leakage of  $z_{i,i+1}$. ( See g 2).

 In the rest of these notes, we consider the positive root 
system  $$A_r^+=\{(e^i-e^j), 1\leq i< j\leq (r+1)\}.$$ The system 
$A_r^+$ spans the vector space $$E_r=\{a=a_1 e^1+ a_2 
e^2+...+a_{r+1}e^{r+1}, \sum_{i=1}^{r+1}a_i=0\}.$$ The   cone  
$C(A_r^+)\subset E_r$  generated by positive roots is the cone 
$a_1\geq 0$, $a_1+a_2\geq 0$,..., $ a_1+a_2+\cdots +a_r\geq 0, 
a_1+a_2+...+a_{r+1}=0$. We denote by 
 $\c^+$ the  open set  of $C(\Delta^+)$ defined by
 $$\c^+=\{a
\in C({A_r}^+)\,\,\text{such that} \,\,a_i>0 ,\,i=1,..,r\}.$$ It 
is a big chamber of our subdivision, and  will be called the nice 
chamber.  

Thus the functions $a\mapsto \vol P_\Phi(a)$ 
  and the Ehrhart function   $ a\mapsto k_\Phi(a)$
  on the chamber $\c^+$ are the mixed volumes and the mixed lattice point enumerators
  of the fountain polytopes $P_k=P_{\Phi}(e^k-e^{r+1})$.

Let  $\Phi$  be a sequence of $N$ elements of $A_r^+$ spanning 
$E_r$. By definition, a flow polytope is a polytope isomorphic to 
a polytope $P_\Phi(a)$, where $a$ is an element of $C(A_r^+)$ and 
a cascade polytope the special class of flow polytopes   
$P_\Phi(a)$, where $a$ is  constrained to be in the nice chamber 
of $C(A_r^+)$. 
 The flow graph associated to $\Phi$ has 
$m_{ij}$ edges from $i$ to $j$, if $(e^i-e^j)$ has multiplicity 
$m_{ij}\geq 0$ in $\Phi$. In fact, any polytope associated to a 
smooth toric variety (i.e. Delzant polytopes and limits) can be 
realized as a flow polytope (Szenes: personal communication).

The set of simple roots $$\Pi=\{(e^1-e^2),(e^2-e^3),..., 
(e^{r-1}-e^{r}),(e^r-e^{r+1})\}$$ is a basis of $A_r^+$ and the 
element  $$f_{\Pi}=\frac{1}{(e^1-e^2)(e^2-e^3)\cdots 
(e^{r-1}-e^{r})(e^r-e^{r+1})}$$ is an element of $S_{A_r}$. The 
dimension of $S_{A_r}$ is $r!$ and a particularly nice basis of 
$S_{A_r}$ is given by the elements $w\cdot f_{\Pi}$ where  $w\in 
\Sigma_r$= permutation group of $r$ elements, acts by permutation 
on the set $\{1,2,...r\}$. By THEOREM \ref{tt}, the linear form 
$f\mapsto \langle \langle \c^+,f\rangle\rangle$ is the ``simplest" 
possible: the iterated residue $Ires_{x=0}$. This is defined as 
follows: a function in $R_{A_r}$ is identified with  rational 
function $f(x_1,...,x_{r})$ (setting $e^{r+1}=0$) with poles on 
the hyperplanes $x_i=x_j$ , or $x_i=0$, and we have  

$$\langle\langle \c^+,f\rangle\rangle= Ires_{x=0} f=Res_{x_1=0} 
Res_{x_2=0}...Res_{x_r=0}f(x_1,x_2,..., x_r).$$

Thus the volume of a cascade polytope is given by an iterated 
residue formula. But problem B still remains to be solved.

 When   $a$ is a root of the system $A_{r+1}$, the calculation  
of the total residue of  the particular rational function  
$$J_\Phi(a)(x)=\frac{\langle a,x\rangle 
^{r(r+1)/2}}{\prod_{\alpha\in A_{r+1}^+}\langle \alpha,x\rangle}$$ 
(as well as   very similar examples) will be done in Section 
\ref{calculation}. Let  $a=e^0-e^{r+1}$. Setting as before 
$e^{r+1}=0$, 

$$ J_\Phi(a)(x_0,..., x_r) 
  =\frac{x_0^{\frac{r(r+1)}{2}}}{x_0x_1.... x_r 
\prod_{0\leq i< j\leq r}(x_i-x_j)}.$$

The ``integration" of this function is very familiar in the 
context of Selberg integral ([Se]). However the  usual factor 
$\prod_{0\leq i< j\leq r}(x_i-x_j)$ here is in the denominator and  
``Integration " means that  we explicitly compute $$ 
\frac{1}{x_1.... x_r \prod_{0\leq i< j\leq r}(x_i-x_j)}$$  
 modulo partial derivatives   $\partial_{x_1}, 
\partial_{x_2},..., \partial_{x_r}$, following Aomoto [A]. 
In THEOREM \ref{morris}, we obtain  the formula

$$Tres_{A_{r+1}^+}(\frac{(e^0-e^{r+1})^{r(r+1)/2}}{\prod_{0\leq i< 
j\leq r+1} (e^i-e^j)})$$ $$=\prod_{i=1}^{r-1} \frac{(2i)!}{i! 
(i+1)!} \sum_{w\in \Sigma_r}\epsilon(w)w\cdot 
\frac{1}{(e^0-e^1)(e^1-e^2)\cdots (e^{r-1}-e^{r})(e^r-e^{r+1})}$$ 
where $\Sigma_r$ acts by permutations on $\{1,2,...,r\},$ . 
Similar formulae are also obtained. As pointed out by Zeilberger 
[Z], these calculations are mere reformulations of  Morris 
Identity [M]. The total residue formula replaces here the iterated 
constant term. It would be interesting  to  generalize them 
 to other root systems.

The Chan-Robbins polytope is described as a face of the polytope 
of doubly stochastic matrices in Section \ref{CR}, as 
Chan-Robbins-Yuen originally described it. We apply then THEOREM 
\ref{ct}, to prove a formula for the volume of $CRY_n$ , and  
PROPOSITION \ref{division}  for the volume  of a  particular face 
of it.

In Section \ref{changeof}, we prove a formula of Change of 
variables in total residue. 

As a corollary, the $S_\Delta$- valued polynomial function 
$K_\Phi(a)$ for $\Phi $ a  sequence of $N$ elements of $A_r^+$ is 
deduced from the function $J_\Phi(a)$,  in a way which will  be 
explained now.  Define $$t_j^{\Phi}=\sum_{k=j+1}^{r+1}m_{jk}-1,$$ 
$$s_j^{\Phi}=1-\sum_{k=1}^{j-1}m_{kj}.$$

 We write  the vector valued polynomial $J_\Phi$, which is homogeneous of degree $(N-r),$
   $$J_\Phi(a_1e^1+a_2 e^2+...+a_r e^r-(a_1+a_2+...+a_{r})e^{r+1})$$ 
   $$=\sum _{i_1+i_2+\cdots +i_r=N-r} 
 \frac{a_1^{i_1}}{i_1!}\frac{a_2^{i_2}}{i_2!}\cdots 
 \frac{a_{r-1}^{i_{r-1}}}{i_{r-1}!} \frac{a_r^{i_{r}}}{i_{r}!}f_\Phi( {\bf i}) $$
where $f_\Phi( {\bf i}) $ are  elements of $S_\Delta$.

 Then  we obtain THEOREM  \ref{vw}: 
 $$K_\Phi(a_1e^1+a_2 e^2+...+a_r e^r-(a_1+a_2+...+a_{r})e^{r+1})$$ 
 $$=\sum_{{|\bf i|}=N-r} 
 \binom{a_1+t_1^{\Phi}}{i_1} 
\binom{a_{2}+t_2^{\Phi}}{i_{2}}\cdots  
\binom{a_{r-1}+t_{r-1}^{\Phi}}{i_{r-1}} 
\binom{a_{r}+t_r^{\Phi}}{i_{r}}f_\Phi({\bf i}).$$

  Applying FORMULAE (1)  and (2)  for the volume and Ehrhart polynomials, we see that 
 the Ehrhart polynomial on a big chamber $\c$ is immediately deduced from the 
  polynomial function $\vol P_\Phi(a)$ (i.e; its highest degree component)
   by replacing the monomial 
  $\frac{a_k^{i_k}}{i_k!}$
  by the function  $\binom{a_{k}+t_k^{\Phi}}{i_{k}}$ (with same leading term).
   The function 
  $a\mapsto \binom{a}{k}$ is more adequate than the monomial $a^k/k!$
  in the integral context, as it takes integral values on integers.

 THEOREM  \ref{vw} is a generalization of  B.V. Lidskii formula for 
Kostant partition function [Li]. It was also proven by Stanley for 
general flow polytopes. R. Stanley  suggested to look for a proof 
via residues.

In Section \ref{pitman}, we apply the ``nice formula" to a 
polytope considered by Pitman and Stanley [Pi-S].

Finally, in Section \ref{sym}, using the action of the full Weyl 
group $\Sigma_{r+1}$ on $S_{A_r}$, we  list $r-1$  symmetries 
properties of the Kostant polynomial 
$$k(A_r^+,\c^+)(a_1(e^1-e^{r+1})+\cdots 
+a_{r-2}(e^{r-2}-e^{r+1})+a_{r-1}(e^{r-1}-e^{r+1})+a_{r}(e^r-e^{r+1})).$$

(This is the mixed lattice point enumerator  $N(\sum a_k P_k)$ of  
the Chan-Robbins-Yuen polytopes $P_k=P_{A_r^+}(e^k-e^{r+1})$,  
$1\leq k \leq r$.) The function $k(A_r^+,\c^+)(a)$ is a function 
of $(r-1)$ variables $k^+(a_1,a_2,....., a_{r-1})$. This function 
satisfies for example, the properties:

\begin{itemize}
     \item 
    for any values $x_1,..., x_{r}$, we have 
    
    $$k^+(x_1,....,x_{r-2},-( x_1+...+x_{r-2}+x_{r-1}+x_r+2))=$$
    $$-
    k^+(x_1,x_2,..., x_{r-2},x_{r-1})-k^+(x_1,...,x_{r-2},x_{r}-1).$$

  \item for any values $x_1,..., x_{r}$, we have 
$$k^+(x_1,x_2,..., x_{r-3},-(x_1+x_2+...+ x_{r}+3),x_{r-1})=$$ $$ 
-k^+(x_1,x_2,..., x_{r-3}, x_{r-2},x_{r}-1) +k^+(x_1,x_2,..., 
x_{r-3}, x_{r-1}-1 , x_{r-2}+1)$$ $$ -k^+(x_1,x_2,..., x_{r-3}, 
x_{r}-2 , x_{r-2}+1)+k^+(x_1,x_2,..., x_{r-3}, x_{r-1}-1,x_{r}-1). 
$$ 
\end{itemize}

We do hope  to  have shown in  these notes  that  the space 
$S_\Delta$, well known in the study of hyperplanes arrangements, 
is  providing an efficient   tool in the computation of Volumes   
and Ehrhart  polynomials of  polytopes. However, we also learnt in 
writing these notes how difficult and exciting are ``concrete 
polytopes''.

\newpage

\section{Total residues}\label{total}

  Let $V$ be a $r$-dimensional real vector space and
$V^*$ its dual vector space. Let $\Delta\subset V^*$ be a finite 
subset of non-zero linear forms. We assume $\Delta$ symmetric: 
$\Delta=-\Delta$.  Each $\alpha\in \Delta$ determines a linear 
form on $V$ and a  complex hyperplane $\{x\in V_\C; \alpha(x)=0\}$ 
in $V_\C$. Consider the hyperplane arrangement $$ 
\CH_{\C}=\bigcup_{\alpha\in \Delta}\{x\in V_\C,\alpha(x)=0\}. $$   
The ring $R_\Delta$ of rational functions with poles on $\CH_{\C}$ 
is the ring $\Delta^{-1}S(V^*)$ generated by the ring $S(V^*)$ of 
polynomial functions on $V$, together with inverses of the linear 
functions $\alpha\in \Delta$. 

A subset $\sigma$ of $\Delta$ is called a {\bf basis of} $\Delta$, 
if the elements $\alpha\in \sigma$ form a basis of $V^*$. We 
denote by $\CB(\Delta)$ the set of bases of $\Delta$. For $\sigma$ 
a basis of $\Delta$, set $$ 
f_{\sigma}(x):=\frac{1}{\prod_{\alpha\in \sigma}\alpha(x)}. $$ 

\begin{definition}
The subspace $S_{\Delta}$ of $R_{\Delta}$ spanned by the elements 
$f_{\sigma},~ \sigma\in\CB(\Delta),$ will be called the space of 
simple elements of $R_{\Delta}$: $$ S_{\Delta}=\sum_{\sigma\in 
 \CB(\Delta)}\R f_{\sigma}. $$ 
\end{definition}

 The vector space $S_{\Delta}$ is contained in the  homogeneous 
component of degree $-r$  of $R_{\Delta}$.
 If $\Delta$ does not span $V^*$, then the set $\CB(\Delta)$ is empty and  $S_\Delta=0$.

In general, elements $f_\sigma$ are not linearly independent.

 {\bf  Example 1.}
 
Let $V$ be a $2$-dimensional vector space with basis $e_1,e_2$. 
Let $\Delta$ be the  set $$\Delta=\{\pm e^1, \pm e^2, \pm 
(e^1-e^2)\}.$$

Then we have the linear relation 
 $$\frac{1}{e^1 e^{2}}= \frac{1}{e^2 (e^1-e^{2})}- 
\frac{1}{e^1 (e^1-e^{2})}$$
 between  elements $f_{\sigma_1}$,$f_{\sigma_2}$, 
$f_{\sigma_3}$  with $\sigma_1=\{e^1, e^2\}$, $\sigma_2=\{e^1, 
(e^1-e^2)\}$ and  $\sigma_3=\{e^2, (e^1-e^2)\}$ bases of $\Delta$.

 We let elements $v$ of $V$ act on $R_{\Delta}$ by 
differentiation: $$ (\partial(v)f)(x):= 
\frac{d}{d\epsilon}f(x+\epsilon v)|_{\epsilon=0}. $$

Then the following holds ([B-V 2], Proposition 7.) 
\begin{theorem}\label{bv2}
 $$ R_\Delta= \partial(V) 
R_{\Delta}\oplus S_{\Delta}. $$

\end{theorem}

 As a corollary of this decomposition, we can define the 
projection map $$ Tres_{\Delta}: R_{\Delta}\to S_{\Delta}. $$ The 
projection $(Tres_{\Delta}f)(x)$ of a function $f(x)$ is again 
 a function of $x$ that we  called the {\bf total residue} of $f$. By 
definition, this function can be expressed as a linear combination 
of the simple fractions $f_{\sigma}(x)$. {\bf  The main property 
of the map $Tres_{\Delta}$ is that it vanishes on derivatives}. If 
$f\in R_\Delta=\frac{P}{\prod_k\alpha^k}$ $(P\in 
S(V^*),\alpha^k\in \Delta$ ) has a denominator product of linear  
forms $\alpha^k\in \Delta$ which do not generate $V^*$, then it is 
easy to see that $f$ is a derivative and  the total residue of $f$ 
is equal to $0$.

{\bf Example 2:}\label{example}
 Let us do a computation.
 Let $V$ be with basis $e_0,e_1,e_2$ and let 
 $$\Delta=\{\pm e^1, \pm e^2,\pm (e^0-e^1),\pm (e^0-e^2),\pm (e^1-e^2)\}.$$

 We write $x\in V$ as $x=x_0 e_0+x_1 e_1+x_2 e_2$. The reason, for writing $x_i$ instead of $x^i$
 is in order not to misinterpret an upper index with a power.
 Consider the following function $W_3$ of $R_{\Delta}$:
 $$W_3(x_0,x_1,x_2)=\frac{x_0^{2}}{(x_0-x_1)(x_0-x_2)(x_1-x_2)x_1x_2}.$$  
  Then $W_3$ is homogeneous of degree $-3$. To compute the total residue of
  $W_3(x_0,x_1,x_2)$, we write $x_0$ as  a linear combination of linear forms 
  in the denominator of $W_3$, 
  in order to reduce the degree  of the denominator. For example,
  writing  $x_0^2= ((x_0-x_1)+x_1)((x_0-x_2)+x_2)$, we  obtain
  
   $$W_3(x_0,x_1,x_2)=\frac{1}{(x_1-x_2)x_1x_2}+  \frac{1}{(x_0-x_2)(x_1-x_2)x_1}$$  
  $$+\frac{1}{(x_0-x_1)(x_1-x_2)x_2}+  
  \frac{1}{(x_0-x_1)(x_0-x_2)(x_1-x_2)}.$$

  The first and last fractions have denominators with linearly dependent forms,
  so that their  total residue is zero and we obtain:

$$Tres_{\Delta}(W_3(x))=\frac{1}{(x_0-x_1)(x_1-x_2)x_2}+\frac{1}{(x_0-x_2)(x_1-x_2)x_1}.$$
  
(More precisely, in the direct sum decomposition,  
  $$R_\Delta=S_\Delta\oplus (\partial_{x_0}R_\Delta+\partial_{x_1}R_\Delta+
  \partial_{x_2}R_\Delta),
  $$  
 we have,
 
 $$W_3(x_0,x_1,x_2)-Tres_{\Delta}(W_3(x))=U_3(x_0,x_1,x_2),$$
 
 with 
 $$U_3(x_0,x_1,x_2)=-\partial_{x_1}\frac{x_0-2x_2}{x_2(x_0-x_2)(x_1-x_2)}-\partial_{x_2}
 \frac{x_0-2x_1}{x_1(x_0-x_1)(x_1-x_2)}.)$$

 If $V$ is one dimensional, and $\Delta=\{\pm e^1\}$, then 
$R_{\Delta}$ is the ring of Laurent series $$L=\{f(x)=\sum_{k\geq 
-q}a_k x^k\}.$$ The total residue of a function $f(x)\in L$ is  
the function $\frac{a_{-1}}{x}$. The usual residue denoted 
$Res_{x=0}f$ is the constant $a_{-1}$. The constant term $a_0$ of 
the Laurent series $f(x)$ is denoted by $CT_{x=0}f$. We have 
$CT_{x=0} xf(x)=Res_{x=0}f(x)$.

We will  also use the following obvious properties.

\begin{lemma}\label{obvious}

\begin{itemize}

\item
Assume $\Gamma\subset \Delta$ is a subset of $\Delta$. Then 
$$R_\Gamma\subset R_\Delta,$$ $$S_\Gamma\subset S_\Delta.$$ 
Furthermore, if $f\in R_\Gamma$, then $Tres_\Delta(f)$ belongs to
$S_\Gamma$ and $$Tres_\Gamma(f)=Tres_\Delta(f). $$

\item
 Assume that $V=V_1\oplus V_2$ and $\Delta=\Delta_1\cup 
\Delta_2$, with $\Delta_i\subset V_i$, then 
$$R_\Delta=R_{\Delta_1}\otimes R_{\Delta_2},$$
$$S_\Delta=S_{\Delta_1}\otimes S_{\Delta_2},$$
$$Tres_\Delta=Tres_{\Delta_1}\otimes Tres_{\Delta_2}.$$ 

\end{itemize}
\end{lemma}

 Consider  the vector space $\hat S(V^*)$  of 
formal power series on $V$. Define $\hat 
R_{\Delta}=\Delta^{-1}\hat S(V^*)$. As the total residue vanishes 
outside the homogeneous component of degree $-r$ of $R_\Delta$, 
the map $Tres_{\Delta}$ extends  as a map  $$ Tres_{\Delta}: \hat 
R_{\Delta}\to S_{\Delta}. $$

Consider the open set  $U_\Delta= V_\C-\CH_{\C}$,  complement of 
the union of hyperplanes $\{\alpha=0\}$.  Choose a basis $e^i$ of 
$V^*$ . This gives  coordinates  $x^i$ on $V_\C$. Let 
$dx=dx^1\wedge dx^2\wedge\cdots \wedge dx^r.$  Let $f\in 
R_\Delta$. Consider the $r$- holomorphic differential form 
$f(x)dx$, defined  on the open set  $U_\Delta$. Then 
$(f-Tres_\Delta (f))dx$ is  an exact form on this open set : by 
definition of the total residue, the function $f-Tres_\Delta (f)$ 
is in the span of $\partial_k R_\Delta$, so $(f-Tres_\Delta( 
f))dx$ is the differential of some $(r-1)$- holomorphic 
differential form $\sum_{k=1}^r f^k dx^1\wedge 
dx^2\wedge\widehat{dx^k}\wedge \cdots \wedge dx^r$.

\bigskip

{\bf Bibliographical remarks.}

The space $S_\Delta$ is isomorphic to the highest degree component 
of the  Orlik-Solomon algebra $H^*(U_\Delta,\C)$, described in   
[O-T] by generators and relations. The definition of the total 
residue is given in [B-V 2]. It formalizes notions  introduced in 
Jeffrey-Kirwan [J-K].  Proposition \ref{bv2} is proven in [B-V 2]. 

\section{Jeffrey-Kirwan residue formula for volumes of convex polytopes}\label{jeffrey}

Now the vector space $V$  is equipped with a fixed Lebesgue 
measure $dx$. We denote by $da$ the dual  measure on $V^*$. If 
$\sigma$ is a basis of $\Delta$, we denote by 
 $\vol(\sigma)$ the volume of the parallelepiped 
 $\oplus_{\alpha\in \sigma} [0,1]\alpha$, for our Lebesgue measure $da$.

We consider $\R^N$  with  basis $(w^1,\ldots,w^N)$ and  
corresponding Lebesgue measure $dw$. Let $p$ be a surjective 
linear map from $\R^N$ to  the vector space $V^*$. Then the vector 
space $\Ker(p)=p^{-1}(0)$ is of dimension $d:=N-r$.   It is 
equipped with a Lebesgue  measure $dm$ quotient of $dw$ by $da$. 
For $a\in V^*$, $p^{-1}(a)$ is an affine space parallel to 
$\Ker(p)$, 
 thus also equipped with the Lebesgue measure $dm$.
 Volumes of subsets of $p^{-1}(a)$ are computed for this measure.
 
  We set 
$p(w^k):=\alpha^k$ for $1\leq k\leq N$, and 
$$\Phi=(\alpha^1,\ldots,\alpha^N).$$

 We have thus, for $u_i\in \R$, $$p(u_1w^1+u_2w^2+...+ u_Nw^N)= 
u_1 \alpha^1+u_2\alpha^2+\cdots +u_N\alpha^N .$$

 The $N$ elements $\alpha^k$ of 
the sequence $\Phi$ {\bf need not to  be distinct}.  We denote by 
$C_N$ the closed convex cone in $\R^N$ generated by 
$w^1,\ldots,w^N$, and we set $C(\Phi):=p(C_N)$, the cone generated 
by $(\alpha^1,\ldots,\alpha^N).$

We assume that $p^{-1}(0)\cap C_N=\{0\}$. Then  $0$ is not in the 
convex hull of the vectors $\alpha^k$  and $C(\Phi)$ is an acute 
cone.

 \begin{definition}
 Let $a\in V^*$.
 We define $$~P_{\Phi}(a):=p^{-1}(a)\cap C_N~.$$
 \end{definition}

We immediately see  that  the set  $P_{\Phi}(a)$ is the convex 
polytope consisting  of all solutions $(u_1,u_2, ..., u_N),$ in 
non negative real numbers $u_k$, of the equation $$\sum_{k=1}^N 
u_k \alpha^k=a.$$ In particular, the polytope $P_{\Phi}(a)$ is 
empty if $a$ is not in the cone $C(\Phi)$. 

 {\em Remark:  any convex polytope can 
be realized canonically as a polytope $P_\Phi(a)$, see for example 
[B-V 1, Section 4.1].  } 

The following lemma is obvious. 

\begin{lemma}\label{fd}
 For any invertible transformation $w$ of $V^*$, then 
$$ P_{w\Phi}(w\cdot a)= P_\Phi(a).$$  
 
 \end{lemma}

 For $a$ in the interior of 
 $C(\Phi)$, the dimension of $P_\Phi(a)$ is $N-r$. 
The function 
 $a\mapsto \vol P_\Phi(a)$ is a continuous function 
 on $C(\Phi)$, homogeneous of degree $N-r$. 
  This function is locally polynomial. Let us recall the description of  regions where 
we are sure that  this function is given by  a polynomial formula.

We consider the set $\{\Phi\}\subset V^*$ consisting of elements 
of $\Phi$: we mean $\Phi$ and $\{\Phi\}$ are the same sets, but 
$\Phi$ may have multiplicities. Let $\Delta^+$ be a finite subset 
of $V^*$ containing the set ${\Phi}$ and such that the cone 
$C(\Delta^+)$ is acute. Usually, we take $\Delta^+=\{\Phi\}$.  We 
define $$\Delta= \Delta^+\cup -\Delta^+.$$

For any  subset  $\nu^+$  of   $\Delta^+$, 
 we denote by $C(\nu^+)$ the closed cone generated by 
 $\nu^+$.
 We denote by 
 $C(\Delta^+)_{sing}$  the union of the cones $C(\nu^+)$
 where $\nu^+$ is any subset of $\Delta^+$ of cardinal strictly less 
 than $r=\dim(V)$. By definition,
the set $C(\Delta^+)_{reg}$ of $\Delta^+$-regular elements is the 
complement of $C(\Delta^+)_{sing}$. 
 A connected component of $C(\Delta^+)_{reg}$ is called a big chamber.
 If  $\c$ is  a big chamber, and $\sigma^+$ a basis 
 of $\Delta^+$, then either $\c\subset C(\sigma^+)$, or 
 $\c\cap C(\sigma^+)=\emptyset$, as the boundary of $C(\sigma^+)$
 does not intersect $\c$. Thus the closure of the  big chamber 
 $\c$ is the intersection of the simplicial cones  $C(\sigma^+)$, $\sigma^+$
 a basis of $\Delta^+$,
  containing $\c$.

  A wall of $\Delta$ is a (real) hyperplane generated 
 by $r-1$ linearly independent elements of $\Delta$. 
 We denote by $\CH^*$ the union of walls. A small chamber of 
 $C(\Delta^+)$ is a connected component of $C(\Delta^+)-\CH^*$.
 Clearly $C(\Delta^+)-\CH^*$ is contained in $C(\Delta^+)_{reg}$
 and each small chamber is contained in a big chamber. 
 In the appendix, we  describe  big and small chambers $\c$ when 
 $\Delta^+$ is the positive root system of 
 $A_n$, for $n=2,3 $. 
  See Figure 4 in Section \ref{appendix2}.

{\em Remark. When $a$ varies in a big chamber  $\c$, the 
combinatorial nature of $P_\Phi(a )$ remains the same, and the 
family of polytopes $P_\Phi(a )$ have parallel facets. The notion 
of Minkowski sums and of big chambers are intimately related:   if 
${\bf a_1},..., {\bf a_m}$ are vectors in the closure of a big 
chamber $\c$, then the polytope $P_\Phi({\bf a_1}+\cdots+{\bf 
a_m})$ is isomorphic to the  Minkowski sum of the polytopes  
$P_\Phi({\bf a_k})$. This follows for example of ([B-V 1 Section 
3.1] )}

\bigskip

 The function  $a\mapsto \vol P_\Phi(a)$ is given by a polynomial formula, when 
 $a$ varies in a big chamber $\c$. Let us recall this formula.

 The Jeffrey-Kirwan residue ([J-K])  associate to a big 
chamber $\c$ of $C(\Delta^+)$ a linear form 
 $f \mapsto \langle\langle \c,f\rangle\rangle$
 on  the vector space $S_{\Delta}$ of simple fractions.
 To determine the linear  map
  $f\mapsto \langle\langle \c,f\rangle\rangle$, it is enough to determine it on the 
 generating set $f_{\sigma}$, with 
 $\sigma$   a basis of $\Delta$. 
  If $\sigma=\sigma^+\cup \sigma^-$ with $\sigma^+\subset 
\Delta^+$ and $\sigma^-\subset \Delta^-$, then $$f_\sigma= 
(-1)^{|\sigma^-|}f_{(\sigma^+\cup\{-\sigma^-\})}.$$ Thus any 
function $f$ in $S_\Delta$ can be written as a linear combination 
of  functions $f_\sigma$, with $\sigma$ a basis of $\Delta$ 
consisting of positive elements. 

  The  linear form  $\langle\langle \c, f\rangle\rangle$ has the following properties:
 
A) If $\sigma$ is a basis of $\Delta^+$ (it is important to 
require that 
 $\sigma$ consists of positive elements),
 then

 a) If $\c\subset C(\sigma)$,
 $$\langle\langle \c, f_{\sigma}\rangle\rangle=\frac{1}{\vol(\sigma)},$$
 
  b) If $\c\cap C(\sigma)=\emptyset$,
 
 $$\langle\langle \c, f_{\sigma}\rangle\rangle=0.$$

B) If $\sigma=\sigma^+\cup \sigma^-$ with $\sigma^+\subset 
\Delta^+$ and $\sigma^-\subset \Delta^-$, then 

$$\langle\langle 
\c,f_\sigma\rangle\rangle=(-1)^{|\sigma^-|}\langle\langle \c, 
f_{(\sigma^+\cup\{-\sigma^-\})}\rangle\rangle.$$

{\em Remark. We can interpret the linear form $f\to \langle\langle 
\c,f\rangle\rangle$ in terms of Laplace transform (see Section 
\ref{appendix1}) . 
 Assume $\sigma$ is a basis of $\Delta$ contained in 
 $\Delta^+$.
  Then we have:
 
 $$\vol(\sigma)f_{\sigma}(x)=\int_{C(\sigma)}e^{-\langle x,a\rangle}da$$
 whenever  $x$ is in the dual cone of $C(\sigma)$. 
 In particular, a function $f$ in $S_\Delta$,  being a linear combination
 of such functions 
 $f_\sigma$, coincide on the dual cone to 
 $C(\Delta^+)$ with the Laplace 
transform of a function $\hat f$ locally constant on big chambers 
of $C(\Delta^+)$, and the linear form $\langle\langle 
\c,f\rangle\rangle$ consists in evaluating the 
 function $\hat f$ at a point of $\c$.}

\bigskip

 Let $a\in V^*$. Let  $f\in R_\Delta$.
  Consider 
 $$e^a f =\sum_{k=0}^{\infty} \frac{a^k}{k!}f.$$
 
This is an element of  $\hat R_\Delta$ . If $e_k$ is a basis of 
$V$, and we write $a\in V^*=\sum_{i=1}^r a_i e^i$ and $x\in V$ as 
$x=\sum_{k=1}^r x_k e_k$ (here, and in many other local 
calculations, we write $x_i$ instead of $x^i$, 
 as we do not want to 
misinterpret  an upper index as a power), then $e^a f$ is the 
analytic function on $V_\C-\CH_{\C}$, given by $x\mapsto 
e^{a_1x_1+a_2 x_2+\cdots +a_r x_r}f(x)$. 
 Its total residue is 
defined  by $$Tres_{\Delta}(e^a f)=\sum_{k=0}^{\infty} 
Tres_{\Delta}(\frac{a^k}{k!}f).$$ There are only a finite number 
of non-zero terms in this sum. More precisely, if $f$ is 
homogeneous of degree $-q$, then $$Tres_{\Delta}(e^a f)= 
\frac{1}{(q-r)!}Tres_{\Delta}(a^{q-r}f).$$

\begin{definition}\label{J}

Define $$J_\Phi(a):= Tres_{\Delta} 
\left(\frac{e^{a}}{\prod_{\alpha\in \Phi}\alpha}\right).$$ 
\end{definition}

 This  is a polynomial function of $a\in V^*$ with value in 
the vector space $S_{\Delta}$.  More precisely,  we choose a basis 
$e^1,..., e^r$ of $V^*$, and write $a=\sum_{i=1}^ra_i e^i$, 
$x=\sum_{i=1}^r x_i e_i$. Let  ${\bf i}$ be a  sequence 
 $(i_1, i_2,..., i_r)$ of non-negative integers. We write $|{\bf i}|:=i_1+i_2+\cdots +i_r$.
 If $|{\bf i}|=N-r$, then  the function 
 $$\frac{x_1^{i_1}x_2^{i_2}\cdots
 x_r ^{i_r}}{\prod_{\alpha\in \Phi}\alpha(x)}$$
 is of homogeneous degree $-r$.

\begin{definition}\label{f}
Let  $f_\Phi({\bf i})$ be the following element 
 of $S_{\Delta}$:

 $$f_\Phi({\bf i})(x)= Tres_{\Delta}\left(\frac{x_1^{i_1}x_2^{i_2}\cdots
 x_r ^{i_r}}{\prod_{\alpha\in \Phi}\alpha(x)}\right).$$ 
  \end{definition}

 \begin{lemma}
 
 We have 
 $$J_\Phi(a)=\sum _{|{\bf i}|=N-r}
 \frac{a_1^{i_1}}{i_1!}\frac{a_2^{i_2}}{i_2!}\cdots 
 \frac{a_{r-1}^{i_{r-1}}}{i_{r-1}!} \frac{a_r^{i_{r}}}{i_{r}!}f_\Phi({\bf i}).$$
 
 \end{lemma}

 \begin{proof}
 We have 
 $$J_\Phi(a)(x)=Tres_{\Delta}\left(
 \frac{e^{a_1x_1+...+a_r x_r}}{\prod_{\alpha\in \Phi}\alpha(x)}\right)$$ 
 $$=
 \sum _{i_1, i_2,..,i_r}
 \frac{a_1^{i_1}}{i_1!}\frac{a_2^{i_2}}{i_2!}\cdots 
 \frac{a_{r-1}^{i_{r-1}}}{i_{r-1}!} \frac{a_r^{i_{r}}}{i_{r}!}Tres_{\Delta}\left(\frac{x_1^{i_1}x_2^{i_2}\cdots
 x_r ^{i_r}}{\prod_{\alpha\in \Phi}\alpha(x)}\right).$$

The total residue is zero outside the homogeneous component of 
degree $-r$. This gives the formula of the lemma.

\end{proof} QED

 Thus the function $a\mapsto J_\Phi(a)$  is an homogeneous polynomial in 
the variable $a$ of degree $N-r$, with value in the finite 
dimensional vector space $S_\Delta$.

\begin{theorem}\label{jk}(Jeffrey-Kirwan) 
 Let $\c$ be a big chamber of $C(\Delta^+)$. 
 Then, if $a\in \overline{\c}$,
 we have 
 $$\vol P_{\Phi}(a)=\langle\langle \c, J_\Phi(a)\rangle\rangle.$$

\end{theorem}

 We  give the proof of this theorem in the appendix (Section\ref{appendix1}).

\begin{definition}\label{vc}
Let $\c$  be   a big chamber. We denote by $v(\Phi,\c)(a)$ the 
polynomial function on $V^*$ such that 

$$v(\Phi,\c)(a)= \vol P_{\Phi}(a)$$ 

when $a$  varies in $\overline{\c}$. 
\end{definition}

 More explicitly, when  $a$ varies in the closure of a  big chamber $\c$, 
 we have  the polynomial formula for the volume: 

\begin{equation}\label{jj}
v(\Phi,\c)(a)= \sum _{|{\bf i}|=N-r}\langle\langle \c, f_\Phi({\bf 
i})\rangle\rangle 
\frac{a_1^{i_1}}{i_1!}\frac{a_2^{i_2}}{i_2!}\cdots 
 \frac{a_{r-1}^{i_{r-1}}}{i_{r-1}!} \frac{a_r^{i_{r}}}{i_{r}!}.
  \end{equation}

{\em Remark. The function $a\to J_\Phi(a)\in S_\Delta$ is a 
vector-valued polynomial function on $V^*$ of degree $N-r$. It is 
thus determined by $\dim S_\Delta$ scalar  homogeneous polynomial 
functions of degree $(N-r)$ in $r$ variables. 
 As a result of THEOREM \ref{jk}, the collection of 
polynomial functions  $v(\Phi,\c)$ is entirely determined by the 
function $J_\Phi(a)$. It is very difficult to determine the number 
of big chambers $\c$. It is usually much larger than $\dim 
S_\Delta$.  For example, consider $A_n^+$ the positive root system 
of $A_n$.  For $n=1,2,3$, the number of big chambers is 
$1$,$2$,$7$; for $n=4$,$5$,$6$, it is known  that their number is   
greater or equal to than $48$, $820$,$51133$ ([S]). The  dimension 
of $S_\Delta$ is $n!$, so respectively $1,2,6,24,120,720.$ So, 
there are many linear relations between the different linear 
functions $f\to \langle\langle \c,f\rangle\rangle$ on $S_\Delta$, 
thus many linear relations between the  polynomials $v(\Phi,\c)$.}

The relative volume  of $P_\Phi(a)$ compare the volume of  
$P_\Phi(a)$ to the  volume of the standard simplex, so is defined 
by $$\vol_{rel} P_{\Phi}(a) =(N-r)! \vol P_\Phi(a).$$  We have: 
$$\vol_{rel} P_\Phi(a)=\langle\langle \c, 
Tres_{\Delta}(\frac{a^{N-r}}{\prod_{\alpha\in 
\Phi}\alpha})\rangle\rangle.$$ 

\bigskip

{\bf Bibliographical remarks.}

Results on this section are mainly due to Jeffrey-Kirwan [J-K]. 
The stratification on big chambers is introduced in [B-V 1]. The 
fact that the volume is polynomial on big chambers follows from 
the proof given in the Appendix ( Section \ref{appendix 1}).  More 
details on the relations between Jeffrey-Kirwan formulae and 
Laplace transforms are given in [B-V 2].

\section{Flow polytopes }\label{flow}

 Consider a $(r+1)$ dimensional vector space, with basis $e^1, 
e^2,..., e^{r+1}$, and  let $A_r^+$ (the positive root system of 
$A_r$) be

$$A_r^+=\{(e^i-e^j), 1\leq i< j\leq (r+1)\}.$$ 

Let $E_{r}$  be the vector space spanned by the elements 
$(e^i-e^j)$,
 $$E_r=\{ a\in \R^{r+1}, a=a_1 e^1+...+a_{r} e^r+a_{r+1} e^{r+1} \,{\rm with}\, \, 
a_1+a_2+...+a_r+a_{r+1}=0\}.$$

 The vector space $E_{r}$ is of dimension $r$ and is 
provided with the lattice generated by $A_r$, so has a canonical 
measure $d{\bf a}$. We have $d{\bf a}= da_1...da_r$. 
 If $\sigma$ is any basis of $A_r$, then $ 
\vol(\sigma)=1$.    The cone  $C(A_r^+)$  generated by positive 
roots is the cone $a_1\geq 0$, $a_1+a_2\geq 0$,...,$ 
a_1+a_2+\cdots +a_r\geq 0$.

\begin{definition}\label{bf}

If $a=(a_1,a_2,..., a_r)\in \R^r$, we denote by ${\bf a}\in E_r$ 
the element $${\bf a}=a_1 e^1+...+a_{r} e^r-(a_1+\cdots + a_r) 
e^{r+1}.$$ 
\end{definition}

Let $\Phi$ be a {\bf sequence} of  $N$ elements of  $A_r^+$ 
spanning $E_r$.

\begin{definition}\label{mjk}

We define  $m_{ij}(i< j)$  to be the multiplicity of the root 
$(e^i-e^j)$ in $\Phi$ .  Thus $m_{ij}\geq 0$.

\end{definition}

 Associate to $\Phi$ 
 a graph  with $r+1$ vertices $1,2,...,r+1$, and  $m_{ij}$ edges  from $i\to j$
 if $(e^i-e^j)$ is in $\Phi$ (thus $i< j$). 
  This graph is called a flow graph. 
  The graph associated to $A_r^+$ is the complete flow graph (with  $1$
  edge $i\mapsto j$  for any $i< j$).

 \begin{definition}

\begin{itemize}
  \item A flow polytope is a polytope isomorphic to a polytope 
 $P_\Phi({\bf a})$ where $\Phi$ is a sequence of elements of 
 $A_r^+$ and  $a=(a_1,a_2,..., a_r)\in \R^r$.

  \item The polytope $P_\Phi(e^1-e^{r+1})$  will be called 
 the  fountain polytope of shape $\Phi$.
 
 \item  
A flow polytope  $P_\Phi({\bf a})$ with $a_1\geq 0,a_2\geq 
0,...,a_r\geq 0$  will be called a cascade polytope.
  
\end{itemize}

\end{definition}
 
  The flow polytope  $P_\Phi({\bf a})$ is of dimension $N-r$
 when $a_1>0$.
 The polytope $P_{A_r^+}(e^1-e^{r+1})$ associated to the complete flow graph
is the Chan-Robbins-Yuen polytope. It will be described in Section 
\ref{CR}.

 {\em Remark 1. If $a_1,..., a_r$ are greater or equal to $0$ , the polytope 
 $P_\Phi({\bf a})$ is the Minkowski sum of 
 the polytopes $a_k P_\Phi(e^k-e^{r+1})$.
If $k>1$, the vector $e^k-e^{r+1}$ is on the boundary of 
$C(\Delta^+)$ and the polytope $ P_\Phi(e^k-e^{r+1})$ 
 is of smaller dimension than $N-r$. It is the fountain  polytope
 of shape  $\Phi'$ where $\Phi'$ is the system where
  we have deleted all $(e^i-e^j)$ with $i< k$.}

{\em Remark 2. 
 Let $d_a$ be the dimension of $P_\Phi({\bf a})$. If 
$s$ is a vertex of $P_\Phi({\bf a})$ and $s+\R^+ v_1$, $s+\R^+ 
v_2$,..., $s+\R^+ v_{d_a}$ are $d_a$ edges  passing through $s$, 
directed by integral vectors $v_i$ with minimal length, then the 
volume of the parallelepiped spanned by $v_k,1\leq k\leq d_a$ is 
equal to $1$. Such polytope is associated to an ample line bundle 
on a smooth toric variety: a limit of Delzant polytopes. 
Reciprocally, any polytope associated to a smooth toric  variety 
can be realized as a  flow polytope $P_\Phi({\bf a})$ when $\Phi$ 
is a sequence of elements of the positive root system $A_r^+$ 
(Szenes: personal communication).}

In the rest of this article, we study volumes and Ehrhart 
polynomials of  flow polytopes. Our method is based on the study 
of the vector space $R_{A_r}$ and of the vector space $S_{A_r}$. 

Consider the system 
 $A_{r-1}=\{(e^i-e^j), 1\leq i< j \leq r\}$.
 
\begin{proposition}\label{bi}
The  map from 
 $$\sum_{i=1}^r S_{A_{r-1}}\otimes \frac{1}{(e^i-e^{r+1})}\to S_{A_r}$$
 given by 
 $$\sum_{i=1}^r f_i \otimes  
 \frac{1}{(e^i-e^{r+1})}\mapsto \sum_{i=1}^r f_i \frac{1}{(e^i-e^{r+1})}$$
  is a bijection.  
 \end{proposition}

 \begin{proof}
 If 
 $f$ is in 
 $S_{A_{r-1}}$, then for any 
 $1\leq i\leq r$, the element 
 $f  \frac{1}{(e^i- e^{r+1})}$
 is in 
 $S_{A_r}$, and it is easy to see that the map above is injective.

 To prove that it is surjective, observe first that 
 for any set of elements $A\subset \{1,...,r\}$: 
  $$\prod_{i\in A}\frac{1}{(e^{r+1}-e^i)}=
  \sum_{i\in A}\frac{1}{\prod_{j\in A; i\neq j}(e^i-e^j)}
   \frac{1}{(e^{r+1}-e^i)}.$$
  If $\nu$ is a subset of $A_{r-1}$, we denote by 
 $f_\nu=\frac{1}{\prod_{\alpha\in \nu}\alpha}$.
 Then  $$\prod_{i\in A}\frac{1}{(e^{r+1}-e^i)}f_\nu=
 \sum_{i\in A}\frac{1}{\prod_{j\in A; i\neq j}(e^i-e^j)}f_\nu \frac{1}{(e^{r+1}-e^i)}.$$

 If $\nu\cup \{(e^i-e^{r+1}),i\in A\}$ is a basis of 
 $A_r$, then for every $i\in A$, 
 $$\nu\cup  \{(e^j-e^{r+1}),j\in A, j\neq i\}\cup \{(e^i-e^{r+1})\}$$
 is a basis of $A_{r}$. The lemma follows. 
 
 \end{proof} QED

It follows from the lemma above that the dimension of $S_{A_r}$ is 
$r!$. 
 We denote by $\Sigma_r$ the set of permutations of 
$\{1,2,...r\}$.  As seen from PROPOSITION \ref{bi} above, a 
particularly nice  basis  of $S_{A_r}$ is given by the elements 
$$f_w=w\cdot\frac{1}{(e^1-e^2)(e^2-e^3)\cdots 
(e^{r-1}-e^{r})(e^r-e^{r+1})}$$ where $w\in \Sigma_r$ acts by 
permutation on the set $\{1,2,...r\}$. 

We define as usual the character $\epsilon(w)=\pm 1$ of $\Sigma_r$
with value $-1$ on symmetries.

\bigskip

{\bf  Bibliographical remarks.} Polytopes associated  to subsets 
of  $A_r^+$ are related to graphs and  called flow polytopes by  
Stanley [S]. The representation of the Weyl group in the space  
$S_{\Delta}$, for $\Delta$ any root system, has been studied 
extensively by G. Lehrer [L]. In particular, for $A_r$, the space 
$S_{\Delta}$ carries the regular representation of the subgroup 
$\Sigma_r$ of the Weyl group $\Sigma_{r+1}$. We employ the 
corresponding basis $f_w$ ($w\in \Sigma_r$) in these notes. 

\section{Chambers and Iterated residues for $A_r$.}\label{chambers}

Let  $$V=\{\sum_{i=1}^{r+1}x_ie_i,x_{r+1}=0\}.$$ (as before, we 
write $x_i$ instead of $x^i$).  The vector space $V$ has basis 
$e_1,..., e_r$. An element $a=\sum_{i=1}^{r+1}a_ie^i$ of $E_r$ 
gives the   linear form $\sum_{i=1}^{r}a_i x_i$ on $V$. This 
identifies $E_r$ with $V^*$.  In this identification  
 $$A_r^+=\{ x_i, 1\leq i\leq r, (x_i-x_j), 1\leq 
i< j\leq r\}.$$ A function in $R_{A_r}$ is  thus a rational 
function  $f(x_1,x_2,...,x_r)$ on $V_\C$, with  poles on the 
hyperplanes $x_i=x_j$ or $x_i=0$. The basis $f_w$  considered in 
the preceding section is given by the elements 
$$f_w(x_1,...,x_r)=w\cdot\frac{1}{(x_1-x_2)(x_2-x_3)\cdots 
(x_{r-1}-x_{r})x_r}$$ where $w\in \Sigma_r$ acts by permutation on 
the set $\{1,2,...r\}$.

The following  linear form $f\mapsto Ires_{x=0}f$  defined by 
$$Ires_{x=0}f$$ $$=Res_{x_1=0}Res_{x_2=0}\cdots 
Res_{x_r=0}f(x_1,x_2,..., x_r)$$ is a linear form on $R_{A_r}$ 
which vanishes on the vector space of derivatives $\sum_{i=1}^r 
\partial_i R_{A_r}$. It will be called the iterated residue.
It provides a linear form on $S_{A_r}$. If we compute it on the 
basis $f_w$ of $S_{A_r}$ indexed by the symmetric group, we have 
$Ires_{x=0}(f_w)=\delta_w^1$. The iterated residue depends on the 
order on variables. Permuting the variables by the group 
$\Sigma_r$, we obtain $r!$ linear forms 
 on $S_{A_r}$, dual to the basis $f_w$. Precisely set
$$Ires_{x=0}^{\sigma}f=Res_{x_{\sigma(1)}=0}Res_{x_{\sigma(2)}=0}\cdots 
Res_{x_{\sigma(r)}=0}f(x_1,x_2,..., x_r)=$$
$$Res_{x_1=0}Res_{x_2=0}\cdots 
 Res_{x_r=0}f(x_{\sigma^{-1}(1)},x_{\sigma^{-1}(2)},..., x_{\sigma^{-1}(r)}),$$
  it is not difficult to check that $\sigma\cdot Ires_{x=0}f_w=  Ires_{x=0}\sigma^{-1}f_w= 
  Ires_{x=0}^{\sigma} f_w=\delta_w^{\sigma}$ while 
$Ires_{x=0}^1=Ires_{x=0}.$

Iterated residues are particularly easy to calculate, thus it is 
important to express the linear forms associated to big chambers 
in terms of iterated residues.

Consider the set $\CH^*$ of hyperplanes in $E_r$ spanned by 
$(r-1)$ linearly independent  vectors of $A_r$.  
\begin{lemma}

 An element $H $ of $\CH^*$ is the kernel of a linear form 
$ \sum_{i\in A} a_i $ where $A$ is a subset of  ${1,2,...,r}$. 
\end{lemma}

\begin{proof}
Let $\nu$ be a subset of $A_r^+$ spanning a $(r-1)$-dimensional 
vector space.  If $\nu$ is contained in  the set $ \pm {(e^i-e^j), 
1\leq i< j\leq r\ }$, the hyperplane spanned by $\nu$ is the 
hyperplane $\sum_{i=1}^r  a_i=0$. If not, then  $\nu $ contains a  
vector   $(e^i-e^{r+1})$, and  we conclude by induction. 

\end{proof} QED

 \begin{definition}\label{c}
 
 We denote by 
 $\c^+$ the  open set  $$\c^+=\{a
\in C({A_r}^+)\,\,\text{such that} \,\,a_i>0 ,\,i=1,..,r\}.$$ 
\end{definition}

\begin{lemma}

 The set $\c^+$ is a  small and a big chamber for $A_r^+$.
 {\bf It will be called the nice chamber.} 

 \end{lemma}
 
 \begin{proof}
The set $\c^+$ is a  small chamber  since it doesn't  meet any 
hyperplane. But $\overline{\c^+}$ is the simplicial cone generated 
by  the elements  $(e^i-e^{r+1}),\  1\leq i\leq r,$ so that $\c^+$  
is also a   big  chamber. 
 \end{proof} QED

 Let $w\in \Sigma_r$ and $n(w)$ be the number of elements such 
that $w(i)>w(i+1)$. 
 We denote by $C_w^+\subset C({A_r}^+)$ the 
simplicial cone generated by the vectors   $$\epsilon(1) 
(e^{w(1)}-e^{w(2)}),\epsilon(2) 
(e^{w(2)}-e^{w(3)}),...,\epsilon(r-1)  (e^{w(r-1)}-e^{w(r)}), 
(e^{w(r)}-e^{r+1}),$$ where $\epsilon(i)$ is 1 or -1 depending 
whether $w(i)<w(i+1)$ or not.  When $w=1$, then $C_1=C(A_r^+)$.  
The dual basis for the vectors generating $C_w^+$  is given by 
$$\epsilon(1)e_{w(1)}  ,\epsilon(2) 
(e_{w(1)}+e_{w(2)}),...,\epsilon(r-1)  (e_{w(1)}+..+e_{w(r-1)}), 
(e_{w(1)}+..+e_{w(r)}). $$ 

If we write $a=\sum_{j=1}^{r+1} a_je^j$ in $E_r$, then the cone 
$C_w^+$ is  given by the following system of inequalities 
$\sum_{j=1}^i  a_w(j)\geq 0$, for all $i $ such that $w(i)<
w(i+1)$, but $\sum_{j=1}^i a_w(j)\leq 0$ if $w(i)>w(i+1)$.

\begin{lemma}\label{resw}

 Let $\c$ be a big chamber. Consider the set of elements $w\in \Sigma_r$ 
such that $\c\subset C_w^+$. Then, for $f\in S_{A_r}$,

$$\langle\langle \c,f\rangle\rangle=\sum_{w\in \Sigma_r, \c\subset 
C_w^+}(-1)^{n(w)}  Ires_{x=0} w^{-1} f.$$ 
\end{lemma}

\begin{proof}
 Let $f\in S_{{A_r}}$ and write $f=\sum_{w\in \Sigma_r} 
a_w^{f}f_w$. 
 Using the dual base, we 
compute $a_s^{f}={Ires_{x=0}}^s (f)$ and we find  that the linear 
form $f\mapsto \langle\langle \c,f\rangle\rangle$ is the sum over 
all elements   $w\in \Sigma_r$  such that $\c\subset C_w^+$  of 
the linear form $(-1)^{n(w)} Ires_{x=0}^w=(-1)^{n(w)}  
Ires_{x=0}^w$. 
\end{proof} QED

\begin{theorem}\label{tt}

The linear form \  $f\mapsto \  \langle\langle 
\c^+,f\rangle\rangle$ coincide with the iterated residue 
$Ires_{x=0}$. \end{theorem} 

\begin{proof}
As we just  observed, the linear form $f \mapsto \langle\langle 
\c^+,f\rangle\rangle$ is the sum over all elements $w$ such that 
$\c^+\subset C_w^+$  of the linear form $(-1)^{n(w)} 
Ires_{x=0}^w$. So, to prove the lemma, it is enough to prove 
 that if $\c^+ \subset C_w^+$, then $w=1$. Assume that $w\neq 1.$
We let $k$ be the smallest integer such that $w(k)>w(k+1)$.  In 
particular if $a\in C_w^+$ then $$a_{w(1)}\geq 0, 
a_{w(1)}+a_{w(2)}\geq 0,\ldots 
,a_{w(1)}+a_{w(2)}+\ldots+a_{w(k-1)}\geq 0,$$ 
$$a_{w(1)}+a_{w(2)}+\ldots +a_{w(k)}\leq 0$$ which forces 
$a_w(k)\leq 0$, and therefore, if $w\neq 1$, then $\c^+$ is not 
contained in $C_w^+.$ 

\end{proof} QED

\bigskip

{\bf Bibliographical remarks.} 

Bases of $S_\Delta^*$ given by iterated residues are constructed 
for any system $\Delta$ by Szenes [Sz].  The linear forms 
$f\mapsto \  \langle\langle \c,f\rangle\rangle$ are called 
Jeffrey-Kirwan residues. It is usually difficult to express them 
in functions of iterated residues.

\section{Volumes of flow polytopes}\label{volume}
 
Let  $\Phi$ be a sequence  of $N$  vectors in $A_r^+$ spanning 
$E_r$. The set of regular elements for the system $A_r^+$ is  
clearly contained in the set of regular elements for the smaller 
system $\{\Phi\}$.  Recall the map $a=(a_1,a_2,...,a_r)\to 
{ \bf a}$ given in DEFINITION \ref{bf}. We will sometimes identify 
implicitly $E_r$ with $\R^r$ using this map.   
 Thus if $\c$ is a big chamber for $A_r^+$, the 
function $a\mapsto \vol P_\Phi({\bf a})$ is polynomial on $\c$. We 
denote it by $v(\Phi,\c)$. We have (see  Section \ref{jeffrey}, 
Formula \ref{jj}):

$$v(\Phi,\c)(a) = \langle\langle \c, J_\Phi(a)\rangle\rangle.$$

 The function $v(\Phi,\c^+)$, attached to the nice chamber 
$\c^+$ is particularly important, due to the following lemma.

\begin{lemma}\label{www}
Let $\c$ be a big chamber. Let  $\Sigma_r(\c)$ be the set of 
elements $w\in \Sigma_r$ such that $\c\subset C_w^+$. Then 
$$v(\Phi,\c)(a)=\sum_{w\in \Sigma_r(\c)}(-1)^{n(w)}  
v(\Phi,\c^+)(w^{-1} a).$$ 
 
 \end{lemma}

\begin{proof}
Due to THEOREM \ref{tt}, the linear form $f\mapsto \langle\langle 
\c^+,f\rangle\rangle$ is the iterated residue $Ires_{x=0}$. Thus 
$v(\phi,\c^+)(a)=Ires_{x=0}J_\Phi(a)$. By LEMMA \ref{resw}, the 
form $\langle\langle \c,f\rangle\rangle$ is a signed sum of 
iterated residues $Ires_{x=0}^w$ over the  elements $w$ in 
$\Sigma_r(\c)$ and we obtain the lemma. 
\end{proof} QED

We write 

$$v(\Phi,\c)(a_1,..., a_r) = \sum _{i_1+i_2+\cdots 
+i_r=N-r}v(\Phi,\c,{\bf i}) 
\frac{a_1^{i_1}}{i_1!}\frac{a_2^{i_2}}{i_2!}\cdots 
 \frac{a_{r-1}^{i_{r-1}}}{i_{r-1}!} \frac{a_r^{i_{r}}}{i_{r}!}$$
 where $$v(\Phi,\c,{\bf i})=\langle\langle \c,f_\Phi({\bf i})\rangle\rangle.$$

\begin{lemma}

The coefficient $v(\Phi,\c^+,{\bf i})$ is the  non negative 
integer defined by: $$v(\Phi,\c^+,{\bf i})=Ires_{x=0} 
\frac{x_1^{i_1}...x_r^{i_r}}{\prod_{1\leq i< j\leq 
r}(x_i-x_j)^{m_{ij}}\prod_{1\leq i\leq r} x_i^{m_{i,r+1}}}.$$ 
\end{lemma}

\begin{proof}
It is immediate to check that  $$v(\Phi,\c^+,{\bf i})=Ires_{x=0} 
\frac{x_1^{i_1}...x_r^{i_r}}{\prod_{1\leq i< j\leq 
r}(x_i-x_j)^{m_{ij}}\prod_{1\leq i\leq r} x_i^{m_{i,r+1}}}$$ is a 
non negative integer. We will give  a combinatorial  
interpretation of this integer in LEMMA \ref{ppp} in Section 
\ref{eh}.

\end{proof} QED

\begin{proposition}\label{petit}

The polynomial   $v(\Phi,\c^+)(a_1,a_2,..., a_r)$ is homogeneous 
of degree $|\Phi|-r$. It  is divisible by $a_1^{(\sum_{k=2}^{r+1} 
m_{1k}-1)}$. It is of degree less or  equal to  $(m_{r,r+1}-1)$ in 
$a_r$. 

\end{proposition}

\begin{proof}
We have $$J_{\Phi}(a)(x)= \frac{e^{a_1 x_1+...+a_{r-1}x_{r-1}+ a_r 
x_r}}{\prod_{1\leq i< j\leq r}(x_i-x_j)^{m_{ij}}\prod_{1\leq i\leq 
r-1} x_i^{m_{i,r+1}} x_r^{m_{r,r+1}}}$$

$$ =\frac{e^{a_1 x_1+...+a_{r-1}x_{r-1}+ a_r 
x_r}}{x_1^{m_{1,r+1}}\prod_{ 1< j\leq 
r}(x_1-x_j)^{m_{1j}}\prod_{2\leq i< j\leq r}(x_i-x_j)^{m_{ij}} 
\prod_{2\leq i\leq r} x_i^{m_{i,r+1}}}.$$

The function  $v(\Phi,\c^+)(a_1,a_2,..., a_r)$  is  the iterated 
residue of $J_\Phi(a)(x)$.
 If we start by taking the 
residue in $x_r=0$, (use the first expression), we have to develop 
the term $e^{a_r x_r}$ up to order $m_{r,r+1}-1$, so we obtain the 
second 
 property.

On  the other hand, $x_1$ is considered as generic up to the last 
step. Consider  the function $$F(x_1, a_2,..., a_r)=$$ $$ 
Ires_{x_2=0,...,x_r=0} \frac{e^{a_2 x_2+...+a_{r-1}x_{r-1}+ a_r 
x_r}}{\prod_{2\leq j\leq r}(1-x_j/x_1)^{m_{1j}}\prod_{2\leq i< 
j\leq r}(x_i-x_j)^{m_{ij}} \prod_{2\leq i\leq r} x_i^{m_{i,r+1}} 
}.$$ 

This function is  of the form $\sum_{k\geq 0} C_k(a_2,..., a_r) 
x_1^{-k}$, and  $$Ires_{x=0}J_\Phi(a)(x)= Res_{x_1=0}\frac{e^{a_1 
x_1}}{ x_1^{N_1}}F(x_1,a_2,..., a_r),$$ with $N_1=\sum_{k=2}^{r+1} 
m_{1k}$. 

 Thus, 
expanding the exponential, we obtain  

$$v(\Phi,\c^+)(a_1,a_2,...., a_r)=\sum_{k\geq 0}C_k(a_2,..., 
a_r)\frac{a_1^{N_1+k-1}}{(N_1+k-1)!}.$$ 

This establishes the first property . 

\end{proof} QED

\begin{proposition}\label{division}( Schmidt-Bincer)

Let $\c^+$ be the nice chamber. Let $v_r^+=v(A_r^+,\c^+)$. The 
function $v_r^+(a_1,a_2,..., a_r)$ is independent of $a_r$. It is 
of homogeneous degree $r(r-1)/2$,  of degree less than $1$ in the 
variable $a_{r-1}$ and is divisible by $a_1^{r-1} (
a_1+a_2+a_3+...+a_{r-2}+3a_{r-1}) $. 
 
 More precisely, we have:

$$3v_r^+(a_1,a_2,..., a_r)$$ $$= (a_1+a_2+...+a_{r-2}+3a_{r-1} )  
v(A_r^+ \, {\rm minus}\,(e^{r-1}-e^r),\c^+)(a_1, a_2,..., 
a_{r-2},0,0).$$

 \end{proposition}

\begin{proof}
Let $$J(a,x)= \frac{e^{a_1 x_1+...+a_{r-1}x_{r-1}+ a_r 
x_r}}{\prod_{1\leq i< j\leq r-1}(x_i-x_j)\prod_{1\leq i\leq 
r-1}(x_i-x_r) \prod_{1\leq i\leq r} x_i}.$$ 

We have 

$$v_r^+(a_1,a_2,..., a_r)=Ires_{x=0}J(a,x).$$

We first  take the residue in $x_r=0$. We obtain 

$$Res_{x_r=0}J(a,x) =\frac{e^{a_1 x_1+...+a_{r-1} 
x_{r-1}}}{\prod_{1\leq i\leq r-1} x_i^2\prod_{1\leq i< j\leq 
r-1}(x_i-x_j)} $$ 

$$ =\frac{e^{a_1 x_1+...+a_{r-1} x_{r-1}}}{\prod_{1\leq i\leq r-2} 
x_i^3\prod_{1\leq i< j\leq r-2}(x_i-x_j)\prod_{1\leq i\leq r-2}(1- 
x_{r-1}/x_i)x_{r-1}^2}. $$ 

This shows already that $v_r^+(a_1,a_2,..., a_r)$ is independent 
of $a_r$.We proceed to take the residue in $x_{r-1}=0$. There is a 
double pole in $x_{r-1}$, so that the dependence in $a_{r-1}$ is 
of degree at most $1$. More precisely, 

$$Res_{x_{r-1}=0}[\frac{e^{a_1 x_1+...+a_{r-1} 
x_{r-1}}}{\prod_{1\leq i\leq r-2} x_i^3\prod_{1\leq i< j\leq 
r-2}(x_i-x_j)\prod_{1\leq i\leq r-2}(1- x_{r-1}/x_i)x_{r-1}^2}] $$  

$$=\frac{e^{a_1 x_1+...+a_{r-2}x_{r-2}}}{\prod_{1\leq i\leq r-2} 
x_i^3\prod_{1\leq i< j\leq 
r-2}(x_i-x_j)}(a_{r-1}+\sum_{i=1}^{r-2}\frac{1}{x_i})$$

$$=(a_{r-1}+\frac{1}{3}(a_1+a_2+...+a_{r-2}))\frac{e^{a_1 
x_1+...+a_{r-2}x_{r-2}}}{\prod_{1\leq i\leq r-2} x_i^3\prod_{1\leq 
i< j\leq r-2}(x_i-x_j)}$$ 
$$-\frac{1}{3}(\partial_1+\partial_2+\cdots 
+\partial_{r-2})\frac{e^{a_1 x_1+...+a_{r-2}x_{r-2}}}{\prod_{1\leq 
i\leq r-2} x_i^3\prod_{1\leq i< j\leq r-2}(x_i-x_j)},$$ 

as $(\partial_1+\partial_2+\cdots +\partial_{r-2})$  annihilates 
functions $x_i-x_j$.

Residues vanishes on derivatives, so that we obtain  
$$3v_r^+(a_1,a_2,..., a_r)= (3a_{r-1}+ 
a_1+a_2+a_3+...+a_{r-2})\times$$ $$ Res_{x_1=0}...Res_{x_{r-2}=0} 
\frac{e^{a_1 x_1+...+a_{r-2}x_{r-2}}}{\prod_{1\leq i\leq r-2} 
x_i^3\prod_{1\leq i< j\leq r-2}(x_i-x_j)} .$$

On the other hand, the residue computation of  
 $$v(A_r^+
\, {\rm minus}\,(e^{r-1}-e^{r}),\c^+)(a_1, a_2,..., a_{r-2},0,0)$$ 
gives $$v(A_r^+  \, {\rm minus}\,(e^{r-1}-e^{r}),\c^+)(a_1, 
a_2,..., a_{r-2},0,0)$$ $$=Res_{x_1=0}...Res_{x_{r-2}=0} 
\frac{e^{a_1 x_1+...+a_{r-2}x_{r-2}}}{\prod_{1\leq i\leq r-2} 
x_i^3\prod_{1\leq i< j\leq r-2}(x_i-x_j)} $$ as the  step 
$x_{r}=0$ as well as the step $x_{r-1}=0$ involves only simple 
poles, and we obtain the divisibility property announced.

 \end{proof} QED

{\bf Example.} 
  
The dimension of the polytope  $P_{A_2^+}(a)$ is  $1$, 
 of $P_{A_3^+}(a)$ is   $3$ and of $P_{A_4^+}(a)$ is $6$. Thus the corresponding 
polynomials $v(A_2^+,\c^+), v(A_3,\c^+), v(A_4,\c^+)$ are 
homogeneous of degrees $1,3, 6$.

 We have :
 
 $$v(A_2^+,\c^+)(a_1,a_2)=a_1,$$

$$v(A_3^+,\c^+)(a_1,a_2,a_3)=\frac{1}{3!}a_1^3+\frac{1}{2}a_1^2a_2=\frac{1}{6}a_1^2(a_1+3a_2),$$

$$v(A_4^+,\c^+)(a_1,a_2,a_3,a_4)=\frac{1}{120}a_1^3(a_1+a_2+3a_3)(a_1^2+5a_1a_2+10 
a_2^2).$$ 

 More generally, we have the following proposition, with same proof. 
  
 \begin{proposition}\label{ddd}
 Let $\Phi$ be a sequence of  $N$ vectors in 
 $A_r^+$, generating $E_r$. Assume $m_{r,r+1}=1$ and  $m_{r-1,r+1}+m_{r-1,r}=2$.
 Furthermore, assume that $$\frac{m_{j,r+1}+m_{j,r}+m_{j,r-1}}{m_{j,r-1}}=c$$
 is independent of $j$ for $1\leq j\leq (r-2)$, then 
 
 $$ v(\Phi,\c^+)(a_1,..., a_{r-1},a_r)= v(\Phi,\c^+)(a_1,..., a_{r-1},0)$$
 $$=
 (\frac{a_1+\cdots +a_{r-2}}{c}+ a_{r-1})v(\Phi \, {\rm minus}\,(e^{r-1}-e^{r}),\c^+)
 (a_1,a_2,..., a_{r-2},0,0).$$
\end{proposition}

Consider the permutation $w_0: [1,2,...,r, r+1]\mapsto [r+1,r,..., 
2,1]$. Then $-w_0$ preserves  $A_r^+$ and the vector 
$e^1-e^{r+1}$.   It transforms $\Phi=A_r^+ \, {\rm minus}\,
(e^{r-1}-e^{r})$ in $A_r^+ \, {\rm minus}\,(e^{2}-e^{3})$. Thus we 
obtain from LEMMA \ref{fd} 
\begin{corollary} 
We have $$3 v(A_r^+,\c^+)(e^1-e^{r+1})=v(A_r^+\,{\rm minus}\, 
(e^{r-1}-e^{r}),\c^+ )(e^1-e^{r+1})$$  $$=v(A_r^+\,{\rm minus}\, 
(e^{2}-e^{3}),\c^+)(e^1-e^{r+1}).$$

\end{corollary}

\bigskip

{\bf Bibliographical remarks.}

Proposition \ref{division} is due to [S-B].

 \section {Calculation of some total residues  for the system  $A_{r+1}$}\label{calculation}

Consider a  $(r+2)$  dimensional real vector space, with basis 
$e^0, e^1,e^2,...., e^r, e^{r+1}$. In this Section, we consider 
$A_{r+1}$  as the collection of elements $(e^i-e^j)$ with $0\leq i 
\leq r+1,$\,$ 0\leq j\leq r+1$, and $i\neq j$. 

 Let $$\Pi =\{(e^0-e^1),(e^1-e^2),...,(e^{r}-e^{r+1})\}$$ be the  set of 
simple roots. Let $$f_{\Pi}=\frac{1}{(e^0-e^1)(e^1-e^2)\cdots 
(e^{r}-e^{r+1})}.$$ 

This is an element of $S_{A_{r+1}}$ (in particular  is homogeneous 
of degree $ -(r+1)$).

 As explained in the introduction,  we are particularly 
interested in the  $\Sigma_r$   anti-invariant function of  
homogeneous degree $-(r+1)$ given by  $$W_{r+1}= 
\frac{(e^0-e^{r+1})^{r(r+1)/2}}{\prod_{0\leq i< j\leq 
r+1}(e^i-e^j)}.$$  

As seen by the Example 2 of Section \ref{total}, this function is 
not in the space $S_\Delta$. However, its projection on $S_\Delta$ 
is particularly nice. 

We prove in this section the following Theorem

\begin{theorem}\label{morris}

We have
 $$Tres_{A_{r+1}^+}(\frac{(e^0-e^{r+1})^{r(r+1)/2}}{\prod_{0\leq 
i< j\leq r+1} (e^i-e^j)})$$ $$=\prod_{i=1}^{r-1} \frac{(2i)!}{i! 
(i+1)!} \sum_{w\in \Sigma_r}\epsilon(w)w\cdot 
\frac{1}{(e^0-e^1)(e^1-e^2)\cdots (e^{r-1}-e^{r})(e^r-e^{r+1})}.$$ 

\end{theorem}

In fact, we will  prove  more general identities, which are 
reformulation of Morris identity.

Consider the  group $\Sigma_r$ of permutations of $\{1,..., r\}$. 
Let $0\leq \ell\leq r$ and denote by $P_{\ell,r}$ the 
$\Sigma_r$-invariant polynomial 
 $$P_{\ell,r}=\sum_{w\in \Sigma_r} w\cdot[(e^1-e^{r+1})(e^2-e^{r+1}) 
...(e^\ell-e^{r+1})].$$ 

In particular  $$P_{0,r}=r! \hspace{2cm} P_{r,r}=r!\prod_{j=1}^r 
(e^j-e^{r+1}).$$ 

 When $r$ is fixed, we will write $P_\ell$ for $P_{\ell,r}$, 
when $\ell>0$.

We consider the rational function  given by 

$$\phi_{r+1}(\ell,k_1,k_2,k_3)$$ $$=\frac{P_{\ell}}  
{(\prod_{j=1}^r(e^j-e^{r+1}))^{k_1} 
 (\prod_{j=1}^ r(e^0-e^j))^{k_2} (\prod_{1\leq i< j\leq r}(e^i-e^j))^{k_3}}.$$

In particular  

$$\phi_{r+1}(0,k_1,k_2,k_3)$$ $$=r!\frac{1}  
{(\prod_{j=1}^r(e^j-e^{r+1}))^{k_1} 
 (\prod_{j=1}^ r(e^0-e^j))^{k_2} (\prod_{1\leq i< j\leq r}(e^i-e^j))^{k_3}}.$$

\bigskip

 {\bf Here, $k_1$, $k_2$ and $k_3$ are non negative 
integers}, so that $\phi_{r+1}(\ell, k_1,k_2,k_3)$ is an element 
of $R_{A_{r+1}}$  of homogeneity degree  
 $\ell-(k_1+k_2)r-k_3\frac{r(r-1)}{2}$.
 If $k_3$ is odd, this function is anti-invariant  under the group 
$\Sigma_r$ of permutations of $\{1,..., r\}$. If $k_3$ is even, 
this function is invariant.

  Remark that if $k_1\geq 1$,
$$\phi_{r+1}(r, k_1,k_2,k_3)=\phi_{r+1}(0,k_1-1,k_2,k_3).$$

Let  $k_1,k_2,k_3\geq 0$, $0\leq \ell\leq r$. Let $D=(k_1+k_2)r+ 
k_3\frac{r(r-1)}{2}-\ell$. 

 Then the function $(e^0-e^{r+1})^{D-(r+1)}\phi_{r+1}(\ell, k_1,k_2,k_3) $
 is of homogeneity degree equal to $-(r+1)$. 

We have in particular  
$$W_{r+1}=\frac{1}{r!}(e^0-e^{r+1})^{r(r+1)/2}\phi_{r+1}(0,1,1,1).$$

 \begin{theorem}\label{ct}
 Let  $k_1,k_2,k_3\geq 0$, 
$0\leq \ell\leq r$. Let $D=(k_1+k_2)r+ k_3\frac{r(r-1)}{2}-\ell$. 

 Then the function $(e^0-e^{r+1})^{D-(r+1)}\phi_{r+1}(\ell, k_1,k_2,k_3) $
 is of homogeneity degree equal to $-(r+1)$, and we have 

\begin{itemize}
  \item 
  If $k_3$ is odd,
  $$Tres_{A_{r+1}^+}((e^0-e^{r+1})^{D-(r+1)}\phi_{r+1}(\ell, k_1,k_2,k_3))$$
  $$ =C_{r+1}(\ell, k_1,k_2,k_3)
  [ \sum_{w\in \Sigma_r}\epsilon(w) 
w.f_{\Pi}] .$$

    \item 
   If $k_3$ is even,
  $$Tres_{A_{r+1}^+}((e^0-e^{r+1})^{D-(r+1)}\phi_{r+1}(\ell, k_1,k_2,k_3))$$
  $$=C_{r+1}(\ell,k_1,k_2,k_3)[\sum_{w\in 
\Sigma_r} w.f_{\Pi} ].$$ 
   
\end{itemize}

The constants $C_{r+1}(\ell,k_1,k_2,k_3)$ are determined uniquely 
by the relations:

\begin{itemize}

  \item  for $1\leq \ell\leq r,$ 
 
$$(k_1+k_2-2+\frac{k_3}{2}(2r-\ell-1))C_{r+1}(\ell,k_1,k_2,k_3)$$ $$=(k_1-1+\frac{k_3}{2}(r-\ell)) 
 C_{r+1}(\ell-1,k_1,k_2,k_3), $$ 

  \item $$C_{r+1}(r,k_1,k_2,k_3)= C_{r+1}(0,k_1-1,k_2,k_3),$$

  \item   $$C_{r+1}(r-1,1,k_2,k_3)=C_r(0,k_3,k_2,k_3).$$

\item $$C_{r+1}(0,k_1,k_2,k_3)=C_{r+1}(0,k_2,k_1,k_3).$$

\item $$C_{r+1}(0,1,1,0)=r!.$$

\item If $k_1$ or $k_2=0$,
 $$C_{r+1}(\ell,k_1,k_2, k_3)=0.$$

\end{itemize}

\end{theorem}

{\em Remark 1. The function 
$(e^0-e^{r+1})^{D-(r+1)}\phi_{r+1}(\ell, k_1,k_2,k_3)$ is 
invariant or anti-invariant  under the group $\Sigma_r$  depending 
on the parity of $k_3$, then its  total residue must be an element  
of $S_{A_{r+1}}$ which is invariant or anti-invariant by 
$\Sigma_r$. There are $(r+1)$ linearly  independent such 
functions. Let us consider the basis $w\cdot f_{\Pi}$ of  
$S_{A_{r+1}}$ with $w$ a permutation of $\{0,1,2,..., r\}$. For 
homogeneity reasons, it is easy to see that the iterated residue  
$Ires^{\sigma}_{x=0}$ of the function  
$x_0^{D-(r+1)}\phi_{r+1}(\ell, k_1,k_2,k_3)(x_0,x_1,...., x_r,0)$ 
is equal to $0$, if the permutation  $\sigma$ of $\{0,1,....r\}$ 
does not leave $0$ fixed. Thus the total residue of  
$(e^0-e^{r+1})^{D-(r+1)}\phi_{r+1}(\ell, k_1,k_2,k_3) $ belongs to 
the vector space spanned by elements $f_{w}$ with $w\in \Sigma_r$ 
, and the total residue of the function  $$ 
(e^0-e^{r+1})^{D-(r+1)}\phi_{r+1}(\ell, k_1,k_2,k_3))$$ is 
proportional to either $[ \sum_{w\in \Sigma_r}\epsilon(w) 
w.f_{\Pi}] $ or  $[ \sum_{w\in \Sigma_r} w.f_{\Pi}] $.  
 The calculation of the constant of proportionality is thus 
equivalent to the Morris iterated constant term identity. However, 
we will give here a direct proof.}

\bigskip

 The  recurrence formula above determines entirely the 
constants $$C_{r+1}(\ell, k_1,k_2,k_3).$$ Indeed, assume first  
$k_3>0$. Then, in  the first relation, when $k_3>0$ 
and $r>1$, all constants $(k_1-1)+\frac{k_3}{2}(r-\ell)$ are 
strictly positive, so if $k_1>1$, we can increase $\ell$ to 
$\ell=r$, then using the second relation,  we decrease $k_1$ to 
$k_1-1$.  Now, if $k_1-1> 1$, then using one we can again 
increase 
 $\ell=0$ to $\ell=r$, and using the second decrease to $k_1-2$.  
 In conclusion, we may determine
using alternatively one and two, the constant $C_{r+1}(\ell, 
k_1,k_2,k_3)$ from $C_{r+1}(\ell,1,k_2,k_3).$ If $k_1=1$, $k_3>0$  
and $r>1$,   we can increase $\ell$ up to $(r-1)$, using  the 
first relation. Then using the third, we decrease $r+1$ to $r$.  
In conclusion we determine $C_{r+1}(\ell, k_1,k_2,k_3)$, if 
$k_3>0$ 
 from the value of $C_2(\ell, 
k_1,k_2,k_3)$. But if $r=1$, there is no factor corresponding to 
$k_3$, so that $C_2(\ell,k_1,k_2,k_3)=C_2(\ell, k_1,k_2,0).$ Now, 
check that constants $C_{r+1}(\ell, k_1,k_2,0)$ 
 are uniquely determined  by the recurrence relations above.
 The first one reads, for $1\leq \ell\leq r$, 
 $$(k_1+k_2-2)C_{r+1}(\ell,k_1,k_2,0)=(k_1-1) C_{r+1}(\ell-1,k_1,k_2,0). $$ 
  In the same way, using alternatively one and two,
  we compute $C_{r+1}(\ell, k_1,k_2,0)$
  from the value of $C_{r+1}(\ell, 1,k_2,0)$.
  Relation  above shows that $C_{r+1}(\ell, 1,k_2,0)=0$,
  if   $\ell$ is not equal to $0$.
  It remains to compute $C_{r+1}(0,1,k,0)$.
    By  the symmetry relation,  we can also assume $k=1$.
  We are finally reduced  to  $C_{r+1}(0,1,1,0)$.

\begin{corollary} Assume $r>1$.

\begin{itemize}
  \item    
 If $k_3>0$, or if $k_1+k_2>2$, then for $1\leq \ell\leq r$,  
  $$C_{r+1}(\ell,k_1,k_2,k_3)=\prod_{j=1}^{\ell}
\frac{k_1-1+(r-j)\frac{k_3}{2}}{k_1+k_2-2+(2r-j-1)\frac{k_3}{2}} 
C_{r+1}(0,k_1,k_2,k_3).$$

  \item 
  
  If $k_1+k_2\geq 2$,
   $$C_{r+1}(0,k_1,k_2,k_3)= r!\prod_{j=0}^{r-1}\frac{\Gamma(1+\frac{k_3}{2})\Gamma(k_1+k_2-1+(r+j-1)\frac{k_3}{2})} 
{\Gamma(1+(j+1)\frac{k_3}{2})\Gamma(k_1+j\frac{k_3}{2})\Gamma(k_2+j\frac{k_3}{2})} 
.$$ 
 \end{itemize}
\end{corollary}

\begin{corollary}\label{catalan}

We have 
\begin{itemize}
  \item $$C_{r+1}(\ell,1,1,1)=\prod_{j=1}^{\ell}
\frac{(r-j)}{(2r-j-1)}C_{r+1}(0,1,1,1).$$ 

  \item $$C_{r+1}(0,1,1,1)= r!\prod_{i=1}^{r-1} C_i,$$  
where $C_i=\frac{(2i)!}{i! (i+1)!}$ is the $i$-th Catalan number.

 \item 
$$ 
C_{r+1}(0,k,1,1)=r!\prod_{i=k-1}^{r+k-3}\frac{1}{2i+1}\binom{r+k+i-1}{2i}.$$

\item 

$$ C_{r+1}(\ell, k,1,1)= \prod_{j=1}^{\ell} 
\frac{2(k-1)+(r-j)}{2(k-1)+(2r-j-1)}C_{r+1}(0, k,1,1).$$

\end{itemize}

\end{corollary}

The second corollary is of course a consequence of the first, 
using several  times  the duplication formula for the Gamma 
function, but it is somewhat easier to use directly the recurrence 
formulas in $k_1,k_2,k_3$, with $k_3=1$ as  the value of $k_3$ 
remains constant and equal to $1$, through the recurrence.

THEOREM \ref{morris} is then a corollary of THEOREM \ref{ct} and  
COROLLARY \ref{catalan}.

Remark that  $$ 
\frac{1}{r!}C_{r+1}(0,1,1,1)=\frac{1}{(r-1)!}C_{r}(0,2,1,1)= 
\prod_{i=1}^{r-1} 
C_i=\prod_{i=1}^{r-2}\frac{1}{2i+1}\binom{r+i}{2i}.$$

Let us first verify the corollaries , assuming THEOREM \ref{ct}.

 To verify the first corollary, we verify the recurrence 
relations, the first being obvious, we check the second:
$$C_{r+1}(r,k_1,k_2,k_3)= C_{r+1}(0,k_1-1,k_2,k_3).$$

We write  $$C_{r+1}(r,k_1,k_2,k_3)$$ $$= r!\prod_{j=1}^{r} 
\frac{k_1-1+(r-j)\frac{k_3}{2}}{k_1+k_2-2+(2r-j-1)\frac{k_3}{2}}$$ 
$$\times\prod_{j=0}^{r-1}\frac{\Gamma(1+\frac{k_3}{2})\Gamma(k_1+k_2-1+(r+j-1)\frac{k_3}{2})} 
{\Gamma(1+(j+1)\frac{k_3}{2})\Gamma(k_1+j\frac{k_3}{2})\Gamma(k_2+j\frac{k_3}{2})}. 
$$ In the first product, we change $j$ in $(r-j)$, in the second  
we use $\Gamma(z+1)=z\Gamma(z)$ and we obtain:

$$ r! \prod_{j=0}^{r-1} 
\frac{k_1-1+j\frac{k_3}{2}}{k_1+k_2-2+(r+j-1)\frac{k_3}{2}} 
\frac{\Gamma(1+\frac{k_3}{2})(k_1+k_2-2+(r+j-1)\frac{k_3}{2})}{\Gamma(1+(j+1)\frac{k_3}{2})(k_1-1+j\frac{k_3}{2})} 
$$ $$\times 
\prod_{j=0}^{r-1}\frac{\Gamma(k_1-1+k_2-1+(r+j-1)\frac{k_3}{2})} 
 {\Gamma(k_1-1+j\frac{k_3}{2})\Gamma(k_2+j\frac{k_3}{2})} $$ 

$$= C_{r+1}(0,k_1-1,k_2,k_3).$$

 We verify the third condition. 

We write

$$C_{r+1}(r-1,1,k_2,k_3)$$ $$= r! \prod_{j=1}^{r-1} 
\frac{(r-j)\frac{k_3}{2}}{k_2-1+(2r-j-1)\frac{k_3}{2}} 
\prod_{j=0}^{r-1}\frac{\Gamma(1+\frac{k_3}{2})\Gamma(k_2+(r+j-1)\frac{k_3}{2})} 
{\Gamma(1+(j+1)\frac{k_3}{2})\Gamma(1+j\frac{k_3}{2})\Gamma(k_2+j\frac{k_3}{2})} 
$$

$$= r! \prod_{j=1}^{r-1} 
\frac{(r-j)\frac{k_3}{2}}{k_2-1+(2r-j-1)\frac{k_3}{2}} 
\prod_{j=0}^{r-1}\frac{\Gamma(1+\frac{k_3}{2})} 
{\Gamma(1+(j+1)\frac{k_3}{2}) \Gamma(k_2+j\frac{k_3}{2})}$$ 
$$\times 
\Gamma(k_2+(r-1)\frac{k_3}{2})\prod_{j=1}^{r-1}\frac{\Gamma(k_2+(r+j-1)\frac{k_3}{2})} 
{\Gamma(1+j\frac{k_3}{2})} $$

In the first product, we change $j$ in $(r-j)$, in the last we use 
$\Gamma(z+1)=z\Gamma(z)$,  and we obtain  after simplification 
that 
 $C_{r+1}(r-1,1,k_2,k_3)$ is equal to $$  r!
\prod_{j=0}^{r-1}\frac{\Gamma(1+\frac{k_3}{2})} 
{\Gamma(1+(j+1)\frac{k_3}{2})\Gamma(k_2+j\frac{k_3}{2})} 
\Gamma(k_2+(r-1)\frac{k_3}{2})$$ $$\times 
\prod_{j=1}^{r-1}\frac{\Gamma(k_2-1+(r+j-1)\frac{k_3}{2})} 
{\Gamma(j\frac{k_3}{2})} $$

 $$= r! \prod_{j=0}^{r-2}\frac{\Gamma(1+\frac{k_3}{2})} 
{\Gamma(1+(j+1)\frac{k_3}{2})\Gamma(k_2+j\frac{k_3}{2})}
\frac{\Gamma(1+\frac{k_3}{2})}{\Gamma(1+r\frac{k_3}{2})} 
\prod_{j=1}^{r-1}\frac{\Gamma(k_2-1+(r+j-1)\frac{k_3}{2})} 
{\Gamma(j\frac{k_3}{2})} $$

$$= (r-1)! \prod_{j=0}^{r-2}\frac{\Gamma(1+\frac{k_3}{2})} 
{\Gamma(1+(j+1)\frac{k_3}{2})\Gamma(k_2+j\frac{k_3}{2})} 
\frac{\Gamma(\frac{k_3}{2})}{\Gamma(r\frac{k_3}{2})} 
\prod_{j=1}^{r-1}\frac{\Gamma(k_2-1+(r+j-1)\frac{k_3}{2})} 
{\Gamma(j\frac{k_3}{2})} $$

while

$$C_{r}(0,k_3,k_2,k_3)=  (r-1)! 
\prod_{j=0}^{r-2}\frac{\Gamma(1+\frac{k_3}{2})\Gamma(k_3+k_2-1+(r-1+j-1)\frac{k_3}{2})} 
{\Gamma(1+(j+1)\frac{k_3}{2})\Gamma(k_3+j\frac{k_3}{2})\Gamma(k_2+j\frac{k_3}{2})} 
.$$ 
 
 $$=  (r-1)! \prod_{j=0}^{r-2}\frac{\Gamma(1+\frac{k_3}{2})} 
{\Gamma(1+(j+1)\frac{k_3}{2})\Gamma(k_2+j\frac{k_3}{2})} 
\prod_{j=0}^{r-2}\frac{\Gamma(k_2-1+(r+j)\frac{k_3}{2})} 
{\Gamma((j+2)\frac{k_3}{2})}$$ 

It remains to verify 

$$ \frac{\Gamma(\frac{k_3}{2})}{\Gamma(r\frac{k_3}{2})} 
\prod_{j=1}^{r-1}\frac{\Gamma(k_2-1+(r+j-1)\frac{k_3}{2})} 
{\Gamma(j\frac{k_3}{2})} 
=\prod_{j=0}^{r-2}\frac{\Gamma(k_2-1+(r+j)\frac{k_3}{2})} 
{\Gamma((j+2)\frac{k_3}{2})}.$$ 

which is true. 

The remaining properties are obvious.

We now prove  THEOREM \ref{ct} by induction on $r$. 

\begin{proof}
If $k_1=0$, the remaining roots $(e^i- e^{j})$ occurring in the 
denominator  of $\phi_{r+1}(\ell,0,k_2,k_3)$ are contained in the 
hyperplane $\sum_{i=0}^{r}e^i=0$. So the total residue of 
$\phi_{r+1}(\ell,0,k_2,k_3)$ is $0$.  The same argument 
shows that  $\phi_{r+1}(\ell,k_1,0,k_3)$ is $0$ . 
 
 We thus may assume that $k_1, k_2 >0$.
 We first  show that the function $(e^0-e^{r+1})\phi_{r+1}(\ell-1,k_1,k_2,k_3)$ 
is proportional to the function $\phi_{r+1}(\ell,k_1,k_2,k_3)$ 
modulo $\sum_{i=1}^r \partial_i R_{A_{r+1}^+}$. Thus the  
functions $ (e^0-e^{r+1})^{D-(r+1)}\phi_{r+1}(\ell-1,k_1,k_2,k_3)$ 
and  $(e^0-e^{r+1})^{D-r}\phi_{r+1}(\ell-1,k_1,k_2,k_3)$ will  be 
proportional too, modulo the vector space $\sum_{i=1}^r \partial_i 
R_{A_{r+1}^+}$, their total residues 
 will be proportional and we 
will get the first recursive relations for the constant $C_{r+1}.$

In the following, we use $x_i$ instead of $e^i$, etc... as it is a 
more familiar notation for computing derivatives.

We write $U=(\prod_{j=1}^r(x_j-x_{r+1}))^{-k_1} 
 (\prod_{j=1}^ r(x_0-x_j))^{-k_2} (\prod_{1\leq i< j\leq r}(x_i-x_j))^{-k_3}$
so that $P_{\ell} U=\phi_{r+1}(\ell,k_1,k_2,k_3)$. 

 We compute $$\partial_1 
[(x_0-x_1)(x_1-x_{r+1})(x_2-x_{r+1})\cdots (x_\ell-x_{r+1})U]$$

 This is equal to  

$$-(1-k_2)(x_1-x_{r+1}) (x_2-x_{r+1})... 
(x_\ell-x_{r+1})U$$ $$+(1-k_1)(x_0-x_1) (x_2-x_{r+1})... 
(x_\ell-x_{r+1})U $$ $$-k_3(x_0-x_1) (x_1-x_{r+1})(x_2-x_{r+1})... 
(x_\ell-x_{r+1})\sum_{j=2}^r\frac{1}{x_1-x_j} U.$$ Using 
$(x_0-x_1)=(x_0-x_{r+1})+(x_{r+1}-x_1)$, this is also equal to

$$=(k_1+k_2-2)(x_1-x_{r+1}) (x_2-x_{r+1})... (x_\ell-x_{r+1})U$$ 
$$+(1-k_1)(x_0-x_{r+1}) (x_2-x_{r+1})... (x_\ell-x_{r+1})U $$ 
$$-k_3(x_0-x_1)(x_1-x_{r+1})(x_2-x_{r+1})... 
(x_\ell-x_{r+1})\sum_{j=2}^r\frac{1}{x_1-x_j} U.$$

Assume first that $k_3$ is odd, so that $U$ is  anti-invariant by 
the group $\Sigma_r$.

Let us antisymmetrize  over permutations. We obtain

$$\sum_{w\in \Sigma_r}\epsilon(w) w\cdot (\partial_1 
((x_0-x_1)(x_1-x_{r+1}) (x_2-x_{r+1})... (x_\ell-x_{r+1})U))= $$ 

$$(k_1+k_2-2)P_{\ell}U+(1-k_1)(x_0-x_{r+1}) P_{\ell-1}\,U $$ 
$$-k_3\sum_{w\in \Sigma_r}w\cdot\left((x_0-x_1)(x_1 
-x_{r+1})(x_2-x_{r+1})... 
(x_\ell-x_{r+1})\sum_{j=2}^r\frac{1}{x_1-x_j}\right) U.$$ 

To compute $\sum_{w\in \Sigma_r}w\cdot((x_0-x_1)(x_1 
-x_{r+1})(x_2-x_{r+1})... (x_\ell-x_{r+1})\frac{1}{x_1-x_j})$, we 
first  sum over the transposition $(j,1)$. 

If $2\leq j\leq \ell$, we use 
$$\frac{(x_0-x_1)(x_1-x_{r+1})(x_j-x_{r+1})}{x_1-x_j}+\frac{(x_0-x_j)(x_j-x_{r+1})
(x_1-x_{r+1})}{x_j-x_1}  $$ $$   =-(x_1-x_{r+1})(x_j-x_{r+1}).$$

If $j>\ell$, we use $$\frac{(x_1-x_{r+1})(x_0-x_1)}{x_1-x_j}+ 
\frac{(x_j-x_{r+1})(x_0-x_j)}{x_j-x_1}$$ 
$$=(x_0-x_1)+(x_{r+1}-x_j)=(x_0-x_{r+1})+(x_{r+1}-x_1)+(x_{r+1}-x_j).$$ 
We obtain  that $$2\sum_{w\in \Sigma_r}w\cdot((x_0-x_1)(x_1 
-x_{r+1})(x_2-x_{r+1})... 
(x_\ell-x_{r+1})\sum_{j=2}^r\frac{1}{x_1-x_j})$$ is equal to 
$$(-(\ell-1)- 2(r-\ell)) P_\ell+ 
(x_0-x_{r+1})(r-\ell)P_{\ell-1}.$$

 Thus finally, we obtain 
 
 $$\sum_{w\in \Sigma_r}\epsilon(w) w\cdot (\partial_1 ((x_0-x_1)(x_1-x_{r+1}) 
(x_2-x_{r+1})... (x_\ell-x_{r+1})\phi_{r+1}(0,k_1,k_2,k_3))= $$ 

$$(k_1+k_2-2+\frac{k_3}{2}(2r-\ell-1))\phi_{r+1}(\ell,k_1,k_2,k_3)$$ 
$$+ (x_0-x_{r+1})(-k_1+1-\frac{k_3}{2}(r-\ell))\phi_{r+1}(\ell-1, 
k_1,k_2,k_3).$$

If $k_3$ is even, we also obtain 
 $$\sum_{w\in \Sigma_r} w\cdot (\partial_1 ((x_0-x_1)(x_1-x_{r+1}) 
(x_2-x_{r+1})... (x_\ell-x_{r+1}))\phi_{r+1}(0,k_1,k_2,k_3))= $$ 
 
$$(k_1+k_2-2+\frac{k_3}{2}(2r-\ell-1))\phi_{r+1}(\ell,k_1,k_2,k_3)$$ 
$$+(x_0-x_{r+1})(-k_1+1-\frac{k_3}{2}(r-\ell))\phi_{r+1}(\ell-1,k_1,k_2,k_3).$$

Thus we see that $$(x_0-x_{r+1})\phi_{r+1}(\ell-1,k_1,k_2,k_3) 
\text{ is proportional to} \  \phi_{r+1}(\ell,k_1,k_2,k_3),$$ 
modulo derivatives with respect to $x_1$, $x_2$,..., $x_r$. In 
particular the total residue of the function 
$(x_0-x_{r+1})^{D-(r+1)}\phi_{r+1}(\ell, k_1,k_2,k_3)$ is 
proportional to the total residue of the function  
$(x_0-x_{r+1})^{D-(r+1)+1}\phi_{r+1}(\ell-1,k_1,k_2,k_3)$.  This 
proves the first property.

 We proceed to the proof of the second property.
We return to the notation $x_i=e^i$.

To avoid confusion in the following argument we will write 
explicitly the dependence by the parameters of $D$, that is we 
will write, whenever necessary,
$$D=D_{A_{r+1}}(\ell,k_1,k_2,k_3)=(k_1+k_2)r+ 
k_3\frac{r(r-1)}{2}-\ell.$$

We now  compute the total residue of 
$$(e^0-e^{r+1})^{D-(r+1)}\phi_{r+1}(r-1, 1,k_2,k_3)$$  
 with $D=D_{A_{r+1}}(r-1,1,k_2,k_3)=(1+k_2)r+k_3\frac{r(r-1)}{2}-(r-1)$ .

 We have: $$r!(e^0-e^{r+1})^{D-(r+1)}\phi_{r+1}(r-1, 
1,k_2,k_3)$$ 
$$=(e^0-e^{r+1})^{D-(r+1)}\frac{P_{r-1}}{(e^1-e^{r+1}) 
(e^2-e^{r+1})... (e^r-e^{r+1})}\phi_{r+1}(0,0,k_2,k_3).$$

 with  

 $$\frac{P_{r-1}}{(e^1-e^{r+1}) (e^2-e^{r+1})... (e^r-e^{r+1})} $$
$$ \frac{ \sum_{w\in \Sigma_r}w\cdot[(e^1-e^{r+1}) 
(e^2-e^{r+1})...(e^{r-1}-e^{r+1})] }{(e^1-e^{r+1}) 
(e^2-e^{r+1})... (e^r-e^{r+1})} 
 $$
$$=(r-1)!\sum_{j=1}^r\frac{1}{(e^j-e^{r+1})}.$$

Consider the subgroup $Cyc_r$ generated by the circular 
permutation of $1,...,r$.   Then  

$$\sum_{j=1}^r\frac{1}{(e^j-e^{r+1})}=\sum_{w\in Cyc_r}w\cdot 
\frac{1}{(e^r-e^{r+1})}.$$

 Assume $k_3$ odd. Thus 
 
 $$\frac{r!}{(r-1)!}\phi_{r+1}(r-1,1,k_2,k_3)=
[ \sum_{w\in Cyc_r} w\cdot  \frac{1 }{(e^r-e^{r+1})}] 
 \phi_{r+1}(0,0,k_2,k_3)$$
 
$$=\sum_{w\in Cyc_r}\epsilon(w)w\cdot [\frac{1 
}{(e^r-e^{r+1})}\phi_{r+1}(0,0,k_2,k_3)]$$ as  
$\phi_{r+1}(0,0,k_2,k_3)$ is anti-invariant under $\Sigma_r$.

 Remark that: $$\phi_{r+1}(0,0,k_2,k_3) =r 
\frac{1}{(e^0-e^{r})^{k_2}}\phi_r(0,k_3,k_2,k_3)$$ 

so that 
$$(e^0-e^{r+1})^{D-(r+1)}\frac{1}{(e^r-e^{r+1})}\phi_{r+1}(0,0,k_2,k_3)$$ 

$$=r  (e^0-e^{r+1})^{D-(r+1)}\frac{1}{(e^r-e^{r+1})} 
\frac{1}{(e^0-e^{r})^{k_2}} \phi_r(0,k_3,k_2,k_3).$$

It follows that we have 

\label{look}

$$ (e^0-e^{r+1})^{D-(r+1)}\phi_{r+1}(r-1,1,k_2,k_3)$$ 
$$=\sum_{w\in Cyc_r}\epsilon(w)w\cdot 
[(e^0-e^{r+1})^{D-(r+1)}\frac{1}{(e^r-e^{r+1})}\frac{1}{(e^0-e^{r})^{k_2}} 
\phi_r(0,k_3,k_2,k_3)].$$

 We now use LEMMA \ref{obvious} to compute the total residue of 
the last term in the equality. 

 We write the vector space $E_{r+1}$ as $E_r\oplus 
\R(e^r-e^{r+1})$, and we consider $\Delta'=A_r\cup 
\{(e^r-e^{r+1})\}$ 

 Using the decomposition 
$(e^0-e^{r+1})=(e^r-e^{r+1})+(e^0-e^{r}),$ we  write   

$$(e^0-e^{r+1})^{D-(r+1)}=\sum_{i\geq 0, j\geq 0,
i+j=D-(r+1)}c_{ij}(e^r-e^{r+1})^i(e^0-e^{r})^j.$$ 

Thus 
$$(e^0-e^{r+1})^{D-(r+1)}\frac{1}{(e^r-e^{r+1})}\frac{1}{(e^0-e^{r})^{k_2}}  
\phi_r(0,k_3,k_2,k_3)$$ $$=\sum_{i\geq 0, j\geq 0, 
i+j=D-(r+1)}c_{ij}\frac{(e^r-e^{r+1})^i} {(e^r-e^{r+1})}  
\frac{(e^0-e^r)^j}{(e^0-e^{r})^{k_2}} \phi_r(0,k_3,k_2,k_3)$$ 
 belongs to the vector space  $ R_{\{(e^r-e^{r+1})\}}  \otimes R_{A_r}$ 
 and we can easily compute  the total residue using  LEMMA \ref{obvious}
 of Section \ref{total}, as well as the obvious calculation  for a one dimensional space.
 Precisely
$$Tres_{A_{r+1}}[(e^0-e^{r+1})^{D-(r+1)}\frac{1}{(e^r-e^{r+1})}\frac{1}{(e^0-e^{r})^{k_2}}  
\phi_r(0,k_3,k_2,k_3)]$$

$$=\sum_{i\geq 0, j\geq 0, i+j=D-(r+1)}c_{ij}Tres_{(e^r-e^{r+1})} 
(\frac{(e^r-e^{r+1})^i} {(e^r-e^{r+1})})$$  $$\times 
Tres_{A_r}[\frac{(e^0-e^r)^j}{(e^0-e^{r})^{k_2}} 
\phi_r(0,k_3,k_2,k_3)].$$ Only the term $i=0$ gives a non zero 
residue, so we obtain 
$$Tres_{A_{r+1}}[(e^0-e^{r+1})^{D-(r+1)}\frac{1}{(e^r-e^{r+1})}\frac{1}{(e^0-e^{r})^{k_2}}  
\phi_r(0,k_3,k_2,k_3)]$$ $$=\frac{1} {(e^r-e^{r+1})}  
Tres_{A_r}[(e^0-e^r)^{D-(r+1)-k_2} \phi_r(0,k_3,k_2,k_3)].$$ 

Now 
$$D-(r+1)-k_2=D_{A_{r+1}}(r-1,1,k_2,k_3)-(r+1)-k_2=D_{A_r}(0,k_3,k_2,k_3)-r.$$
So  we obtain  
$$Tres_{A_{r+1}}[(e^0-e^{r+1})^{D-(r+1)}\frac{1}{(e^r-e^{r+1})} 
(e^0-e^{r})^{k_2}\phi_r(0,k_3,k_2,k_3)]$$ 
$$=\frac{1}{(e^r-e^{r+1})} 
Tres_{A_r}[(e^0-e^{r})^{D_{A_r}(0,k_3,k_2,k_3)-r} 
\phi_r(0,k_3,k_2,k_3)].$$

We apply induction hypothesis on $r$. We have  

$$Tres_{A_r}[(e^0-e^{r})^{D_{A_r}(0,k_3,k_2,k_3)-r} 
\phi_r(0,k_3,k_2,k_3)]$$ $$= C_r(0,k_3,k_2,k_3)\sum_{w'\in 
\Sigma_{r-1}}\epsilon(w') 
w'\cdot[\frac{1}{(e^0-e^1)(e^1-e^2)\cdots (e^{r-1}-e^{r})}].$$  

By FORMULA \ref{look}, as the total residue commutes with the 
action of $W$, we obtain:

 $$ Tres_{A_{r+1}} (e^0-e^{r+1})^{D-(r+1)}\phi_{r+1}(r-1,1,k_2,k_3)$$ 
$$=C_{r}(0,k_3,k_2,k_3)\sum_{w\in Cyc_r}\epsilon(w)w\cdot 
[\frac{1}{(e^r-e^{r+1})}] $$ $$\times \sum_{w'\in \Sigma_{r-1}}\epsilon(w') 
w'\cdot[\frac{1}{(e^0-e^1)\cdots (e^{r-1}-e^{r})}].$$

But  $$ \sum_{w\in Cyc_r}\epsilon(w)w\cdot 
[\frac{1}{(e^r-e^{r+1})}] \sum_{w'\in \Sigma_{r-1}}\epsilon(w') 
w'\cdot[\frac{1}{(e^0-e^1)(e^1-e^2)\cdots (e^{r-1}-e^{r})}]$$   

 $$= \sum_{ww'\in 
Cyc_r\Sigma_{r-1}=\Sigma_{r}}\epsilon(ww') ww'\cdot 
[\frac{1}{(e^r-e^{r+1})} \frac{1}{(e^0-e^1)(e^1-e^2)\cdots 
(e^{r-1}-e^{r})}]$$ $$= \sum_{w\in \Sigma_r}\epsilon(w) 
w.f_{\Pi} . $$ 

Thus we obtain the second relation. The case $k_3$ even is 
completely analogous, so the proof  of the first and second 
relation is complete. The symmetry property in $k_1,k_2$ is 
obvious.

 Let us check   $$C_{r+1}(0,1,1,0)=r!.$$

 More precisely, we have  the following exact formula ( 
without projection on $S_\Delta$).

\begin{lemma} 

$$\frac{(e^0-e^{r+1})^{(r-1)}}{\prod_{i=1}^r (e^0-e^i)(e^i-e^r)} 
=\sum_{w\in \Sigma_r} w\cdot f_{\Pi}.$$ 

\end{lemma}

Indeed, by reduction to the same denominator, the right hand side 
can be written as 

$$\frac{P}{\prod_{i=1}^r (e^0-e^i)(e^i-e^r)\prod_{1\leq i< j\leq 
r}(e^i-e^j)}$$ 

>From invariance consideration, $P$ has to be anti-invariant under 
$\Sigma_r$ , so is divisible by $\prod_{1\leq i< j\leq 
r}(e^i-e^j).$ From degree consideration, we obtain the desired 
equality, and $C_{r+1}(0,1,1,0)=r!$.

\end{proof} QED

\bigskip

{\bf Bibliographical remarks.}

For the proofs of Theorem 26 and 27, we followed the indications 
of  Zeilberger [Z]. Indeed, as explained in the remark following 
Theorem \ref{ct}, these formulae are closely related to 
calculations of Selberg integral [Se], and to   Morris identity 
[M]. We follow here closely Aomoto's [A0] proof of the Selberg  
formula .

\section{Chan-Robbins-Yuen polytope}\label{CR}

Let $B_n$ be the polytope of $n\times n$ doubly stochastic 
matrices; that is the set of $n\times n$ matrices with non 
negative entries, and such that the sum of entries in any row or 
column is equal to $1$. The vector space  of $n\times n$ matrices 
is equipped with the lattice of matrices with integral 
coefficients. The intersection of this lattice with the affine 
space spanned by  $B_n$ is an affine lattice, thus determines an 
unity of volume. The volume of $B_n$ is defined with respect to 
this unity. It has been computed only up to $n=8$ ([CY]). In 
[C-R-Y], some conjectures on the volume of some faces of $B_n$ 
were stated. We show that they follow from Jeffrey-Kirwan formula 
for volumes, together with the explicit calculation of the 
preceding  paragraph.

  \begin{definition}
   Define the vector space   $Y_n$ as the vector space  of $(n\times n)$  matrices, with 
 $x_{ij}=0$ , if $j\geq (i+2)$. 
The Chan-Robbins-Yuen polytope is the polytope $CRY_n= B_n\cap 
Y_n$. 
\end{definition}

 Consider the vector  space of $(n+1)\times (n+1)$  lower triangular matrices
  such that first and last coefficients in the diagonal are equal to $0$. 
 Clearly this vector space is isomorphic to 
 $Y_n$, by adding to a matrix $X\in Y_n$  a row above
 identically equal to $0$, and a last column identically equal to $0$.

{\bf  Example $n=5$.} 

 The vector  space  $Y_5$ consists of  $5\times 5$ matrices 
$X=(x_{ij})$ such that 
 
$$X=\begin{pmatrix} 
  x_{11} & x_{12} & 0 & 0& 0 \\  x_{21} & x_{22} & x_{23} & 0& 0\\
  x_{31} & x_{32} & x_{33} & x_{34} & 0 \\
  x_{41} & x_{42} & x_{43} & x_{44} & x_{45} \\
  x_{51} & x_{52} & x_{53} & x_{54} & x_{55} 
\end{pmatrix}$$

The augmented matrix $aug(X)$ is the $6\times 6$   lower 
triangular matrix 

$$aug(X)=\begin{pmatrix} 0&0&0&0&0&0 \\ 
  x_{11} & x_{12} & 0 & 0& 0&0 \\
  x_{21} & x_{22} & x_{23} & 0& 0&0\\
  x_{31} & x_{32} & x_{33} & x_{34} & 0 &0\\
  x_{41} & x_{42} & x_{43} & x_{44} & x_{45}&0 \\
  x_{51} & x_{52} & x_{53} & x_{54} & x_{55}&0 
\end{pmatrix}$$

 It is clear that the transformation  $X\mapsto aug(X)$
 transforms the polytope $CRY_n$ 
to the polytope  $C'_n$  of $(n+1)\times (n+1)$ lower triangular 
matrices  with non negative entries,  first row and last column 
identically equal to $0$, and other rows and columns summing to 
$1$. 

{\bf  Example $n=5$.} The polytope $C'_5$ consists of matrices 

 $$A=\begin{pmatrix} 0&0&0&0&0&0 \\ 
  a_{21} & a_{22} & 0 & 0& 0&0 \\
  a_{31} & a_{32} & a_{33} & 0& 0&0\\
  a_{41} & a_{42} & a_{43} & a_{44} & 0 &0\\
  a_{51} & a_{52} & a_{53} & a_{54} & a_{55}&0 \\
  a_{61} & a_{62} & a_{63} & a_{64} & a_{65}&0 
\end{pmatrix}$$

such that $a_{ij}\geq 0$, and such that

$$ a_{21} +a _{31}+ a_{41} + a_{51}+a_{61}=1,$$ 

 $$a_{21} + a_{22} = 1=a_{22}+a_{32}+ a_{42}+ a_{52}+a_{62},$$
 \hspace{5cm}.......
$$ a_{51} + a_{52} + a_{53} + a_{54} + a_{55}=1= a_{55}+a_{65},$$ 

$$ 
 a_{61} + a_{62} + a_{63} +a_{64} + a_{65}=1.$$

Let $\n$ be the vector space of strictly lower triangular 
$(n+1)\times (n+1)$- matrices, with basis $E_{ij},(j< i)$ with 
unique non zero entry at the $i$-th row and $j$ column. 
 Consider  a $(n+1)$  dimensional vector space 
 with basis $e_i,1\leq i\leq (n+1)$ and dual basis 
 $e^i$. Let 
$A_n=\{(e^k-e^\ell); k\neq \ell\}$ with positive root  system  
$\{(e^k-e^\ell);1\leq k< \ell\leq (n+1)\}$. Let $p$ be the map 
sending the basis $E_{ji}$ of the vector space $\n$ to 
$(e^j-e^i)$, for $1\leq j< i\leq (n+1)$. The span of $A_n$ is the 
$n$-dimensional vector space $E_n$ consisting of elements 
$a=\sum_{i=1}^{n+1} a_i e^i$, with $\sum_{i=1}^{n+1} a_i=0$

\begin{lemma}
The polytope $CRY_n$ is isomorphic to $P_{A_n^+}(e^1-e^{n+1})$. 
\end{lemma}

We give the proof for $n=5$, as it is most transparent on an 
example. By definition the polytope $P_{A_5^+}(e^1-e^{6})$ 
consists of all $(6\times 6)$ lower triangular matrices 
  $$B=\begin{pmatrix} 0&0&0&0&0&0 \\ 
  b_{21} &0 & 0 & 0& 0&0 \\
  b_{31} & b_{32} &0 & 0& 0&0\\
  b_{41} & b_{42} & b_{43} & 0& 0 &0\\
  b_{51} & b_{52} & b_{53} & b_{54} & 0&0 \\
  b_{61} & b_{62} & b_{63} & b_{64} & b_{65}&0 
\end{pmatrix}$$
with non negative coefficients $b_{ij}$ and such that $$\sum 
b_{ij}(e^j-e^i)= e^1-e^6.$$ This gives the  $6$ equations: 

$$ b_{21} +b _{31}+ b_{41} + b_{51}+b_{61}=1,$$ 
 $$-b_{21}+b_{32}+ b_{42}+ b_{52}+b_{62}=0,$$
 $$-b_{31}-b_{32}+b_{43}+ b_{53}+b_{63}=0,$$
\hspace{5cm}....... 

$$- b_{61} - b_{62} - b_{63}- b_{64} - b_{65}=-1.$$

Define $$b'_{22}=1-b_{21},\hspace{5mm} b'_{33}= 
1-(b_{31}+b_{32})$$ $$b'_{44}=1-(b_{41} + b_{42}+b_{43}) $$  $$ 
b'_{55}=1-(b_{51} + b_{52}+b_{53}+b_{54}) .$$ 

These coefficients are positive: indeed the total sum of 
coefficients of the entries of the matrix $B$ in the rectangle 
$R_{ij}$  consisting on elements in a  row  of index greater than 
$i$ and on a  column of index less than $j$ is identically equal 
to $1$, so that any partial sum is less then $1$. For example, 
adding first $3$ inequalities, we obtain that 

$$ b_{41} + b_{42} +b_{43}+$$ 
 $$ b_{51} + b_{52} + b_{53}+$$
  $$  b_{61} + b_{62} + b_{63} $$
is equal to $1$, so that $ b_{41} + b_{42} +b_{43}$ is less than 
$1$.

 Then the matrix  
 $$B'=\begin{pmatrix} 
  b_{21} &b'_{22} & 0 & 0& 0 \\
  b_{31} & b_{32} &b'_{33} & 0& 0\\
  b_{41} & b_{42} & b_{43} & b'_{44} & 0 \\
  b_{51} & b_{52} & b_{53} & b_{54} & b'_{55}\\
  b_{61} & b_{62} & b_{63} & b_{64} & b_{65}
\end{pmatrix}$$
is in the polytope $CRY_5$.

 We now show that  the conjectures of Chan-Robbins-Yuen
 on their polytope $CRY_n$ and related polytopes are consequences of THEOREM \ref{morris}
 and PROPOSITION \ref{division}

\begin{theorem}\label{Z}(Zeilberger)
The relative volume of $CRY_n$ is equal to 
$$\prod_{i=1}^{n-2}\frac{1}{i+1}\binom{2i}{i}.$$ 
\end{theorem}

\begin{proof} 
The polytope $CRY_n$ is isomorphic to the polytope 
$P_{A_n^+}(e^1-e^{n+1})$. The  linear isomorphism  described above 
using the map $aug$ preserves the volume, as it preserves the 
corresponding lattices. We thus compute the relative volume of  
$P_{A_n^+}(e^1-e^{n+1})$.   Consider $$J_n=\frac{ 
(e^1-e^{n+1})^{n(n-1)/2}} 
 { \prod_{1\leq i< j\leq (n+1)}(e^i-e^j)}.$$
 It is an element of $R_{A_n}$ of homogeneous degree $-n$. We 
choose the big chamber $\c^+$ in $E_n$ consisting of elements 
$a=\sum_{i=1}^{n+1}a_i e^i$, with  $a_i>0$, for $i=1, 2,..., n$. 
The element $a=e^1-e^{n+1}$ is in the closure of this big chamber. 
Jeffrey-Kirwan formula is:

 $$\vol_{rel}P_{A_n^+}(e^1-e^{n+1})=\langle\langle \c^+, 
Tres_{A_n} (J_n)\rangle\rangle.$$

The function $J_n$ up to renumbering has been introduced in 
Section 5. This is the function $W_n$. 
 Recall the formula for its total residue given in THEOREM \ref{morris}.

Let $$f_{\Pi}=\frac{1}{(e^1-e^2)(e^2-e^3)\cdots 
(e^n-e^{n+1})}$$ be the particular element of $S_{\Delta}$ 
obtained by taking the inverse of the product of simple roots. 
Then 

$$Tres_{A_n}(J_n)=\prod_{i=1}^{n-2}\frac{1}{i+1}\binom{2i}{i} 
[\sum_{w\in \Sigma[2,...,n]}\epsilon(w) w.f_{\Pi}.]$$

The linear form $\langle\langle \c^+, f\rangle\rangle$ has been 
computed in Section 4. We write 
$$j_n(x)=\frac{x_1^{n(n-1)/2}}{x_1...x_n \prod_{1\leq i< j\leq n 
}(x_i-x_j)}.$$ It is the function $J_n$ restricted at elements 
$\sum_{i=1}^n x_i e_i$ . Then $$\langle\langle \c^+, 
J_n\rangle\rangle=Ires_{x=0}j_n(x).$$

As iterated residues form a dual basis to the basis $w.f_{\Pi}$ 
of $S_\Delta$ (with $w\in \Sigma[1,2,....,n]$), we obtain 
$\langle\langle 
\c^+,J_n\rangle\rangle=\prod_{i=1}^{n-2}\frac{1}{i+1}\binom{2i}{i}$.

\end{proof} QED

Consider the polytope $CRY_{n}(\widehat {(2,2)})$, consisting on 
elements of $CRY_n$ with the entry $(2,2)$ equal to $0$.   It is a 
face of codimension $1$ of $CRY_n$, isomorphic to the polytope 
$P_{A_n^+\,{\rm minus }\,(e^2-e^3)}(e^1-e^{n+1})$.  Taking in 
account the fact that the dimension of $CRY_n$ is $\binom{n}{2}$, 
we obtain from PROPOSITION  \ref{division} the following 
corollary.

\begin{lemma}\label{Con}(Conjecture 4 of [C-R-Y])

We have $$\binom{n}{2}\vol_{rel} CRY_{n}(\widehat {(2,2)})=3 
\vol_{rel} CRY_{n}.$$ 

\end{lemma}

\bigskip

{\bf Bibliographical remarks.}

The Chan-Robbins-Yuen polytope is described in [C-R-Y].
 Its volume
is  computed by  Zeilberger [Z],  using  a formula of Postnikov- 
Stanley [S]  or of Chan-Robbins-Yuen [C-R-Y]  for the volume  of 
$CRY_n$ , as a particular  value of Kostant partition function 
(see Lemma \ref{st} in next Section \ref{eh}). Our proof follows  
directly from  the Jeffrey-Kirwan residue formulae for the volume, 
but we use also  the crucial identity of Theorem \ref{ct}. 
 The conjectured relations between volumes of flow polytopes and the volume of a 
particular face is deduced from Proposition \ref{division}. 
However, the volume of the entire boundary of the polytope $CRY_n$ 
is unknown .

\section{Ehrhart polynomial of a flow polytope.}\label{eh}

We consider  the positive root system of $A_r$ realized in 
$E_r\subset \R^{r+1}$.  We  will use sometimes implicitly 
the map  $$a=(a_1,a_2,..., a_r)\to {\bf a}=a_1 e^1+...+a_{r} 
e^r-(a_1+\cdots + a_r) e^{r+1}$$ to identify $E_r$ and $\R^r$.

Let $M=r(r+1)/2$. Let ${\bf a}\in \Z^{r+1}\cap C(A_r^+)$. Then 
$a_1\geq 0$, $a_1+a_2\geq 0$, ..., $a_1+a_2+\cdots +a_r\geq 0, 
a_1+a_2+\cdots +a_{r+1}=0$. 
  Let 
$k_{A_r^+}(a)$ be the number of  solutions $(u_1,u_2, ..., u_M)$, 
in non negative integers  $u_m$, of the equation $$\sum_{m=1}^M 
u_m \alpha^m={\bf a}.$$ Here $\alpha^m$ runs though the $r(r+1)/2$ 
positive roots  $(e^i-e^j)$ of $A_r$. The function $k_{A_r^+}(a)$ 
is the number of integral points in the polytope $P_{A_r^+}({\bf 
a})$ and is called the {\bf Kostant partition function} (for 
$A_r$).

 Consider the permutation $w_0: [1,2,...,r, r+1]\mapsto 
[r+1,r,..., 2,1]$. Then $-w_0$ preserves  $A_r^+$. A solution   
$(u_1,u_2, ..., u_M)$, in non negative integers  $u_m$, of the 
equation $$\sum_{m=1}^M u_m \alpha^m={\bf a}$$ gives a solution of 
the equation 
 $$\sum_{m=1}^M u_m (-w_0\cdot \alpha^m)={-w_0\cdot \bf a}.$$
 So  for 
$a_1+a_2+...+a_{r+1}=0$, we have $$k_{A_{r}^+}(a_1e^1+ a_2 
e^2+\cdots +a_r e^r+a_{r+1}e^{r+1})$$  $$ =k_{A_{r}^+}(-a_{r+1}e^1- 
a_{r} e^2-\cdots -a_2 e^r-a_{1}e^{r+1}).$$ 

 Almost by definition, the function $k_{A_r^+}$ is given by 
an iterated constant term.

 \begin{lemma}\label{res}
 Let ${\bf a}= a_1e^1+a_2e^2+...+a_re^{r}+a_{r+1}e^{r+1}$,
 with   $a_1+a_2+...+a_r+a_{r+1}=0$.  
The value  $k_{A_r^+}({\bf a})$ is given by the iterated constant 
term formula: $$k_{A_r^+}({\bf a}) = CT_{z_1=0} 
CT_{z_{2}=0}...CT_{z_{r+1}=0}\left(\frac{z_1^{a_1}z_2^{a_2}... 
z_{r+1}^{a_{r+1}}} 
 {\prod_{1\leq i< j\leq r+1}(1-z_j/z_i)}\right).$$
 \end{lemma}

Here  $\frac{1}{(1-z_{r+1}/z_i)}$ is expanded as a power series in 
$z_{r+1}$. The constant term  at $z_{r+1}=0$ of  
$\frac{z_{r+1}^{a_{r+1}}} 
 {\prod_{1\leq i\leq r+1}(1-z_{r+1}/z_i)}$
 is in the ring  
$R[z_1^{-1}, z_2^{-1},..., z_r^{-1}]$ and we reiterate. 

 Thus, some 
explicit values for Kostant partition function   can be obtained 
from the total residue formula of Section \ref{calculation}, as 
iterated constant term  of a function $f$  coincide with the 
iterated residue   of this function divided by $z_1 z_2...z_r$. 
 
 For example
 \begin{lemma}\label{kk}
For $r\geq 2$, 
$$k_{A_r^+}((1+d),(2+d),...,(r+d))=\prod_{i=d+1}^{r+d-1}\frac{1}{2i+1}\binom{r+d+i+1}{2i},$$   

$$k_{A_r^+}(1,2,...,r)=\prod_{i=1}^r C_i.$$ 
\end{lemma}

 Let $\Phi=\{\alpha^1,\alpha^2,...,\alpha^N\}$ be
 a sequence of elements of $A_r^+$ generating $E_r$.

 Consider the surjective linear map $p:\R^N\to E_r$  such that $p(w^i)=\alpha^i$.
  Then, for $a\in \Z^r$, the polytope 
$P_{\Phi}({\bf a})\subset \R_+^N\cap p^{-1}({\bf a})$ has  
vertices with integral coordinates.

\begin{definition}\label{rere}

For  $a\in \Z^r$,  define $$k_{\Phi}(a)=| P_{\Phi}({\bf a})\cap 
\Z^N|.$$

Thus $k_\Phi(a)$ is the number of  solutions $(u_1,u_2, ..., 
u_N)$, in non negative integers  $u_m$, of the equation 
$$\sum_{m=1}^N u_m \alpha^m={\bf a}.$$

 Define  $k'_\Phi(a)$ as the number of  solutions $(u_1,u_2, ..., 
u_N)$, in strictly positive  integers  $u_m$, of the equation 
$$\sum_{m=1}^N u_m \alpha^m={\bf a}.$$

 \end{definition}

 {\bf Remark.}
 Let $a=\sum_{i=1}^{r+1}a_i e^i$.
If $k'_\Phi(a)$ is non zero, then necessarily $a_1\geq 
\sum_{j=2}^{r+1}m_{1j}$.

We consider the system $-w_0(\Phi)$. Then clearly 
$k_\Phi(-w_0(a))=k_{-w_0(\Phi)}(a)$. 
 
\begin{lemma}\label{itc} Let ${\bf a}= a_1e^1+a_2e^2+...+a_re^{r}+a_{r+1}e^{r+1}$,
 with   $a_1+a_2+...+a_r+a_{r+1}=0$.  
The value  $k_\Phi({\bf a})$ is given by the iterated constant 
term formula: $$k_\Phi({\bf a})=CT_{z_1=0} 
CT_{z_{2}=0}...CT_{z_{r+1}=0}\left(\frac{z_1^{a_1}z_2^{a_2}... 
z_{r+1}^{a_{r+1}}} 
 {\prod_{1\leq i< j\leq r+1}(1-z_j/z_i)^{m_{ij}}}\right).$$
\end{lemma}

We also consider the number $k'(\Phi,\c)(a)$ of solutions in 
strictly positive integers $u_m$  of the equation $$\sum_{m=1}^N 
u_m \alpha^m={\bf a}.$$

Let $\c$ be a big chamber for the system $\Delta^+=\{\Phi\}$. It 
follows from Ehrhart results on polytopes with integral vertices 
that the function $k_\Phi(a)$ is given by a polynomial  formula 
when $a$ varies in $\overline \c \cap \Z^r$. This polynomial is 
called {\bf the Ehrhart polynomial} of the family of polytopes 
$P_{\Phi}({\bf a})$, {\bf on the big chamber $\c$}.

For $\alpha\in \Phi$,  the rational function on $E_r^*$ 
$$\frac{\alpha}{1-e^{-\alpha}}=1+\frac{1}{2}\alpha 
+\frac{1}{12}\alpha^2+...$$ is analytic at the origin, and can be 
expanded as a Taylor series. Thus 
$$\frac{1}{1-e^{-\alpha}}=\frac{1}{\alpha}(\frac{\alpha}{1-e^{-\alpha}})$$ 
as an element of $\hat R_{\Delta}$, and for any ${\bf a}\in E_r$, 
the function $\frac{e^{\bf a}}{\prod_{i=1}^N (1-e^{-\alpha^i})}$ 
is an element of $\hat R_{\Delta}$. 

Define 
 
 $$K_\Phi({\bf a})= Tres_{\Delta} \left(\frac{e^{{\bf a}}}{\prod_{\alpha\in 
\Phi}(1-e^{-\alpha})}\right).$$

 As  $\frac{1}{\prod_{i=1}^N 
(1-e^{-\alpha^i})}=\frac{1}{\prod_{i=1}^N \alpha^i}+ 
\sum_{k>-N}f_k,$ 
  we have  $$Tres_ {\Delta}(\frac{e^{{\bf a}}}{\prod_{i=1}^N 
(1-e^{-\alpha_i})})$$ $$=\frac{1}{(N-r)!} Tres_{\Delta}( 
\frac{{\bf a}^{N-r}}{\prod_{i=1}^N \alpha^i})+  
\sum_{q<(N-r)}\frac{1}{q!} Tres_{\Delta} ({\bf a}^qf_{-r-q}) .$$ 
Thus the  function $$a\mapsto Tres_ {\Delta}(\frac{e^{{\bf 
a}}}{\prod_{i=1}^N 1-e^{-\alpha^i}})$$ is a polynomial function of 
degree $(N-r)$  with value in $S_\Delta$. 
   Recall DEFINITION \ref{J} of the 
function $J_\Phi$. We see that the homogeneous component of degree 
$(N-r)$ of $K_\Phi(a)$ is   
 $$J_\Phi(a)=Tres_{\Delta}(\frac{e^{{\bf a}}}{\prod_{i=1}^N\alpha^i})$$
 $$
 =\sum_{|{\bf i}|=N-r}f_\Phi({\bf i}) \frac{a_1^{i_1}}{i_1!}\frac{a^{i_2}}{i_2!}\cdots 
 \frac{a_{r-1}^{i_{r-1}}}{i_{r-1}!} \frac{a_r^{i_{r}}}{i_{r}!}.$$

 Khovanski-Pukhlikhov formula  [KP] says that the locally polynomial function 
 $k_\Phi(a)$ is deduced (as for any partition family of 
polytopes, see [B-V 1, Section 3.5]) from the locally polynomial 
function $\vol P_\Phi(a)$ via a Todd operator. Using 
Jeffrey-Kirwan residue formula for the volume, we obtain (see [B-V 
1 , Section 3.5]) the following residue formula for the Ehrhart 
polynomials  $k_\Phi(a)$ and $k'_\Phi(a)$ . 

\begin{theorem}\label{mrt}(Multidimensional residue theorem)

Let ${\bf a}\in C(\Delta^+)$. Let $\c\subset C(\Delta^+)$ be a big 
chamber such that ${\bf a}\in \overline{\c}$, then $$k_\Phi(a)= 
\langle\langle \c, K_\Phi(a)\rangle\rangle.$$

 If $a\in \c$, then $$k'_{\Phi}(a)=(-1)^{|\Phi|}
\langle\langle \c, K_\Phi(-a)\rangle\rangle.$$

\end{theorem}

For example, in case   $\Delta=A_1$ and $\Phi$  with 
multiplicity $m$,   the 
 preceding theorem asserts  that for $a\geq 0$ 
$$k_\Phi(a e^1-ae^2)=CT_{z=0}\frac{z^{-a}}{(1-z)^m}=Res_{x=0} 
\frac{e^{ax}}{(1-e^{-x})^m}.$$

 The first  expression   is the residue of the $1$-form 
$\frac{z^{-a}}{(1-z)^m}\frac{dz}{z}$ at $z=0$. The second is the 
residue at $z=1$ of this $1$-form (we choose $z=e^{-x}$ as local 
coordinate, near $x=0$ ). As   $a+m>0$, the rational function 
$\frac{z^{-a}}{(1-z)^m}\frac{dz}{z}$ has no poles at $z=\infty$ 
and we obtain the result from the one dimensional residue theorem.

\begin{definition}
We denote by $k(\Phi,\c)(a)$ the polynomial such that 
$k(\Phi,\c)(a)=k_\Phi(a)$ when $a\in \overline{\c}$ 
\end{definition}

  Thus $$k(\Phi,\c)(a)= \langle\langle \c, 
K_\Phi(a)\rangle\rangle.$$ This is  a polynomial of degree $N-r$ 
and the homogeneous component of degree $N-r$ of $k(\Phi,\c)$ is 
the polynomial $$v(\Phi,\c)(a) = \langle\langle \c, 
J_\Phi(a)\rangle\rangle.$$

As in LEMMA \ref{www}, the knowledge of $k(\Phi,\c^+)(a)$ 
determines (in principle) the other functions $k(\Phi,\c)$.

\begin{lemma}\label{wwww}
Let $\c$ be a big chamber. Let  $\Sigma_r(\c)$ be the set of 
elements $w\in \Sigma_r$ such that $\c\subset C_w^+$. Then 
$$k(\Phi,\c)(a)=\sum_{w\in \Sigma_r(\c)}(-1)^{n(w)} 
k(\Phi,\c^+)(w^{-1} a).$$ 
 
 \end{lemma}

We write $$v(\Phi,\c)(a_1,..., a_r) = \sum _{i_1+i_2+\cdots 
+i_r=N-r}v(\Phi,\c,{\bf i}) 
\frac{a_1^{i_1}}{i_1!}\frac{a_2^{i_2}}{i_2!}\cdots 
 \frac{a_{r-1}^{i_{r-1}}}{i_{r-1}!} \frac{a_r^{i_{r}}}{i_{r}!}.$$

 \begin{definition}\label{ts}
Define $$t_j^{\Phi}=\sum_{k=j+1}^{r+1}m_{jk}-1,$$ 
$$s_j^{\Phi}=1-\sum_{k=1}^{j-1}m_{kj}.$$ 

\end{definition}

We consider the nice chamber $\c^+$. Then we   can express the 
coefficients   of the function $v(\Phi,\c)(a)$  as values at very 
particular points depending on $\bf{i}$ of the function 
$k_{\Phi'}$, where $\Phi'$ is the system where we have deleted all 
roots $e^i-e^{r+1}$.

The system $\Phi'$ is a sequence of vectors in the positive root 
system of $A_{r-1}$, which spans the vector space $E_{r-1}=\{a= 
\sum_{i=1}^r a_i e^i\,,{\rm with}\, \sum_{i=1}^r a_i=0\}$. Remark 
that if $i_1+i_2+...+i_r=N-r$, then the vector $(i_1-t_1) 
e^1+(i_2-t_2) e^2+ ...+(i_r-t_r) e^r$ is in $E_{r-1}$.

 \begin{proposition}(Postnikov-Stanley)\label{ppp}
 
 For the big chamber $\c^+$ of $C(A_r^+)$
 we have 
 $$v(\Phi,\c^+,{\bf i})=k_{\Phi'}((i_1-t_1) e^1+(i_2-t_2) e^2+ ...+(i_r-t_r) e^r).$$
\end{proposition}

\begin{proof}
We use the  iterated residue formula for $f \mapsto \langle\langle 
\c^+,f\rangle\rangle$ and the iterated constant term formula 
(LEMMA \ref{itc}) for $k_{\Phi'}.$ Indeed $$v(\Phi,\c^+,{\bf i}) 
$$ $$= Res_{x_1=0} Res_{x_2=0}\cdots Res_{x_r=0} 
\left(\frac{x_1^{i_1}x_2^{i_2}\cdots 
 x_r ^{i_r}}{x_1^{m_{1,r+1}}...x_r^{m_{r,r+1}}\prod_{1\leq i< j\leq r}(x_i-x_j)^{m_{ij}}}
 \right)$$
 $$=CT_{x_1=0} 
CT_{x_2=0}\cdots CT_{x_r=0} 
\left(\frac{x_1^{i_1+1-m_{1,r+1}}x_2^{i_2+1-m_{2,r+1}}\cdots 
 x_r ^{i_r+1-m_{r,r+1}}}{\prod_{1\leq i< j\leq r}(x_i-x_j)^{m_{ij}}}\right)$$

$$=CT_{x_1=0} CT_{x_2=0}\cdots CT_{x_r=0} 
\left(\frac{x_1^{i_1-t_1}x_2^{i_2-t_2}\cdots 
 x_r ^{i_r-t_r}}{\prod_{1\leq i< j\leq }(1-x_j/x_i)^{m_{ij}}}\right)$$ 

$$=k_{\Phi'}((i_1-t_1)e^1+ (i_2-t_2)e^2+...+ (i_r-t_r) e^r).$$

\end{proof} QED

Consider now the system $\Phi=A_r^+$. We have $t_i^{\Phi}=(r-i)$. 
The system $\Phi'$ is $A_{r-1}^+$. 
 We have thus 
 $$v(A_{r}^+,\c^+,{\bf i})=k_{A_{r-1}^+}((i_1-(r-1))e^1+(i_2+(r-2))e^2+\cdots +
 +i_{r-1}e^{r-1}+i_r e^r)$$

In the right hand side,  if $i_r>0$, the element $ 
(i_1-(r-1))e^1+\cdots 
 +i_{r-1}e^{r-1}+i_r e^r$ cannot be in the  cone generated  by  $(e^i-e^j)$ with 
$1\leq i< j\leq r$, as seen by looking at the component on $e^r$ 
 which  would be negative. Furthermore, if $i_r=0$, we 
see that the coefficients of  roots $e^i-e^r$ in 
 a solution of the equation $$\sum_{\alpha^m\in A_{r-1}^+} 
x_m \alpha^m=(i_1-(r-1))e^1+(i_2+(r-2))e^2+\cdots 
 +i_{r-1}e^{r-1}$$
(with $x_m\geq 0$)  are necessarily $0$.  
 So, consider the system $A_{r-2}^+$ 
realized as $(e^i-e^j), 1\leq i< j \leq (r-1)$. Then we see that 

$$v(A_r^+,\c^+;(i_1,i_2,..., i_{r-1},i_r))=0$$ if $i_r>0$ and  
$$v(A_r^+,\c^+;(i_1,i_2,..., i_{r-1},0))$$ $$= k_{A_{r-2}^+}( 
(i_1-(r-1))e^1+(i_2-(r-2))e^2+ \cdots +(i_{r-1}-1)e^{r-1}).$$

 In particular, we obtain  for monomials,  the relation between the 
volume of the CRY polytope $P_{A_r}(e^1-e^{r+1})$ and Kostant 
Partition function. Indeed, the only ${\bf i}$ to consider is 
$i_1=r(r-1)/2, i_2=0,..., i_r=0$. The point $(i_1-(r-1))e^1 -(r-2) 
e^2+...-2 e^{r-2}-e^{r-1}$ is flipped to the point 
$e^1+2e^2+...+(r-2)e^{r-2}+(r-1)-i_1) e^{r-1}$ under the 
transformation $-w_0$ of the system $A_{r-2}^+$. We then obtain:

\begin{corollary}\label{st}(Postnikov-Stanley)

$$\vol P_{A_r^+}(e^1-e^{r+1})=k_{A_{r-2}^+}(1,2,3,4,..., (r-2)).$$ 

\end{corollary}

Similarly, for any sequence $\Phi$, the  relative volume of the  
flow polytope $P_\Phi(e^1-e^{r+1})$ is an integer, given by the 
Kostant partition function $k_{\Phi'}$ at a particular point. 
Combinatorists are happy of this result, only if they can explain 
this by giving an explicit simplicial decomposition of the 
corresponding  flow polytope.

\bigskip

 {\bf Bibliographical remarks.}
 
 Lemma \ref{kk}, as well as many other results on Kostant Partition function is 
 in [Ki]. Theorem \ref{mrt} follows from [K-P] and [J-K], at least for generic values. 
 A  generalization to any rational polytope is given in [B-V 1].
 Proposition \ref{ppp}  and corollary \ref{st}
 are due to Postnikov- Stanley [ P-S].

 \section{Change of variables for  the total residue}\label{changeof}

We return to the notations of Section 2. 

 Let $F:V_\C\to V_\C$ be an analytic map, such that $F(0)=0$ and 
preserving the open set $U_\Delta=V_\C-\CH_{\C}$ of $V_\C$. If 
$f\in \hat R_\Delta$, the function $(F^*f)(x)=f(F(x))$ is  again 
in $\hat R_\Delta$.  Let $Jac(F)$ be the Jacobian of the map $F$. 
We assume $Jac(F)(0)\neq 0$. We write $F(x)=L(F)(x)+r(x)$ where 
$L(F)$ is a linear invertible map, and $r(x)$ vanishes at $0$ at 
order $2$.  Thus $L(F)$ permutes the hyperplanes $\{\alpha=0 \}$. 
If $f$ is in $S_\Delta$, the function $ L(F)^*(f)$ is again in 
$S_\Delta$.

 \begin{theorem}\label{change}
 For any $f$ in 
 $\hat R_\Delta$, we have the equality in $S_\Delta$:
  $$Jac(F)(0) L(F)^* (Tres_\Delta(f))=Tres_\Delta(Jac(F)(F^*f)).$$
  \end{theorem}

\begin{proof}
 Let $f\in \hat R_\Delta$.
Then $(f-Tres_\Delta (f))dx$  is the differential of some 
$(r-1)$-form $\sum_{k=1}^r f_k dx_1\wedge 
dx_2\wedge\widehat{dx_k}\wedge \cdots \wedge dx_r$, with $f_k\in 
\hat R_{\Delta}$. The vector space $\sum_{k=1}^r f_k dx_1\wedge 
dx_2\wedge\widehat{dx_k}\wedge \cdots \wedge dx_n$ with $f_k\in 
\hat R_{\Delta}$ is stable by the action of $F^*$ on differential 
forms.  As $F^*$ commutes with $d$, $(F^*f) F^*dx- F^*(Tres_\Delta 
(f)) F^*dx$ is   the differential of some $(r-1)$-form 
$\sum_{k=1}^r \Phi_k dx_1\wedge dx_2\wedge\widehat{dx_k}\wedge 
\cdots \wedge dx_r$, with $\Phi_k\in \hat R_\Delta$. 
  We have 
$(F^*f) F^*dx =Jac(F) (F^*f)dx.$ 

 We now analyze the 
differential form $(F^*f)(F^* dx)= (F^*f) Jac(F)dx$ with  $f\in 
S_\Delta$. 
 The function $f$ is the inverse of the 
product of $r$ linear forms $\alpha^k$. As $F$ preserves 
$U_\Delta$, we have $\alpha^k(F(x))=\beta^{k}(x)g^k(x)$ where 
$\beta^k(x)$ is another linear form in the system $\Delta$ and 
$g^k(x)$ is holomorphic at $0$. Thus 
$(F^*f)(x)=\frac{1}{\prod_{k}\beta^{k}(x) }\prod_k (g^k(x))^{-1}$ 
is again in $\hat R_\Delta$. 

Furthermore, we see that $(F^*f)-(L(F)^*f)$ is an element of 
$R_\Delta$ of degree $\rangle-r$. Thus $(F^*f)J(F)-(L(F)^*f) 
J(F)(0)$ is an element of $\hat R_\Delta$ of degree $>-n$.  It 
follows that $(F^*f)F^*dx- L(F)^*(f) J(F)(0)dx$  is  the 
differential of some $(r-1)$-form $\sum_{k=1}^r h_k(x) dx_1\wedge 
dx_2\wedge\hat{dx_k}\wedge \cdots \wedge dx_r$, with $h_k\in \hat 
R_{\Delta}$. 

Adding  these two informations, we see that for any $f\in \hat 
R_\Delta$, 
 $$Jac(F)(F^*f)dx -J(F)(0) (L(F)^*Tres_\Delta(f))dx$$ 
 is the 
differential of some $r-1$-form $\sum_{k=1}^r m_k(x) dx_1\wedge 
dx_2\wedge\hat{dx_k}\wedge \cdots \wedge dx_r$, with $m_k\in \hat 
R_{\Delta}$. This implies the formula of the lemma. 
\end{proof} QED

\section{A nice formula for Kostant restricted partition function\label{nice}}

Let $\Phi$ be a  sequence of $N$ vectors  in $A_r^+\subset E_r$ 
generating $E_r$. In the same spirit that Lidskii formula for 
Kostant partition function ([L]), we will give a closed formula 
for $k_\Phi(a)$ in function of $\vol P_\Phi(a)$.  
 In fact, we will express the $S_\Delta$-valued polynomial 
 $K_\Phi(a)$ in function of the 
 $S_\Delta$-valued polynomial $J_\Phi(a)$ of Section \ref{some}.

  Let $$\left(\binom{u}{k}\right)=\frac{u(u+1)(u+2)\cdots (u+(k-1))}{k!}.$$

 \begin{theorem}\label{vw}
Let  $$J_\Phi(a)=\sum _{i_1+i_2+\cdots +i_r=N-r}f_\Phi({\bf i}) 
 \frac{a_1^{i_1}}{i_1!}\frac{a_2^{i_2}}{i_2!}\cdots 
 \frac{a_{r-1}^{i_{r-1}}}{i_{r-1}!} \frac{a_r^{i_{r}}}{i_{r}!}.$$ 
 
 Then 
 $$K_\Phi(a)=\sum_{{|\bf i|}=N-r} 
f_\Phi({\bf i}) \binom{a_1+t_1^{\Phi}}{i_1} 
\binom{a_{2}+t_2^{\Phi}}{i_{2}}\cdots  
\binom{a_{r-1}+t_{r-1}^{\Phi}}{i_{r-1}} 
\binom{a_{r}+t_r^{\Phi}}{i_{r}}.$$  
 
 We have as well 
 
 $$K_\Phi(a)=\sum_{{|\bf i|}=N-r} f_\Phi({\bf i}) 
\left(\binom{a_1+s_1^{\Phi}}{i_1}\right) \cdots  \left(
\binom{a_{r-1}+s_{r-1}^{\Phi}}{i_{r-1}}\right)\left(
\binom{a_{r}+s_r^{\Phi}}{i_{r}}\right).$$ 
 \end{theorem}

\begin{proof}
We realize $A_r$ in $\R^r$ as the system $(e^i-e^j), 1\leq i< 
j\leq r,  e^i, 1\leq i\leq r$.

Then  $$K_{\Phi(a)}$$  $$=Tres_{\Delta} [\frac{e^{a_1 
x_1}e^{a_2 x_2}\cdots e^{a_r 
x_r}}{\prod_{i=1}^r(1-e^{-x_i})^{m_{i,r+1}}\prod_{1\leq i< j\leq r 
}(1-e^{-(x_i-x_j)})^{m_{ij}}}]. $$ 

 We use  LEMMA \ref{change} for  changing  variables in 
residues. 
 Let $x=\sum_{i=1}^r x_i e_i$ in $V$. Define $F(x)=\sum_{i=1}^r 
(e^{x_i}-1)e_i$. This change of variables preserves the vector 
space $R_\Delta.$ Indeed $(e^{x_i}-e^{x_j})$ is divisible by 
$(x_i-x_j)$. The differential of $F$ at the origin is the 
identity. We have $Jac(F)=e^{x_1}e^{x_2}\cdots e^{x_r}.$ 

We write $t_k,s_k$ instead of $t_k^{\Phi},s_k^{\Phi}$. Let $a_1, 
a_2,..., a_r$ be integers and let $$f(a_1,a_2,..., 
a_r)(x)=\frac{(1+x_1)^{a_1+t_1}(1+x_2)^{a_2+t_2}\cdots 
(1+x_r)^{a_r+t_r}}{x_1^{m_{1,r+1}}x_2^{m_{2,r+2}}\cdots 
x_r^{m_{r,r+1}} \prod_{1\leq i< j\leq r}(x_i-x_j)^{m_{ij}}}.$$ 
Then the function $(F^*f) Jac(F)$ is equal to $$ \frac{e^{a_1 
x_1}e^{a_2 x_2}\cdots e^{a_r 
x_r}}{(1-e^{-x_1})^{m_{1,r+1}}(1-e^{-x_2})^{m_{2,r+1}}\cdots 
(1-e^{-x_r})^{m_{r,r+1}}}$$$$\times \frac{1}{\prod_{1\leq i< j\leq r 
}(1-e^{-(x_i-x_j)})^{m_{ij} }}. $$ 

 Thus, we obtain  from LEMMA \ref{change}:

  $$Tres_{\Delta}(\frac{e^{\langle 
a,x\rangle}}{\prod_{\alpha\in \Phi }(1-e^{-\langle \alpha, 
x\rangle})})= Tres_{\Delta} (\frac{(1+x_1)^{a_1+t_1}\cdots 
(1+x_r)^{a_r+t_r}}{\prod_{\alpha\in \Phi} \alpha(x)}).$$ 

  To compute the total residue of the last expression, 
 as the denominator $\prod_{\alpha\in \Phi} \alpha(x)$ 
is  homogeneous of degree $N$, we have to seek for the term in the 
numerator which is homogeneous of degree $N-r$, thus we seek the 
coefficient of  each term of the form $x_1^{i_1}\cdots x_r^{i_r}$ 
with $i_1+i_2+\cdots +i_n=N-r$. We thus obtain the first part of 
Proposition \ref{vw}.

The second is proved the same way, using 
$F(x)=\sum_{i=1}^r(1-e^{-x_i}) e_i$ and the function 
$$\frac{(1+x_1)^{-(a_1+s_1)}\cdots(1+x_r)^{-(a_2+s_r)}}{x_1^{m_{1,r+1}}x_2^{m_{2,r+1}}\cdots 
x_r^{m_{r,r+1}} \prod_{1\leq i< j\leq r}(x_i-x_j)^{m_{ij}}} $$ 

 \end{proof} QED

  Thus we obtain the following  formula for $k_\Phi(a)$.

  \begin{theorem}\label{ch}
  
  Let $\c$ be a big chamber and $a\in \overline \c$.
  Then 
  $$k_\Phi(a)=\sum_{{|\bf i|}=N-r} 
f_\c({\bf i}) \binom{a_1+t_1^{\Phi}}{i_1} 
\binom{a_{2}+t_2^{\Phi}}{i_{2}}\cdots  
\binom{a_{r-1}+t_{r-1}^{\Phi}}{i_{r-1}} 
\binom{a_{r}+t_r^{\Phi}}{i_{r}}.$$  
 
 We have as well 
 
 $$k_\Phi(a)=\sum_{{|\bf i|}=N-r} f_\c({\bf i}) 
\left(\binom{a_1+s_1^{\Phi}}{i_1}\right) \cdots  \left(
\binom{a_{r-1}+s_{r-1}^{\Phi}}{i_{r-1}}\right)\left(
\binom{a_{r}+s_r^{\Phi}}{i_{r}}\right).$$ \end{theorem}

  Comparing with formula \ref{jj} for the volume, we see that 
  the function $k_\Phi(a)$ is immediately deduced from the 
  polynomial function $\vol P_\Phi(a)$ (i.e; its highest degree component)
   by replacing the monomial 
  $\frac{a_k^{i_k}}{i_k!}$
  by the function  $\binom{a_{k}+t_k^{\Phi}}{i_{k}}$ (with same leading term).

 \begin{corollary}
 Let $q^{\Phi}=(\sum_{k=2}^{r+1} m_{1k})-1$. 
 The polynomial $k(\Phi,\c^+)$ is divisible by 
 $(a_1+1)(a_1+2)...(a_1+q^\Phi)$.
 \end{corollary}
 \begin{proof}
 Indeed, on the big chamber $\c^+$, we have seen in PROPOSITION \ref{petit}
  that that multiindices ${\bf i}$ such that
   $f_\Phi(\bf{i})$ is non zero are such that $i_1\geq q^{\Phi}$. 
   We have $s_1^{\Phi}=1$, thus the corollary follows.
   \end{proof} QED

{\bf Remark} This follows also from the remark following 
definition \ref{rere}, and Theorem \ref{mrt}.

If $\Phi=A_r^+$, we then obtain the relation for Kostant partition 
$k_\c$  on a big chamber $\c$ .   Define $$f_\c({\bf 
i})=\langle\langle \c, 
Tres_{\Delta}\left(\frac{x_1^{i_1}x_2^{i_2}\cdots 
 x_r ^{i_r}}{x_1 x_2...x_r \prod_{i< j}(x_i-x_j)}\right)\rangle\rangle.$$

  Then we have

\begin{proposition}\label{lid}(Lidskii)

We have:

$$k(A_r^+,\c)(a)$$

$$= \sum_{|{\bf i}|=\binom{r}{2}} f_\c({\bf i}) 
\binom{a_1+r-1}{i_1} \binom{a_{2}+r-2}{i_{2}}\cdots  
\binom{a_{r-1}+1}{i_{r-1}} \binom{a_{r}}{i_{r}}.$$ 

We have as well:

 $$k(A_r^+,\c)(a)$$
 $$= \sum_{|{\bf i}|=\binom{r}{2}} f_\c({\bf i}) 
\left(\binom{a_1+1}{i_1}\right) 
\left(\binom{a_{2}}{i_{2}}\right)\cdots  
\left(\binom{a_{r-1}+3-r}{i_{r-1}}\right) 
\left(\binom{a_{r}+2-r}{i_{r}}\right).$$ 

\end{proposition}

\bigskip

{\bf Bibliographical remarks.}

Proposition \ref{lid} is due to  B.V. Lidskii [L].  Theorem 
\ref{ch} is due to Stanley, who suggested  to look for a proof via 
residues, as given here.

\section{ Volumes and Ehrhart polynomials of the Stanley-Pitman polytope}\label{pitman}

 Section \ref{nice} reduces the computation of
 the Ehrhart polynomial to the computation of mixed volumes.

Here is an example. Let $$\Pi_r(a)=\{y\in \R^r; y_i\geq 0, 
y_1+\cdots +y_i\leq a_1+a_2+\cdots +a_i\}$$ 
 with $a_i\geq 0$. 

Let $\Phi$ be the following sequence of $2r$ elements 
$\{\beta^i,\gamma^j\}, 1\leq i\leq r, 1\leq j\leq r\}$ of $A_r^+$ 

$$\Phi= \{\beta^1=(e^1-e^{r+1}),....,\beta^r=(e^r-e^{r+1})\} $$ 
$$\cup\{ \gamma^1=(e^1-e^2), \gamma^2=(e^2-e^3), \cdots 
\gamma^{r-1}=(e^{r-1}-e^r), \gamma^r=(e^{r}-e^{r+1})\}$$ 

The multiplicity of $(e^r-e^{r+1})$ in $\Phi$ is $2$. We write 
$t_i=t_i^{\Phi}$ . Thus 

$t_1=t_2=...=t_{r-1}=t_r=1$, while $s_1=1$, $s_2=s_3=...= s_r=0$. 

 Let 
$\n_\Phi=\R^{2r}$ with basis $v_{\beta^k}, w_{\gamma^j}$ . Let 
$C_{2r}^+$ be the standard cone $$\oplus_{k=1}^r 
y_{i,r+1}v_{\beta^k} \oplus_{k=1}^r z_{k,k+1} w_{\gamma^k}$$ with 
$y_{i,r+1}\geq 0$, $z_{k,k+1}\geq 0$  
 
 Let  $p: \R^{2r}\to E_r$ be the map sending 
 $v_\beta$ to $\beta$ and 
 $w_\gamma$ to $\gamma$.
 Then for $a=a_1 e^1+\cdots +a_r e^r$, with $a_i\geq 0$
the polytope $\Pi_r(a)$ is isomorphic 
 to 
 $P_\Phi(a)= p^{-1}({\bf a})\cap C_{2r}^+$.
Indeed the point $\sum_{k=1}^r y_{k,r+1}v_{\beta^k}+ \sum_{k=1}^r 
z_{k,k+1} w_{\gamma^k}$ with $y_{i,r+1}\geq 0$, $z_{i,i+1}\geq 0$   
is in $P_\Phi(a)$ if and only if

 $$\sum_{k=1}^r y_{k,r+1} (e^k-e^{r+1})+\sum_{i=1}^{r} 
z_{i,i+1} (e^i-e^{i+1})$$ $$=a_1 e^1+ \cdots +a_r 
e^r-(a_1+a_2+...+a_r) e^{r+1}.$$ 
This gives 

\begin{eqnarray*}
\\y_{1,r+1}+z_{1,2}=a_1,\\
y_{2,r+1}+z_{2,3}=a_2+z_{1,2}\\
y_{3,r+1}+z_{3,4}=a_3+z_{1,2}+z_{2,3}
\end{eqnarray*}
$$......$$
  
 so that  
 
 $$y_{1,r+1}\leq a_1,$$
$$y_{1,r+1}+y_{2,r+1}\leq a_1+a_2,$$ 
$$y_{1,r+1}+y_{2,r+1}+y_{3,r+1}\leq a_1+a_2+a_3$$  $$\cdots 
$$  and  the point $(y_{1,r+1},y_{2,r+1},..., y_{r,r+1})\in 
\Pi_r(a)$. 

We compute the volume 

$$\vol \Pi_r(a)=\sum_{i_1+i_2+\cdots +i_r=r} f_{\Phi}({\bf i}) 
\frac{a_1^{i_1}}{i_1!}\frac{a_2^{i_2}}{i_2!}\cdots 
 \frac{a_{r-1}^{i_{r-1}}}{i_{r-1}!} \frac{a_r^{i_{r}}}{i_{r}!}$$
As  $a\in \c^+$, we have, by  LEMMA  \ref{ppp}, $$f_{\Phi}({\bf 
i})=k_{\Phi'}((i_1-1)e^1+(i_2-1) e^2+....+(i_r-1)e^r).$$

The system $\Phi'$ is the set of simple roots 
$(e^1-e^2),...,(e^{r-1}-e^r)$.  They are linearly independent and 
generate a simplicial cone $C(\Phi')$. Thus the function 
$k_{\Phi'}$ is identically $1$ on the cone $C(\Phi')\cap \Z^r$. 
 Thus we obtain $f(i_1, i_2, i_3,..., i_r)=0$ or $1$. It is $1$ if 
and only $(i_1-1)e^1+(i_2-1) e^2+....+(i_r-1)e^r$ is in the cone 
$C(A_{r-1}^+)$. We thus need 

$$(i_1,i_2,i_3,...,i_r)\in K_r$$ where

$$K_{r}=\{(i_1,i_2,i_3,...,i_r), i_1\geq 1, i_1+i_2\geq 2,..., 
i_1+i_2+\cdots +i_{r}=r\}.$$

We obtain the formula for the volume of $\Pi_{r}(a)$ and its 
Ehrhart polynomial, given in  [Pi-S] .

\begin{proposition} (Pitman-Stanley)

$$\vol \Pi_r(a)=\sum_{{\bf i}\in K_r} 
\frac{a_1^{i_1}}{i_1!}\frac{a_2^{i_2}}{i_2!}\cdots 
 \frac{a_{r-1}^{i_{r-1}}}{i_{r-1}!} \frac{a_r^{i_{r}}}{i_{r}!},$$

 $$k_\Phi(a)=\sum_{{\bf i}\in K_r} 
\left(\binom{a_1+1}{i_1}\right) 
\left(\binom{a_{2}}{i_{2}}\right)\cdots  \left(
\binom{a_{r-1}}{i_{r-1}}\right)\left(
\binom{a_{r}}{i_{r}}\right).$$  
 
 \end{proposition}

 The same beautiful occurs for any family of polytopes
associated to a smooth toric variety, giving rather nice formulae 
deduced immediately from the mixed volumes.

\bigskip

{\bf Bibliographical remarks.} Results of this section are due to 
[Pi-S]. 
 
\section{Divisibility  property  of the Kostant partition function}
\label{some}.

Here  we  list some properties of  the polynomial 
$k(A_r^+,\c^+)$. They are similar  to the properties of its 
highest degree term $v(A_r^+,\c^+)$ established in Section 
\ref{volume}. For example, the  following Lemma implies of course 
LEMMA \ref{division}.

\begin{proposition}\label{disionbis}( Schmidt-Bincer)

The function $k(A_r^+,\c^+)(a_1,a_2,..., a_r)$ is independent of 
$a_r$. It is of degree less or equal to $1$ in the variable 
$a_{r-1}$ and is divisible by 
$(a_1+a_2+a_3+...+a_{r-2}+3a_{r-1}+3) $. 
 
 More precisely, we have:
 
$$3k(A_r^+,\c^+)(a_1,a_2,..., a_{r-1},a_r)
=3k(A_r^+,\c^+)(a_1,a_2,..., a_{r-1},0)$$ 
$$=(a_1+a_2+...+a_{r-2}+3a_{r-1}+ 3)k(A_r^+,\c^+)(a_1,a_2,..., 
a_{r-1},0)$$ $$=(a_1+a_2+...+a_{r-2}+3a_{r-1}+3 ) k(A_r^+ \, {\rm 
minus}\, (e^{r-1}-e^r),\c^+)(a_1, a_2,..., a_{r-2},0,0).$$

 \end{proposition}

\begin{proof}
The proof is almost identical to the proof of PROPOSITION  
\ref{division}. Let $$K(a,x)$$ $$= \frac{e^{a_1 
x_1+...+a_{r-1}x_{r-1}+ a_r x_r}}{\prod_{1\leq i< j\leq 
r-1}(1-e^{-(x_i-x_j)})\prod_{1\leq i\leq r-1}(1-e^{-(x_i-x_r)}) 
\prod_{1\leq i\leq r} (1-e^{-x_i})}.$$

We have: 

$$k(A_r^+,\c^+)(a_1,a_2,..., a_r)=Ires_{x=0}K(a,x).$$  

We write  $k(A_r^+,\c^+)=k_r^+$. We first  take the residue of 
$K(a,x)$ in $x_r=0$. We obtain 

$$Res_{x_r=0}K(a,x)=\frac{e^{a_1 x_1+...+a_{r-1} 
x_{r-1}}}{\prod_{1\leq i< j\leq r-1}(1-e^{-(x_i-x_j)})\prod_{1\leq 
i\leq r-1}(1-e^{-x_i})^2}.$$

This shows already that $k_r^+(a_1,a_2,..., a_r)$ is independent 
of $a_r$. We proceed now to take the residue in $x_{r-1}=0$. There 
is a double pole in $x_{r-1}$, so that the dependence in $a_{r-1}$ 
is of degree at most $1$. More precisely, a simple calculation 
shows that 

$$Res_{x_{r-1}=0}\frac{e^{a_1 x_1+...+a_{r-1} 
x_{r-1}}}{\prod_{1\leq i< j\leq r-1}(1-e^{-(x_i-x_j)})\prod_{1\leq 
i\leq r-1}(1-e^{-x_i}){}^2}$$  

  $$=\frac{e^{a_1 x_1+...+a_{r-2} x_{r-2}}}{\prod_{1\leq i< 
j\leq r-2}(1-e^{-(x_i-x_j)})\prod_{1\leq i\leq 
r-2}(1-e^{-x_i})^3}(a_{r-1}+1+\sum_{i=1}^{r-2}\frac{e^{-x_i}}{1-e^{-x_i}})$$

$$=(a_{r-1}+1+\frac{1}{3}(a_1+a_2+...+a_{r-2}))\frac{e^{a_1 
x_1+...+a_{r-2}x_{r-2}}}{\prod_{1\leq i< j\leq 
r-2}(1-e^{-(x_i-x_j)})\prod_{1\leq i\leq r-2}(1-e^{-x_i})^3}$$ 
$$-\frac{1}{3}(\partial_1+\partial_2+\cdots 
+\partial_{r-2})\frac{e^{a_1 x_1+...+a_{r-2}x_{r-2}}}{\prod_{1\leq 
i\leq r-2}(1-e^{- x_i})^3\prod_{1\leq i< j\leq 
r-2}(1-e^{-(x_i-x_j)}) }.$$

As residue vanishes on derivatives, we obtain $$k_r^+(a_1,a_2,..., 
a_{r-1},a_r)$$ 
$$=(a_{r-1}+1+\frac{1}{3}(a_1+a_2+...+a_{r-2}))\times$$ 
$$Res_{x_1=0}...Res_{x_{r-2}=0}\frac{e^{a_1 
x_1+...+a_{r-2}x_{r-2}}}{\prod_{1\leq i< j\leq 
r-2}(1-e^{-(x_i-x_j)})\prod_{1\leq i\leq r-2}(1-e^{-x_i})^3}.$$ 
 On the other hand, 
the residue computation of $$k(A_r^+\,{\rm minus}\,
(e^{r-1}-e^{r}),\c^+)(a_1, a_2,..., a_{r-2},0,0)$$ gives 

  $$k(A_r^+ \, {\rm minus}\, (e^{r-1}-e^{r}),\c^+)(a_1, 
a_2,..., a_{r-2},0,0)$$ 

$$=Res_{x_1=0}...Res_{x_{r-2}=0}\frac{e^{a_1 
x_1+...+a_{r-2}x_{r-2}}}{\prod_{1\leq i< j\leq 
r-2}(1-e^{-(x_i-x_j)})\prod_{1\leq i\leq r-2}(1-e^{-x_i})^3}$$ as 
the  step $x_{r}=0$ as well as the step $x_{r-1}=0$ involves only 
simple poles, and we obtain the divisibility property announced.

 \end{proof} QED
 
 More generally, we have the following lemma, with same proof.
  
 \begin{lemma}
 Let $\Phi$ a sequence of vectors in 
 $A_r^+$ generating $E_r$. Assume $m_{r,r+1}=1$ and  $m_{r-1,r+1}+m_{r-1,r}=2$.
 Furthermore, assume that $$\frac{m_{j,r+1}+m_{j,r}+m_{j,r-1}}{m_{j,r-1}}=c$$
 is independent of $j$ for $1\leq j\leq (r-2)$, then 
 
 $$ k(\Phi,\c^+)(a_1,..., a_{r-1},a_r)=k(\Phi,\c^+)(a_1,..., a_{r-1},0)$$
 $$=
 (\frac{a_1+\cdots +a_{r-2}}{c}+ a_{r-1}+1)v(\Phi\,{\rm minus}\,(e^{r-1}-e^{r}),\c^+)
 (a_1,a_2,..., a_{r-2},0,0).$$
\end{lemma}

\bigskip

{\bf Bibliographical remarks.}

 Proposition \ref{disionbis} is due to 
[S-B].

\section{A ``not so obvious " symmetry property of the volume and Ehrhart polynomial}
\label{sym}

 The full Weyl group $\Sigma_{r+1}$ of the root system 
 $A_r$ acts on $S_{A_r}$ as it 
permutes the elements $(e^i-e^j), 1\leq i<j\leq (r+1)$.  We show 
that the existence of this action implies some constraints on the 
coefficients of the $S_{\Delta}$-valued  polynomial function  
$J_{A_r^+}(a)$ and $K_{A_r^+}(a)$. In particular, the Kostant 
polynomial $k(A_r^+,\c^+)( a)=Ires_{x=0}(K_{A_r^+}(a))$  satisfies 
some symmetries.

\begin{lemma} For any $w\in \Sigma_{r+1}$, we have 
$$J_{A_r^+}(w.a)=\epsilon(w) w\cdot J_{A_r^+}(a).$$ 
\end{lemma}

Let $q$ be a linear form on $S_\Delta$. Then the map $a\mapsto 
\langle q, J_{A_r^+}(a)\rangle$ is a polynomial of degree 
$|A_{r-1}|=r(r-1)/2$ on $E_r$. This way, we obtain  a map from the 
$r!$-dimensional space $S_\Delta^*$ to the space of polynomials of 
degree $|A_{r-1}|=r(r-1)/2$ on $E_r$, commuting with the 
representation of $\Sigma_{r+1}$. 

\bigskip

A subset of $\{1,2,...,N\}$ will be called an interval, if it is 
of the form $\{a, a+1, a+2,...,a+k\}$, for $a$ an integer between 
$1$ and $N$ (and $k$ an integer between $0$ and $ N-a$). Let us 
represent   a permutation  in  $\Sigma_{N}$ as a list of   $N$  
elements .  For example, if we write $w=[Id_{(N-3)}, N,N-2,N-1]$,  
the first $N-3$ indices are fixed by the permutation $w$ , $N-2$ 
is sent to $N$, $N-1$ to $N-2$   and $N$ to $N-1$.    
 \begin{definition}
The subset $B_N$   is the subset of  elements $w$ of $\Sigma_N$  
such that for any $1\leq k\leq N$  the subset  $w^{-1}\{1,2,3,..., 
k\}$ is an interval. 
\end{definition}

\begin{lemma}
The subset $B_N$ is of cardinal $2^{N-1}$. 
\end{lemma}
\begin{proof}
We prove it by induction on $N$. If $v\in B_{N+1}$, the set $ 
v^{-1}\{1,2,3,..., N\}$ is an interval. So it is either equal to 
$\{1,2,3,..., N\}$ or to $\{2,3,..., N,N+1\}$. Thus a 
transformation $w$ of $B_N$ gives rise to two transformations  of 
$B_{N+1}$ namely $[w, N+1]$ and $[N+1, w]$.  
\end{proof} QED

\begin{definition}
For $1\leq i\leq r$, we define $W^{i,r}$ as the subset of 
$\Sigma_r$ consisting of elements $w$ such that: 

\begin{itemize}
  \item $ w(k)=k$,  $\, 1\leq k\leq i-1.$ 
  \item  For  $0\leq s\leq (r-i)$,  the set  
 $w^{-1}\{i,r, r-1,r-2, r-3,..., r-s\}$ is an interval {\rm (}of $\{i,..., r\}${\rm )}.
 \end{itemize}
 \end{definition}

Clearly the subset $W^{i,r}$ is isomophic to $B_{r-(i-1)}$ by 
relabelling indices. Thus the set $W^{i,r}$ is of cardinality 
$2^{r-i}$. A list of $r$ elements  representing an element of 
$W^{i,r}$  is constructed as follows: we start writing $i$ in the 
middle of  a horizontal line, then we write  $r$ either  at the 
immediate left of $i$, or  at the immediate  right. After, we 
write  $(r-1)$ 
  either  at the immediate left of  the unordered set $\{i,r\}$, or  at the immediate 
   right, and we go on until we have written all elements of
    the set  $\{i,r.., i+1\}= \{i, i+1,..., r\}$
   on this  line. At this step,  we  finally write the  ordered set  
   $\{1, 2,..., i-1\}$  at the immediate left of the unordered set 
   $\{i, i+1,..., r\}$.

 {\bf Example.}
 Let us list elements of $W^{i,r}\subset \Sigma_r$ for $i=r, r-1, r-2$.
 
 We have $$W^{r,r}=\{[Id_r]\}$$
 $$W^{(r-1),r}=\{[Id_r], [Id_{r-2},r, r-1]\}.$$
 
 $$W^{(r-2),r}=$$
 $$\{[Id_{r-3},r-2, r,r-1],
  [Id_{r-3},r,r-2,r-1],[Id_{r-3},r-1,r-2,r],[Id_{r-3},r-1,r,r-2]\}.$$

 Let $w\in \Sigma_r$. Consider the element  $\hat r(w)$ of $\Sigma_{r-1}$
 such that  the list representing $\hat r(w)$ is the list $w$ where 
 the element $r$ is omitted.
In other words, if $u$ is such that 
 $w(u)=r$, then 
 $\hat r(w)(k)=w(k)$ if $k<u$, and 
 $\hat r(w)(k)=w(k+1)$ if $k\geq u$. We have the lemma.
 \begin{lemma}\label{ot}
 The element $w\in W^{i,r}$
 if and only if: 
\begin{itemize}
  \item $w^{-1}\{i,r\}$ is an interval.
  \item $\hat r(w)\in W^{i,r-1}$.
\end{itemize}
\end{lemma}

   We denote by $T(i,r+1)\in \Sigma_{r+1}$ the transposition of $i$ and $(r+1)$.
   It acts on $S_\Delta$, thus on $S_\Delta^*$.

\begin{proposition}\label{syme}
Let $1\leq i\leq r$.  Then the following holds: $$T(i,r+1)\cdot 
Ires_{x=0}= \sum_{w\in W^{i,r}} (-1)^{r+1-w^{-1}(i)}w \cdot 
Ires_{x=0}.$$ 
\end{proposition}
 \begin{proof} 
  If $g\in \Sigma_{r+1}$, we have 
 $g\cdot Ires_{x=0}=\sum_{w\in \Sigma_r}c_w^gIres_{x=0}^w$, as 
 elements $Ires_{x=0}^w=w\cdot Ires_{x=0}$ form a basis of linear forms on $S_\Delta$.
Let, for $w\in \Sigma_r$, $$f_w=\frac{1}{\prod_{1\leq p\leq 
r-1}(e^{w(p)}-e^{w(p+1)})(e^{w(r)}-e^{r+1})}.$$  

The coefficient 
 $c_w^g$ is equal to
 $\langle g\cdot Ires_{x=0}, f_w\rangle= Ires_{x=0}(g^{-1}f_w)$. 
  We need  to prove  
  \begin{itemize}
  \item If $w$ is not in $W^{i,r}$, then $$\langle Ires_{x=0}, T(i,r+1)f_w\rangle=0.$$
  
  \item  If $w\in W^{i,r}$, then $$\langle Ires_{x=0}, T(i,r+1) f_w\rangle=(-1)^{r+1-j}$$
  with $j=w^{-1}(i)$.
  
 \end{itemize}

    Let   $\phi=\frac{1}{\prod_{\alpha\in \sigma}\alpha}$ be in $S_{A_r}$.
    The set $\sigma$ consists on $r$ linearly independent elements of $E_r$.
   The partial residue $Res_{x_r=0}\phi$
   is non zero if and only one of the elements $\alpha$ in the denominator of 
   $\phi$ is proportional to  $(e^{r}-e^{r+1})$. Indeed, in computing 
   $Res_{x_r=0}\phi$, we replace $e^i$ by $x_i$, when $i\leq r$, and $e^{r+1}$ by 
   $0$. Thus the only factor creating a pole on $x_r=0$ is $(e^{r}-e^{r+1})$.
   It follows that 
    $Res_{x_r=0}(T(i,r+1)f_w)$ is zero unless  there exists  
  $p$, $1\leq p\leq r$, such that 
    $(e^{w(p)}-e^{w(p+1)})$ is  proportional to  the linear form $T(i,r+1)\cdot (e^r-e^{r+1}).$ 
       
    If $i=r$,  then $T(r,r+1)\cdot (e^r-e^{r+1})=-(e^r-e^{r+1})$ 
    and we see that necessarily  $p=r$ and $w(r)=r$. 
    Then $$Ires_{x=0}T(r,r+1)\cdot f_w$$
    $$=Res_{x_1=0}\ldots 
Res_{x_{r-1}=0}\frac{-1}{\prod_{k=1}^{r-2}(x_{w(k)}-x_{w(k+1)})x_{w(r-1)}} 
=-\delta^1_w$$ since $w(r)=r$ and that the basis $f_{w}$, with 
$w\in \Sigma_{r-1}$ is dual to the basis of iterated residues on 
$S_{A_{r-1}}$, as seen in Section \ref{chambers}. On the other 
hand   $\it W^{r,r}=[Id_{(r)}],$ and thus the case $i=r$ is 
completed.

Let $i<r$. We proceed by induction on $r$.   
   Let $w\in \Sigma_r$ and $j$ such that $w(j)=i$.
       Then, with $g=T(i,r+1)$,
$$g\cdot f_w=\frac{1}{(e^{w(1)}-e^{w(2)})(e^{w(2)}-e^{w(3)})\ldots 
(e^{w(j-2)}-e^{w(j-1})}\times$$$$\times 
 \frac{1}{(e^{w(j-1)}-e^{r+1})(e^{r+1}-e^{w(j+1)})\ldots (e^{w(j+1)}-e^{w(j+2)})
 \ldots (e^{w(r)}-e^{i})}$$
and 
 $$Ires_{x=0}g\cdot 
f_w=-Ires_{x=0}\frac{1}{(x_{w(1)}-x_{w(2)})(x_{w(2)}-x_{w(3)})\ldots 
(x_{w(j-2)}-x_{w(j-1})}$$$$\times 
 \frac{1}{x_{w(j-1)}x_{w(j+1)}(x_{w(j+1)}-x_{w(j+2)})\ldots (x_{w(r)}-x_{i})}.$$
The function $g\cdot f_w$ has only two simple poles  in 
$x_{w(j-1)}$ and in  $x_{w(j+1)},$ thus $Res_{x_r=0}g\cdot f_w\neq 
0$ iff 
 $w(j-1)=r$, or $w(j+1)=r.$
Precisely 
\begin{itemize}
\item if $w(j-1)=r$, then $Res_{x_r=0}g\cdot f_w$

$=\frac{-1}{(x_{w(1)}-x_{w(2)})...(x_{w(j-3)}-x_{w(j-2)}) 
x_{w(j-2)}\cdot 
 x_{w(j+1)}(x_{w(j+1)}-x_{w(j+2)})\ldots (x_{w(r)}-x_{i})}$

\item and if  $w(j+1)=r$, then $Res_{x_r=0}g\cdot f_w$

$=\frac{1}{(x_{w(1)}-x_{w(2)})\ldots 
(x_{w(j-2)}-x_{w(j-1)})x_{w(j-1)}\cdot 
x_{w(j+2)}(x_{w(j+2)}-x_{w(j+3)})\ldots (x_{w(r)}-x_{i})}$ 
\end{itemize}

In particular we see that if $w^{-1}\{i,r\}$ is not  an interval, 
then $Res_{x_r=0}g\cdot f_w$ is zero, so a fortiori 
$Ires_{x=0}g\cdot f_w$ is equal to $0$. 

If $w^{-1}\{i,r\}$ is  an interval, then we check on the 
preeceding formula that 

\begin{itemize}
  \item   if $w(j-1)=r$, then $$Ires_{x=0} T(i,r+1)\cdot 
f_w=Res_{x_1=0}\ldots Res_{x_{r-1}=0}T(i,r)\cdot f_{\hat r(w)}.$$ 

  \item if $w(j+1)=r$, then $$Ires_{x=0} T(i,r+1)\cdot 
f_w=-Res_{x_1=0}\ldots Res_{x_{r-1}=0}T(i,r)\cdot f_{\hat r(w)}.$$ 
\end{itemize}

We conclude by induction, using Lemma \ref{ot}. 
\end{proof}QED

\bigskip

We now consider the volume polynomial 
$v(A_r^+,\c^+)(a)=Ires_{x=0}( J_{A_r^+}(a)).$ 

As the function $J_{A_r^+}(a)$ is anti-invariant under the full 
group $\Sigma_{r+1}$, we obtain from PROPOSITION \ref{syme}  
\begin{proposition}\label{sym2}
Let $a=\sum_{i=1}^{r+1}a_ie^i$ with  $\sum_{k=1}^{r+1}a_{k}=0$.  
Then we have 

$$v(A_r^+,\c^+)(T(i,r+1)\cdot a)=\sum_{w\in \it 
W^{i,r}}\epsilon(w) (-1)^{r-w^{-1}(i)}v(A_r^+,\c^+)(w^{-1}\cdot 
a). $$ 

\end{proposition}

{\bf Example:} We know that the function $v(A_r^+,\c^+)(a)$ is a 
function of $(r-1)$ variables $v^+(a_1,a_2,....., a_{r-1})$. Then, 
for $i=r$, $r-1$, $r-2$, we get  from Proposition \ref{sym2}, and 
the description given of the  corresponding sets $W^{i,r}$, that 
the function $v^+(x)$ satisfies the identities:

\begin{itemize}
  \item $i=r$: $v^+(x)=v^+(x)$.

    \item $i=r-1$:
    for any values $x_1,..., x_{r}$, we have 
    
    $$v^+(x_1,....,x_{r-2},-( x_1+...+x_{r-2}+x_{r-1}+x_r))=$$
    $$-
    v^+(x_1,x_2,..., x_{r-2},x_{r-1})-v^+(x_1,...,x_{r-2},x_{r}).$$
    
  \item $i=r-2$: for any values $x_1,..., x_{r}$, we have 
$$v^+(x_1,x_2,..., x_{r-3},-(x_1+x_2+...+ x_{r}),x_{r-1})=$$ $$ 
-v^+(x_1,x_2,..., x_{r-3}, x_{r-2},x_{r}) +v^+(x_1,x_2,..., 
x_{r-3}, x_{r-1} , x_{r-2})$$ $$ -v^+(x_1,x_2,..., x_{r-3}, x_{r} 
, x_{r-2})+v^+(x_1,x_2,..., x_{r-3}, x_{r-1},x_{r}). $$ 
\end{itemize}

{\bf Remark}: A function $v(x_1,..., x_{r-1})$ of the form $ 
w(x_1,..., x_{r-2})(x_1+x_2+...+x_{r-2}+3 x_{r-1})$ satisfies (2). 
This is in agreement with the divisibility  of $v(A_r^+,\c^+)(a)$ 
by the linear factor $(a_1+a_2+...+a_{r-2}+3 a_{r-1})$ proved in 
Section \ref{volume}.

\bigskip

We now give a stronger result on  the  symmetry for the Kostant 
partition polynomial   $k(A_r^+,\c^+)$ attached to the nice 
chamber $\c^+$.

  Let $\rho=\frac{1}{2}\sum_{\alpha\in 
A_r^+}\alpha$. We have $$\rho=\frac{1}{2}(r e^1+(r-2)e^{2}+\cdots 
-(r-2)e^{r}-re^{r+1}).$$

As usual, we have a shifted symmetry property for $K_{A_r^+}$. 

\begin{lemma}

For any $w\in \Sigma_{r+1}$, the function 
 $$K_{A_r^+}({\bf a})= Tres_{\Delta} \left(\frac{e^{{\bf a}}}{\prod_{\alpha\in 
\Phi}(1-e^{-\alpha})}\right)$$ 
 satisfies the relation: 
  $$w\cdot K_{A_r^+}({\bf a})=\epsilon(w) K_{A_r^+}(w\cdot {\bf a}+w\cdot \rho-\rho).$$
  \end{lemma}

Thus PROPOSITION  \ref{sym} implies the following symmetry 
relation for $k(A_r^+,\c^+)(a)=Ires_{x=0} K_{A_r^+}(a).$ 

\begin{proposition}\label{sym3}

Let $a=\sum_{k=1}^{r+1} a_k e^k$ with $\sum_{k=1}^{r+1}a_{k}=0$. 
Then we have for, every $1\leq i\leq r$, 
$$k(A_r^+,\c^+)(T(i,r+1)\cdot a-(r-i+1)(e^i-e^{r+1}))$$ $$ 
=\sum_{w\in W^{i,r}}\epsilon(w) (-1)^{r-w^{-1}(i)}k(A_r^+,\c^+)(w 
^{-1} a+w^{-1}\cdot \rho-\rho).$$

\end{proposition}

{\bf Example:} We know that the function $k(A_r^+,\c^+)(a)$ is a 
function of $(r-1)$ variables $k^+(a_1,a_2,....., a_{r-1})$. Then, 
for $i=r$, $r-1$, $r-2$, we get  from Proposition \ref{sym3}, that 
the function $k^+(x)$ satisfies the identities:

\begin{itemize}
  \item $i=r$: $k^+(x)=k^+(x)$.

    \item $i=r-1$:
    for any values $x_1,..., x_{r}$, we have 
    
    $$k^+(x_1,....,x_{r-2},-( x_1+...+x_{r-2}+x_{r-1}+x_r+2))=$$
    $$-
    k^+(x_1,x_2,..., x_{r-2},x_{r-1})-k^+(x_1,...,x_{r-2},x_{r}-1).$$
    
  \item $i=r-2$: for any values $x_1,..., x_{r}$, we have 
$$k^+(x_1,x_2,..., x_{r-3},-(x_1+x_2+...+ x_{r}+3),x_{r-1})=$$ $$ 
-k^+(x_1,x_2,..., x_{r-3}, x_{r-2},x_{r}-1) +k^+(x_1,x_2,..., 
x_{r-3}, x_{r-1}-1 , x_{r-2}+1)$$ $$ -k^+(x_1,x_2,..., x_{r-3}, 
x_{r}-2 , x_{r-2}+1)+k^+(x_1,x_2,..., x_{r-3}, x_{r-1}-1,x_{r}-1). 
$$ 
\end{itemize}

{\bf Remark}: A function $k(x_1,..., x_{r-1})$ of the form $ 
w(x_1,..., x_{r-2})(x_1+x_2+...+x_{r-2}+3 x_{r-1}+3)$ satisfies 
(2). This is in agreement with the divisibility  of 
$k(A_r^+,\c^+)(a)$ by the linear factor $(a_1+a_2+...+a_{r-2}+3 
a_{r-1}+3)$ proved in Section \ref{some}.

\bigskip

{\bf Bibliographical remarks.} 

As explained as the beginning, the relations given here are 
transcription of the ''hidden'' action of $\Sigma_{r+1}$ on 
$S(A_r)$.

 \section{Appendix 1: Jeffrey-Kirwan residue formula for  the volume.}\label{appendix1}

We return to the notations  of Section 3. 
 Let
$\Phi=\{\alpha^1,...,\alpha^N\}$ be a sequence of elements of 
$\Delta^+$. Let $p:\R^N\to V^*$ be the map such that 
$$p(u_1w^1+u_2w^2+...+ u_Nw^N)= u_1 \alpha^1+u_2\alpha^2+\cdots 
+u_N\alpha^N .$$.

Consider  {\bf the family } of polytopes $P_\Phi(a)$ where  $a$ 
varies in $C(\Delta^+)\subset V^*$. Let $x\in V$, such that 
$\langle \alpha^k,x\rangle >0$ for all $\alpha^k\in \Phi$. We have 

$$\frac{1}{\prod_{k=1}^N\langle 
\alpha^k,x\rangle}=\int_{\R_+^N}e^{-\sum_{k=1}^Nu_k\langle 
\alpha^k,x\rangle}du_1 \cdots du_N=\int_{\R_+^N}e^{-\langle 
p(w),x\rangle} dw.$$ 

Let $s:V^*\to \R^N$  be a section from $V^*$ to $E$ such that 
$$\R^N=s(V^*)\oplus \Ker(p).$$ We write $w=s(a)+m$, with $p(w)=a$ 
and $p(m)=0$. Our choice of measure is such that $dw=da dm $. By 
Fubini, we obtain $$\int_{\R_+^N}e^{-\langle p(w),x\rangle} 
dw=\int_{C(\Delta^+)}e^{-\langle a,x\rangle} 
(\int_{P_\Phi(a)}dm)da.$$ 
 Thus 
$$\frac{1}{\prod_{k=1}^N\langle 
\alpha^k,x\rangle}=\int_{C(\Delta^+)} e^{-\langle x,a\rangle}\vol 
P_\Phi(a)da,$$ i.e. the Laplace transform of the function $\vol 
P_\Phi(a)da$ is $\prod_{k=1}^N\frac{1}{\alpha^k(x)}$. We need to 
inverse this formula to find  $\vol P_\Phi(a)$.

We denote by $G_{\Delta}$ the subspace of $R_\Delta$ spanned by 
the functions of the form $$\frac{1}{\prod_{\alpha\in 
\sigma}\alpha^{n_\alpha}}$$ where $\sigma$ is a basis of $\Delta$ 
and the $n_\alpha$ are  positive integers.

\begin{lemma}

Let $\kappa$ be a sequence of elements of $\Delta^+$ such that 
$\kappa$ generates $V^*$. Then the function 
$$\frac{1}{\prod_{\alpha\in \kappa }\alpha}$$ belongs to the 
vector space $G_\Delta$. 

\end{lemma}

\begin{proof}
We argue on the cardinal on the underlying set $\{\kappa\}$ to 
$\kappa$. If  the cardinal of $\{\kappa\}$ is  minimum, then 
$\{\kappa\}$ is a basis of $\Delta$. If not, there is  a 
linear relation  $\beta=\sum_{j} c_j\alpha^j$ between elements  
$\beta,\alpha^j$ belonging to  $\kappa$. Then 
$$\frac{1}{\beta\prod_j \alpha^j}=\sum_{j} 
\frac{c_j}{\beta^2\prod_{i\neq j}\alpha^i}.$$ We conclude by 
induction.  
 \end{proof} QED

Let $\sigma=\{\alpha^1,\alpha^2,...,\alpha^r\}$ be a basis of 
$\Delta$ consisting of   elements of $\Delta^+$. Let 
$$f_\sigma=\frac{1}{\alpha^1\cdots \alpha^r}$$ and let 
$[C(\sigma)]$ be the characteristic function of the cone 
$C(\sigma)$. 

Consider the function  $$F(x)=\frac{1}{\langle 
\alpha^1,x\rangle^{k_1+1}\cdots \langle 
\alpha^r,x\rangle^{k_r+1}}$$ where $k_j$ are non negative 
integers. From the one dimensional formula: 
$$\frac{1}{x^{k+1}}=\int_{\R^+} e^{-ax} \frac{a^{k}}{k!}da$$ we 
see   that $$F(x)=\int_{V^*}e^{-\langle 
a,x\rangle}v(a)[C(\sigma)](a) da$$  where $v$ is a polynomial.

Let $\c$ be a big chamber of $\Delta^+$. The  verification of the 
inversion formula: 

\begin{equation}\label{in} 
[C(\sigma)](a) v(a)=\langle\langle \c,Tres_\Delta (e^{ 
a}F)\rangle\rangle 
\end{equation}
 for $a\in \c$ is straightforward. 
 Indeed,  we write 
 $a=\sum_{k=1}^r a_k \alpha^k $, 
then $da=\vol(\sigma) da_1 da_2...da_r$ and 
$v(a)=\frac{1}{\vol(\sigma)} \frac{a_1^{k_1}}{k_1 !}\cdots  
\frac{a_r^{k_r}}{k_r !}$. On the other hand the function $F(x) 
e^{\langle a,x\rangle}$ is equal to  $$\frac{e^{a_1 \langle 
\alpha^1,x\rangle}}{\langle \alpha^1,x\rangle^{k_1+1}}\cdots 
\frac{e^{a_r \langle \alpha^r,x\rangle}}{\langle 
\alpha^r,x\rangle^{k_r+1}}$$ and its total residue is the function 
$$ \frac{a_1^{k_1}}{k_1 !}\cdots  \frac{a_r^{k_r}}{k_r !}\frac{1}{ 
\langle \alpha^1,x\rangle}\cdots \frac{1} {\langle 
\alpha^r,x\rangle} .$$ 

So $$Tres_\Delta(e^{a}F)=\frac{a_1^{k_1}}{k_1 !}\cdots  
\frac{a_r^{k_r}}{k_r !}f_\sigma.$$

Let $\c$ be a big chamber. If $\c\subset C(\sigma)$, then the left 
hand side of  the equation above is $v(a)$, and the right hand 
side is $$\langle\langle 
\c,f_\sigma\rangle\rangle\frac{a_1^{k_1}}{k_1 !}\cdots  
\frac{a_r^{k_r}}{k_r !}=\frac{1}{\vol(\sigma)}\frac{a_1^{k_1}}{k_1 
!}\cdots  \frac{a_r^{k_r}}{k_r !}=v(a).$$  
 If $\c$ is not contained in $C(\sigma)$, then 
 $\c\cap C(\sigma)=\emptyset$ and both sides are equal to $0$.

Thus the inversion formula for any function of $G_\Delta$ is 
established and  we obtain THEOREM \ref{jk}.

\newpage

\section{Appendix 2. Chambers and basis for $A_2$ and
$A_3$}\label{appendix2}

We briefly recall the setting. Following Section \ref{chambers} we 
realize $V^*$ as $E_r$ with basis $\{e^1,e^2,..e^r\}$, so that we 
may write:

 $$A_r=\{\pm e^1,\pm e^2,..,\pm e^r,\pm (e^i-e^j),
1\leq i<j\leq r\}.$$ Then  $\vert A_r^+ \vert =N=r(r+1)/2$ and the 
set of positive roots is $$A_r^+=\{e^1,e^2,..e^r,(e^i-e^j), 1\leq 
i<j\leq r\}.$$ 

 We consider the group $\Sigma_r$ of
permutations of $\{1,2,...,r\}$.  We  denote  a permutation  in 
$\Sigma_r$ as a list of $r$ elements . For example $[231]$ 
represents the permutation in  $3$ elements $\{1,2,3\}$ sending 
$1$ to $2$, $2$ to $3$ and $3$ to $1$. The group $\Sigma_r$ acts 
naturally on $V^*$.

The knowledge of the $S_\Delta$-valued polynomial function 
$$J_{A_r^+}(a)(x)=Tres_{A_r^+}[\frac{e^{\langle 
a,x\rangle}}{\prod_{\alpha\in A_r^+ 0}\langle\alpha,x\rangle}] 
=\frac{1}{(\frac{r(r-1)}{2})!} Tres_{A_r^+}[\frac{\langle 
a,x\rangle^{(\frac{r(r-1)}{2})}} {\prod_{\alpha\in A_r^+ 
}\langle\alpha,x\rangle}]$$ determines the various polynomial 
functions associated to big chambers.

The function $J_{A_r^+}$ is   anti-invariant under the  
group $\Sigma_{r}$: 
 $$J_{A_r^+}(w\cdot a)=\epsilon(w)
J_{A_r^+}(a)$$ thus  we need only to determine the polynomial 
$Ires_{x=0} J_{A_r^+}(a)$.  It is an homogeneous polynomial of 
degree $|A_{r-1}|$. (In fact, $J_{A_r^+}(a)$ is also  
anti-invariant under $\Sigma_{r+1}$ leading to some constraints on 
coefficients of $J_{A_r^+}(a)$.)

Recall that $C(\sigma)$ 
 denote the simplicial cone determined by  a basis $\sigma$ of  $A_r^+.$ 
 If we write $a=\sum_{i=1}^r a_ie^i$ for $a\in V^*$,   the cone generated by $A_r^+$ is $$C(A_r^+)=\{a\in V^* 
\text{such that} \ a_1\geq 0, a_1+a_2\geq 0,...,a_1+..+a_r\geq 
0\}$$

The small chambers, denoted by  $\a$,  are defined as the open 
connected component of $V^*-\CH^*$ in  $C(A_r^+)$. To any $h\in 
C(A_r^+)$ we associate the intersection of all the simplicial 
cones $C(\sigma)$ which contain $h.$ The interior of the maximal 
cones of this subdivision of $C(A_r^+)$ into polyhedral cones are 
called big chambers. The relevance of the big chambers lies in the 
fact that the polynomial volume is the same on all the small 
chambers that make up a big chamber. Recall that the nice chamber 
$\c^+$ is given by $a_i>0$. The nice chamber $\c^+$ is the cone 
$\sum_{i=1}^r a_i e^i$ with  $a_i>0,\ 1\leq i\leq r.$ 
 Using the permutation $w_0\in 
\Sigma_{r+1}$ reversing order on $\{1,2,...,(r+1)\}$, there is the 
``opposite" nice chamber $\c^-$, which is the cone  spanned by the 
roots $\{e^1, (e^1-e^2), (e^1-e^3),..(e^1-e^r)\}$.

If $w\in \Sigma_r$, we denote by $C_w^+\subset C({A_r}^+)$ the 
simplicial cone generated by the vectors   $$\epsilon(1) 
(e^{w(1)}-e^{w(2)}),...,\epsilon(r-1)  (e^{w(r-1)}-e^{w(r)}), 
(e^{w(r)}-e^{r+1}),$$ where $\epsilon(i)$ is 1 or -1 depending 
whether $w(i)<w(i+1)$ or not. 

 The space  $S_{A_r}$ is of dimension $r!$.
 As basis of $S_{A_r}$, we choose   elements  $f_w$ indexed by $w\in
\Sigma_r$, with 
 $$f_w=w\cdot \frac{1}{(e^1-e^2)(e^2-e^3)....(e^{r-1}-e^r)e^r}.$$

 To a big chamber $\c$ is associated a linear form $f\to \langle\langle\c,f\rangle\rangle$
 on $S_{A_r}$.
 By definition, $\langle\langle\c,f_w\rangle\rangle$ is equal to $0$ if $\c$ is not contained in
 $C_w^+$. Otherwise, 
$\langle\langle\c,f_w\rangle\rangle=(-1)^{n(w)}$, where $n(w)$ be 
the number of elements such that $w(i)>w(i+1)$. Thus to compute 
$\langle\langle\c,f\rangle\rangle=(-1)^{n(w)}Ires_{x=0}^wf$, we 
need to determine elements $w$ such that $\c\subset C_w^+$ and 
$n(w)$.

We compute the big chambers forms $\c$, and the corresponding form 
$\langle\langle\c,f\rangle\rangle$, in the $A_2$ and $A_3$ case.  
Only the form associated to $\c^+$ has a simple expression in the 
basis $f_w$. (Of course this basis is somewhat arbitrary).

 For $a\in \overline{\c},$ $\c$ a big chamber we recall the
transmutation formula for the volume and the Kostant partition 
function as given in Sections \ref{nice}: 
 $$\vol  P_{A_r^+}(a)= v(A_r^+,\c)=\sum _{|{\bf
i}|=\binom{r}{2}}f_\c({\bf i}) 
\frac{a_1^{i_1}}{i_1!}\frac{a^{i_2}}{i_2!}\cdots 
 \frac{a_{r-1}^{i_{r-1}}}{i_{r-1}!}
\frac{a_r^{i_{r}}}{i_{r}!}$$ 
  and
$$k(A_r^+,\c)(a)= \sum_{|{\bf i}|=\binom{r}{2}}  f_\c({\bf i}) 
\binom{a_1+r-1}{i_1} \binom{a_{2}+r-2}{i_{2}}\cdots 
\binom{a_{r-1}+1}{i_{r-1}} \binom{a_{r}}{i_{r}}.$$ where 
$\langle\langle\c, f_{A_r^+}({\bf i})\rangle\rangle=f_\c({\bf i})$

 For $A_2$ we may draw the  following picture:

\begin{figure}
 \centering
 \psfrag{f12}{$e_1-e_2$}
 \psfrag{f2}{$e_2$}
 \psfrag{f1}{$e_1$}
 \psfrag{c1}{$\c_2$}
 \psfrag{c2}{$\c_1$}
 \psfrag{f.12}{$a_1+a_2\geq 0$}
 \psfrag{f.1}{$a_1\geq 0$}
 \psfrag{f.2}{$a_2\geq 0$}
 \includegraphics{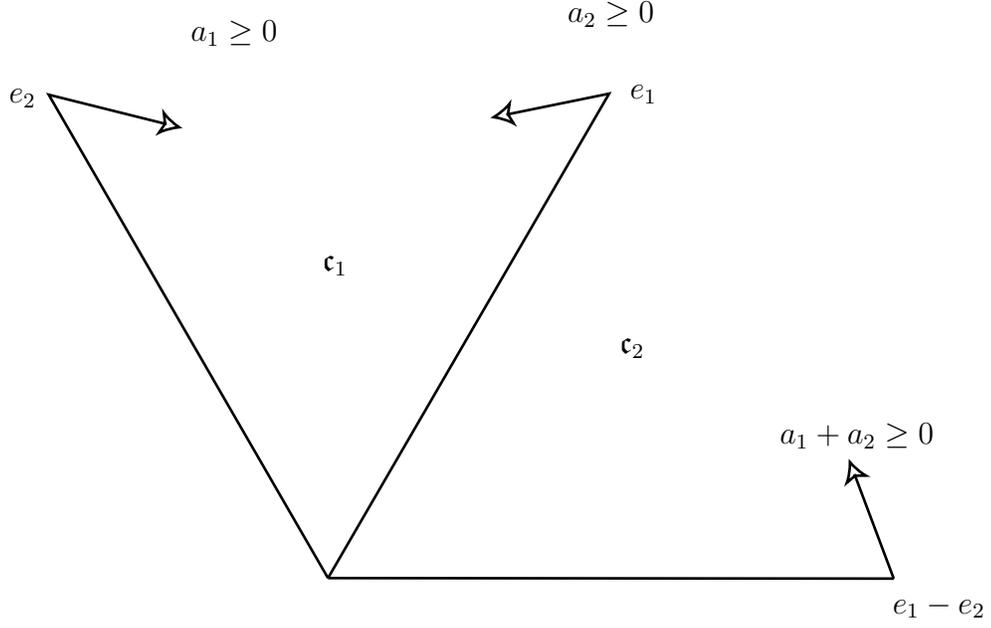}
 \caption{Chambers for $A_2$ with $a=a_1e^1+a_2e^2$.}
 \label{fig:3}
\end{figure}

 The equations defined by the
hyperplanes are: $$a_1=0,a_2=0,a_1+a_2=0.$$ In dimension 2, small 
and big chambers coincide. 

There are 3 bases in    $A_2^+$,$\sigma_1=\{e^1, e^2\}$, 
 $\sigma_2=\{e^1, (e^1-e^2)\}$, and
$\sigma_3=\{e^2, (e^1-e^2)\}$ bases of $\Delta$. The space 
$S_\Delta$ is spanned by $f_{\sigma_2}$,$f_{\sigma_3}$. We have 
the linear relation  $f_{\sigma_1}=f_{\sigma_3}-f_{\sigma_2}$.

 Therefore , if $a =a_1e^1+a_2e^2$,
 then
 $$J_{A_2^+}({\bf a})=a_1 f_{[12]}- a_2 f_{[21]},$$

$$v(A_2^+,\c^+)=a_1,$$ $$k(A_2^+,\c^+)=a_1+1.$$ 

$$\vol P_{A_r^+}(a)= a_1, \,\text{if}\, a\in \c_1$$ and $$\vol 
P_{A_r^+}(a)= a_1+a_2,\,\text{if}\, a\in \c_2.$$ 

As for the Ehrhart polynomials, we have $$k(A_r^+,\c_1)(a)=a_1+1\  \text {and}\ \
 k(A_r^+,\c_2)(a)=a_1+a_2+1.$$

When $r=3$ we may draw the following picture on a plane with 
equation $e_1= \text{constant}.$ 

\begin{figure}
 \centering
 \psfrag{a1}{$e_1$}
 \psfrag{a2}{$e_2$}
 \psfrag{a3}{$e_3$}
 \psfrag{a12}{$e_1-e_2$}
 \psfrag{a23}{$e_2-e_3$}
 \psfrag{a13}{$e_1-e_3$}
 \psfrag{c1}{$\a_1$}
 \psfrag{c2}{$\a_2$}
 \psfrag{c3}{$\a_3$}
 \psfrag{c4}{$\a_4$}
 \psfrag{c5}{$\a_5$}
 \psfrag{c6}{$\a_6$}
 \psfrag{c7}{$\a_7$}
 \psfrag{c8}{$\a_8$}
\psfrag{a.1}{$a_1\geq 0$} \psfrag{a.2}{$a_2\geq 0$} 
 \psfrag{a.3}{$a_3\geq 0$}
 \psfrag{a.12}{$a_1+a_2\geq 0$}
\psfrag{a.13}{$a_1+a_3\geq 0$} \psfrag{a.23}{$a_2+a_3\geq 0$} 
\psfrag{a.123}{$a_1+a_2+a_3\geq 0$} 

 \includegraphics{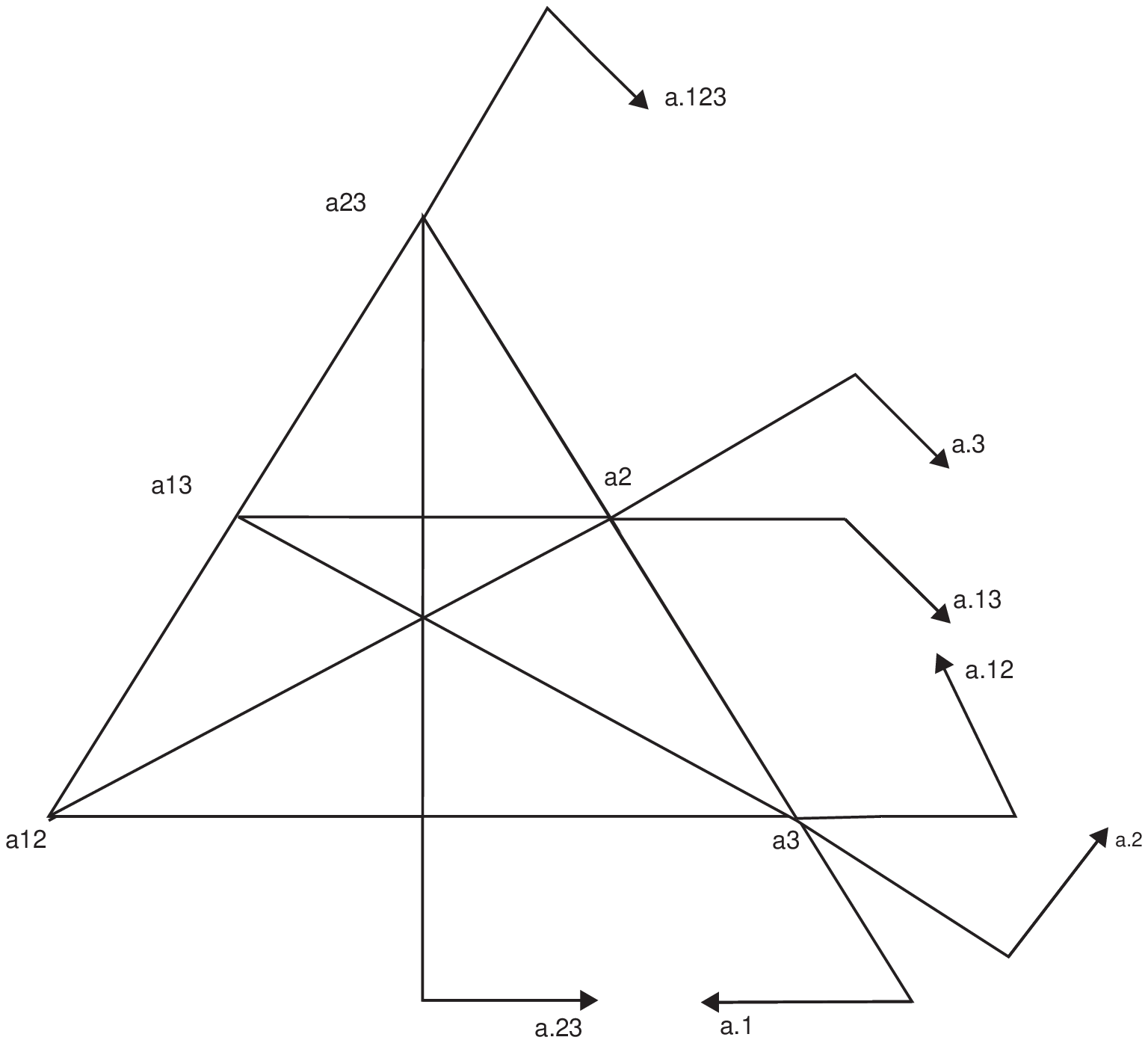}
 \caption{Chambers for $A_3$ with $a=a_1e^1+a_2e^2+a_3e^3$}
 \label{fig:4}

\end{figure}


 We see that there are $8$ small chambers $\a_i$ and $7$ big 
chambers $\c_k$.

 The big chambers are

\begin{tabular}{lll}
 
  $\c_1$&$=\a_1$&$=C(e_1,e_2,e_3).$\\
 $\c_2$&$=\a_2$&$= 
C(e_1,e_1-e_2,e_1-e_3).$\\
$\c_3$&$=\a_3$&$= C(e_1,e_2-e_3,e_2)\cap 
C(e_1,e_2,e_1-e_3).$\\
$ \c_4$&$=\a_4$&$= C(e_1-e_3,e_2-e_3,e_1)\cap C(e_1,e_2,e_1-e_3).$\\
$\c_5$&$=\a_5$&$=C(e_1,e_2-e_3,e_2)\cap C(e_2,e_1-e_3,e_2-e_3).$\\ 
$\c_6$&$=\a_6$&$=C(e_2,e_1-e_3,e_2-e_3)\cap C(e_1,e_2-e_3,e_1-e_3).$\\
$\c_7$&$=\a_7\cup \a_8 $&$=C(e_1-e_2,e_1,e_3).$
\end{tabular}

 For all $w\in \Sigma_3$, we have $n(w)=1$, except if 
$w=[123]$, where $n(w)=0$ and $w=[321]$ where $n(w)=2$. We have 

\begin{tabular}{ll}
$ C^+_{[123]}$&$=\c_1\cup \c_2\cup \c_3\cup \c_4\cup \c_5\cup 
\c_6\cup \c_7 $\\ 
 $C^+_{[213]}$&$=\c_2\cup \c_7 $\\
 $C^+_{[132]}$&$= \c_5\cup \c_6$\\
 $C^+_{[231]}$&$=\c_4\cup
\c_6$\\
 $C^+_{[312]}$&$=\c_2 \cup \c_3\cup \c_4$\\ 
  $ C^+_{[321]}$&$=\c_2\cup \c_4\cup \c_6.$
 \end{tabular}

Thus, we have rather complicated formula for the linear forms 
$\langle\langle\c,f\rangle\rangle$, except on the ``nice" chamber 
$\c^+=\c_1$.

\begin{tabular}{ll}
$ \langle\langle\c_1,f\rangle\rangle$&$=Ires_{x=0}^{[123]}f.$\\ 
$\langle\langle\c_2,f\rangle\rangle$&$=Ires^{[123]}_{x=0}f-Ires_{x=0}^{[213]} 
f-Ires_{x=0}^{[312]}f+Ires_{x=0}^{[321]}f.$\\ 
$\langle\langle\c_3,f\rangle\rangle$&$=Ires^{[123]}_{x=0}f-Ires_{x=0}^{[312]}f.$\\ 
$ \langle\langle\c_4,f\rangle\rangle$&$= 
Ires^{[123]}_{x=0}f-Ires_{x=0}^{[231]}f-Ires_{x=0}^{[312]}f+Ires_{x=0}^{[321]}f.$\\ 
$\langle\langle\c_5,f\rangle\rangle$&$=Ires^{[123]}_{x=0}f-Ires_{x=0}^{[132]}f.$\\ 
$\langle\langle\c_6,f\rangle\rangle$&$= 
Ires^{[123]}_{x=0}f-Ires_{x=0}^{[132]}f-Ires_{x=0}^{[231]}f+Ires_{x=0}^{[321]}f.$\\ 
$\langle\langle\c_7,f\rangle\rangle$&$=Ires^{[123]}_{x=0}f-Ires_{x=0}^{[213]}f.$
\end{tabular}

We easily compute the coefficients $f_{\bf i}$ of the function 
$$J_{A_3^+}(a_1 e^1+a_2 e^2+a_3 e^3)= 
\frac{1}{6}Tres_{A_3^+}(\frac{(a_1 x_1+a_2 x_2+a_3 x_3)^3}{x_1 x_2 
x_3 (x_1-x_2)(x_1-x_3)(x_2-x_3)}).$$

 Let $v_+(a)= \frac{1}{6}a_1^3+\frac{1}{2}a_1^2a_2.$

 We obtain:$$J_{A_r^+}(a_1 e^1+a_2 e^2+a_3 e^3)=
 \sum_{w\in \Sigma_r} \epsilon (w) v_+(w^{-1}\cdot a) f_w.$$

 Thus, we see that

 $$v(A_3^+,\c^+)(a)= \frac{1}{6}a_1^2(a_1+3a_2).$$

 $$k(A_3^+,\c^+)(a)= \frac{1}{6}(a_1+1)(a_1+2)(a_1+3a_2+3).$$

Thus here is the list of formulas for the volume and partition 
functions. We see that the polynomials $v(A_3^+,\c_i)$ are all 
differents so that the big chambers are indeed the minimal domains 
where the function $\vol$ is expressed by a polynomial formula.

 \begin{tabular}{ll}

$ v(A_3^+,\c_1)(a)$&$= \frac{1}{6}a_1^2(a_1+3a_2).$\\

$k(A_3^+,\c_1)(a)$&$= \frac{1}{6}(a_1+1)(a_1+2)(a_1+3a_2+3).$\\
 
 $ v(A_3^+,\c_2)(a)$&$=\frac{1}{6}(a_1+a_2+a_3)^2(a_1+a_2-2a_3).$\\
 
$k(A_3^+,\c_2)(a)$&$=\frac{1}{6}(a_1+a_2+a_3+1)(a_1+a_2+a_3+2)(a_1+a_2-2a_3+3).$\\

$v(A_3^+,\c_3)(a)$&$= 
\frac{1}{6}a_1^3+\frac{1}{2}a_1^2a_2-\frac{1}{2}a_1a_3^2-\frac{1}{6}a_3^3.$\\ 
$k(A_3^+,\c_3)(a)$&$=v(A_r^+,\c_3)(a) 
+a_1^2+\frac{3}{2}a_1a_2+\frac{1}{2}a_1a_3-\frac{1}{2}a_3^2$\\ 
&$+\frac{11}{6}a_1+a_2+\frac{2}{3}a_3+1.$\\ $ v(A_3^+,\c_4)(a)$&$= 
\frac{1}{6}a_1^3+ 
\frac{1}{2}a_1^2a_2-\frac{1}{2}a_1a_3^2-\frac{1}{6}a_2^3-\frac{1}{2}a_2^2a_3-\frac{1}{2}a_2a_3^2 
-\frac{1}{3}a_3^3.$\\ $ k(A_3^+,\c_4)(a)$&$= 
v(A_r^+,\c_4)(a)+a_1^2+\frac{3}{2}a_1a_2+\frac{1}{2}a_1a_3-\frac{1}{2}a_3^2$\\ 
&$+\frac{11}{6}a_1+\frac{7}{6}a_2+\frac{5}{6}a_3+1 .$ 
\\ $v(A_3^+,\c_5)(a)$&$=\frac{1}{6}a_1^2(2a_1+3a_2+3a_3). 
$\\ $ 
k(A_3^+,\c_5)(a)$&$=\frac{1}{6}(a_1+2)(a_1+1)(2a_1+3a_2+3a_3+3).$\\ 
$ v(A_3^+,\c_6)(a)$&$= \frac{1}{6}(a_1+a_2+a_3)^2(2a_1-a_2-a_3) 
.$\\ $ 
k(A_3^+,\c_6)(a)$&$=\frac{1}{6}(a_1+a_2+a_3+1)(a_1+a_2+a_3+2)(2a_1-a_2-a_3+3) 
$\\ $v(A_3^+,\c_7)(a)$&$=\frac{1}{6}(a_1+a_2)^3.$\\ 

$k(A_3^+,\c_7)(a)$&$=\frac{1}{6}(a_1+a_2+1)(a_1+a_2+2)(a_1+a_2+3).$ 
\end{tabular}

 Volume functions vanishes on the boundary of the cone 
$C(\Delta^+)$, thus all functions $v(A_3^+,\c_k)$ for $k\neq 3,4$ 
must have a linear factor. The order of this factor is computed  
as $L-1$, where $L$ is the number of roots  not on the face 
boundaring the chamber ( see for example jump formula in [B-V 2]), 
thus is $2$ for $\c_1,\c_2,\c_5,\c_6$ and $3$ for $\c_7$.

It is also reassuring to check that functions $k(A_3^+,\c_k)$  
define a continuous function on $C(\Delta^+)$ ( polynomials 
$k(A_3^+,\c_{k_1})$ and $k(A_3^+,\c_{k_2})$ agree on $\c_{k_1}\cap 
\c_{k_2}$. For example, if $a_2=a_3=0$, the  five  polynomials 
$k(A_3^+,\c_{k})(xe_1)$ $k=1,2,3,4,7$ restricts to 
$\frac{1}{6}(x+1)(x+2)(x+3).$

{\bf Bibliographical remarks.} Tables for $k(A_3^+,\c)$ are given 
for example in [S-B].

\newpage

{\center \bf REFERENCES}

\begin{itemize}
 
\item {[Ao]} K. Aomoto: Jacobi polynomials associated with Selberg integrals
.SIAM Math analysis (1987), p 545-549 

\item{ [B-S]}  A.Bincer and J.Schmidt: The Kostant partition function
for simple Lie algebras. J.Math.Phys. 25, (1984), p. 2367-2373.

\item {[B-V 1]}  M. Brion and M. Vergne: Residue formulae, vector partition
 functions and Lattice points in rational polytopes.
 J. A. M. S vol 10, (1997), pp 797-833.

 \item {[B-V 2]} M. Brion and M. Vergne: Arrangement of hyperplanes I~Rational 
functions and Jeffrey-Kirwan residue.  Ann. scient. {\'E}c. Norm. 
Sup., vol 32, (1999), p. 715-741.

\item {[C-R]} C.S. Chan and D.P. Robbins:
 On the volume of the polytope of doubly stochastic matrices.
  Electronic preprint (1998). math.CO/9806076. 
  Experimental Mathematics.

 \item  {[C-R-Y]} C.S. Chan, D.P. Robbins, D.S. Yuen:
 On the volume of  a certain  polytope.
  Electronic preprint (1998). math.CO/9810154.

\item  {[E]}. E.  Ehrhart. Sur un probl{\`e}me de g{\'e}om{\'e}trie diophantienne lin{\'e}aire 
II, J. Reine Angew. Math. 226, (1967);p . 1-29.

\item  {[J-K ]} L. C. Jeffrey and F. C. Kirwan: Localization
for non abelian group actions. Topology 34,(1995), p 291-327 . 

\item{[K-P]} G. Khovanskii, A. V. Pukhlikov:
A Riemann-Roch theorem for integrals and sums of quasipolynomials 
over virtual polytopes. St. Petersburg Math. J., vol 4 , (1993), 
p. 789-812.

\item{[Ki1]} A.N. Kirillov: Ubiquity of Kostka polynomials.
 Electronic preprint (1999).arXiv:math.QA/9912094.

\item{[Ki2]} A.N. Kirillov: Kostant partition functions for the root system of type A.
Preprint 1999.

\item{[Le]} G.I.  Lehrer:  On the Poincar{\'e} series associated 
with Coxeter group actions on complement of hyperplanes. J. London 
Math. Society  35(2) , (1987), pp. 275-294.
 
\item{[Li]} B.V. Lidskii: Kostant Function of the root system $A_n$.
Funktsional'nyi Analiz i Ego Prilozhenia, vol 18, (1983), p. 76-77

\item{[M]} W. Morris: Constant term identities for finite and 
affine root systems. (1982).PHD thesis, University of Wisconsin, 
Madison, Wisconsin.  

\item{[O-T]}  P. Orlik and H. Terao: Arrangements of Hyperplanes; Grundlehren
der Mathematischen Wissenscaften. vol 300, Springer-Verlag,  
Berlin 1992. 

\item{[Pi- S]} J. Pitman and R. Stanley: A polytope related to empirical
distribution, plane trees, parking functions, and the 
associahedron. Electronic preprint (1999). math.CO/9908029.

\item{[Po-St]} R.P. Stanley: Acyclic flow polytopes and Kostant's partition function.
Conference.  Transparencies (www-math.mit.edu/~rstan/trans.html).

\item{[Se]} A. Selberg: Bemerkninger om et multipelt integral, Norsk Mat. Tidsskr.,
vol 26, (1994) pp. 71-78

 \item{ [Sz]} A. Szenes:  Iterated residues and multiple Bernouilli
polynomials. hep-th/9707114. International Mathematical Research 
Notices.

\item{[Z]} D. Zeilberger: A conjecture of Chan, Robbins, and Yuen.
 Electronic preprint (1998). math.CO/9811108.

\end{itemize}

 {\small Welleda Baldoni Silva: Dipartimento di Matematica. Universit\'a
degli Studi di Roma Tor Vergata, Via della Ricerca Scientifica. 00133 Roma

E-mail: baldoni@mat.uniroma2.it} 

\bigskip

{\small  Mich{\`e}le Vergne: Centre de Math{\'e}matiques. Ecole 
Polytechnique. Palaiseau Cedex F-91128 France.}  

{E-mail: vergne@math.polytechnique.fr} 

{web-page : http://math.polytechnique.fr} 

\end{document}